\documentclass{ws-m3as}
\usepackage{placeins}
\usepackage{cite}
\usepackage[cal=boondox]{mathalfa}
\usepackage{url}
\usepackage{lineno}
\usepackage[english]{babel}
\usepackage{amssymb}
\usepackage{algpseudocode} 
\usepackage{algorithm} 
\usepackage{graphicx}
\usepackage{caption}
\usepackage{bm}
\usepackage{bbm} 
\usepackage{mathabx}
\usepackage{stmaryrd}
\usepackage{color}
\usepackage{subcaption} 
\usepackage{mathtools}
\usepackage{braket}
\usepackage{nicefrac}
\usepackage{cancel}
\usepackage{booktabs}

\usepackage{csvsimple}
\usepackage{siunitx}
\PassOptionsToPackage{dvipsnames}{xcolor}
\usepackage[customcolors,norndcorners]{hf-tikz}
\usepackage{tikz}
\usepackage{pgfplots}
\pgfplotsset{compat=newest}
\usepgfplotslibrary{patchplots}
\usetikzlibrary{pgfplots.patchplots}
\tikzset{set fill color=yellow, set border color=yellow}
\renewcommand{\theequation}{\arabic{section}.\arabic{equation}}
\numberwithin{equation}{section}
\modulolinenumbers[5]

\DeclareMathOperator{\supp}{supp}
\DeclareMathOperator{\id}{I}
\newtheorem{prop}{Proposition}[section]
\newtheorem{assumption}{Assumption}[section]
\DeclarePairedDelimiter{\abs}{\lvert}{\rvert}		
\DeclarePairedDelimiter{\norm}{\lVert}{\rVert}		
\DeclarePairedDelimiter{\jump}{\llbracket}{\rrbracket}
\DeclarePairedDelimiter{\opnorm}{\interleave}{\interleave}
\DeclareGraphicsExtensions{{.png},{.pdf},{.jpg}}
\renewcommand{\H}{\widehat}
\newcommand{\Sp}{\mathcal{S}}
\newcommand{\nZ}{n_{\mathrm{Z}}}
\newcommand{\Q}{\mathcal{Q}}
\newcommand{\io}{\mathcal{I}_{\Omega}}
\newcommand{\is}{\mathcal{I}_{\Sigma}}
\newcommand{\gl}{\Gamma_\ell}
\newcommand{\iv}{\mathcal{I}_{\mathcal{V}}}
\newcommand{\dn}{\partial_n}
\newcommand{\hh}{\mathcal{H}}
\newcommand{\ic}{\mathcal{I}_{\mathcal{C}}}
\newcommand{\M}{\mathcal{M}}
\newcommand{\F}{\textbf{F}}
\newcommand{\dnl}{\partial_{n}}
\newcommand{\dx}{\mathrm{d}\mathbf{x}}
\newcommand{\ds}{\mathrm{d}\sigma}
\newcommand{\LambdX}{\Lambda}
\newcommand{\LambdY}{\widetilde{\Lambda}}
\newcommand{\ml}{{ \mathcal{p}(\ell)}}
\newcommand{\sell}{ \mathcal{s}(\ell)}
\newcommand{\citeRefs}[1]{{ Refs.~\cite{#1}}}
\newcommand{\citeOneRef}[1]{{ Ref.~\cite{#1}}}
\newcommand{\citeRefDetail}[2]{{ Ref.~\cite{#2}, #1}}

\begin{document}
	
	\markboth{A. Benvenuti, G. Loli, G. Sangalli, T. Takacs}{Isogeometric  $C^1$ mortar method}
	
	%
	%
	
	\title{Isogeometric  $C^1$ mortar method}
	
	\author{A. Benvenuti}
	
	\address{Dipartimento di Matematica ``F. Casorati'',Universit\`a di Pavia, Via A. Ferrata 1\\
			Pavia, 27100, Italy\\
			bnvndr@gmail.com}
	
	\author{G. Loli}
	
	\address{Dipartimento di Matematica ``F. Casorati'',Universit\`a di Pavia, Via A. Ferrata 1\\
		Pavia, 27100, Italy\\
		gabriele.loli@unipv.it}
	
	\author{G. Sangalli}
	
	\address{Dipartimento di Matematica ``F. Casorati'',Universit\`a di Pavia, Via A. Ferrata 1\\
		Pavia, 27100, Italy\\
		giancarlo.sangalli@unipv.it}
	
	\author{T. Takacs}
	
	\address{Johann Radon Institute for Computational and Applied Mathematics, Austrian Academy of Sciences, Altenberger Str. 69\\
		Linz, 4040, Austria\\
		thomas.takacs@ricam.oeaw.ac.at}
	
	\maketitle
	
	
	\begin{abstract}
		We present an isogeometric mortar method for the discretization of the biharmonic equation posed on multi-patch domains. We  assume  only $C^0$-conformity at interfaces and { employ} a mortar approach to weakly enforce $C^1$-continuity across patch interfaces. Discrete inf-sup stability is ensured by selecting a Lagrange multiplier space consisting of splines of degree reduced by two compared to the primal space, with increased smoothness or merged elements near vertices.
		We prove optimal a priori error estimates and confirm the theoretical findings with a series of numerical experiments.
	\end{abstract}
	
	
	
	\section{Introduction}
Isogeometric Analysis (IgA) is a computational framework for the discretization of Partial Differential Equations (PDEs), introduced in~\citeOneRef{cottrell}, that aims to unify Finite Element Analysis (FEA) and Computer-Aided Design (CAD). In IgA, the same functions, such as B-splines and Non-Uniform Rational B-Splines (NURBS), are employed for both the geometric representation of the computational domain and the numerical approximation of PDE solutions. The use of smooth splines for PDE discretization offers several advantages: high-continuity splines can deliver greater accuracy per degree of freedom compared to $C^0$ piecewise polynomials, { see~\citeRefs{Evans_Bazilevs_Babuska_Hughes,Sangalli2018,bressan2018approximation}}, and they enable the direct approximation of higher-order PDEs, such as the Kirchhoff--Love plate and shell models, { cf.~\citeRefs{benson,beirao2,Niiranen2017,Ning2018,kiendl2009isogeometric,kiendl2010bending}}, the Cahn--Hilliard equation, { cf.~\citeOneRef{Gomez2008}}, the Navier--Stokes--Korteweg equation, { cf.~\citeOneRef{gomez2}}, and the approximation of the Optimal Transport solution, { cf.~\citeOneRef{Matthes04112009}}.

In practical applications, complex geometries are typically described by multiple patches rather than a single global parametrization. While high-order continuity is naturally achieved within each patch, enforcing strong coupling across patch interfaces in a multi-patch setting requires careful treatment, { see~\citeRefs{collin,kapl2018construction,KAPL201955,HUGHES2021467,weinmuller2022approximate,FARAHAT2023115706,Farahat2023,grovselj2020super,GROSELJ2024128278}. Recently, approximation error bounds for strongly coupled discretizations over so called analysis-suitable $G^1$ multi-patch parameterizations from~\citeOneRef{collin} were developed in~\citeOneRef{hasanova2026robust}}. Alternatives to strong coupling include non-conforming discretizations or the use of Nitsche’s method to address discontinuities of normal derivatives across interfaces, { as in~\citeRefs{apostolatos2014,weinmuller2021construction,weinmuller2022approximate}}, though these approaches require appropriate regularization parameters. { Alternatively, one can employ discretizations where normal discontinuities are minimized by construction, as in~\citeRefs{seiler2018approximately,Takacs2023}.} A widely adopted technique for coupling subdomains is the mortar method, { cf.~\citeRefs{bernardidec,bernardimor,bernardibas}}, which weakly enforces continuity constraints across interfaces via Lagrange multipliers. The choice of the multiplier space is crucial, as it must provide sufficient approximation properties and satisfy the inf-sup stability condition. In the isogeometric context, mortar methods have been extensively investigated; see, e.g.,~\citeRefs{HESCH2012104,Dornisch2015,dittmann2019weak,APOSTOLATOS2019368,DITTMANN2020112768,ADAM2020113403,Chasapi2020,coradello2021projected,CORADELLO2021114187,Merkel9807359,TONG2022115263,Hesch2022,brivadis,dornisch2021isogeometric}.
{ Related concepts are hybrid formulations, as in~\citeOneRef{horger2019hybrid}, and the IETI method for patch coupling, as in~\citeOneRef{kapl2026isogeometric}.}

In particular, the weak imposition of $C^0$-continuity along multi-patch interfaces within the isogeometric mortar framework is studied in~\citeOneRef{brivadis}. For a primal isogeometric space of degree $p$, several choices for the multiplier 
{ space, of degrees $p$, $p-1$, and $p-2$, are}  considered. The $p/p-1$ pairing is found to be numerically unstable, while the $p/p$ pairing, though stable, requires special boundary modifications (such as degree reduction near vertices) to avoid instabilities. In contrast, the $p/p-2$ pairing is theoretically proven to be stable. Although this pairing is expected to yield suboptimal convergence rates, numerical experiments for certain parametrizations indicate improved convergence behavior beyond the theoretical predictions.

Alternative choices for the multiplier space, analytically proven to satisfy the inf-sup condition, are explored in~\citeOneRef{dornisch2021isogeometric}. Additionally, an isogeometric mortar coupling method designed to ensure tangential continuity in the discretization of the Maxwell eigenvalue problem is analyzed in~\citeOneRef{buffa2020}, where theoretical stability is established for multiplier spaces of degree $p-1$ and numerical stability is demonstrated for degrees $p-k$ with odd $k$.

While previous works have primarily addressed the enforcement of weak $C^0$-continuity, this paper focuses on the imposition of weak $C^1$-continuity constraints across patch interfaces in $C^0$-conforming multi-patch domains. As a model problem, we consider the two-dimensional biharmonic equation. 
We adopt the $p/p-2$ pairing for the primal and Lagrange multiplier spaces, following~\citeOneRef{brivadis}. { While the theory in~\citeOneRef{brivadis} yields suboptimal convergence rates for second order PDEs for this pairing, it is a suitable choice for fourth order PDEs, where the expected optimal rate in the energy norm, i.e., $H^2$, with respect to the mesh size $h$ is $O(h^{p-1})$, which balances correctly with the required approximation properties of the Lagrange multiplier space.} To ensure the correct dimension of the multiplier space and guarantee stability, we introduce a smoothing procedure near vertices: for splines of maximum smoothness, we coarsen the mesh by merging the first two and last two elements; otherwise, we increase the smoothness at the first and last knots by one.

The structure of the paper is as follows. Section~\ref{sec:C1_mortar_formulation} introduces the model problem and its variational mortar formulation. Section~\ref{sec:preliminaries} recalls basic notions of B-splines and isogeometric spaces, particularly those used for discretizing the model problem. Sections~\ref{sec:Xhprop} and~\ref{sec:Vhprop} present approximation results for the global unconstrained space and the weakly $C^1$ space (the mortar kernel), respectively. Section~\ref{sec:error_analysis} provides an a priori error estimate, based on the assumption that a discrete inf-sup condition holds. This central { question} is addressed in Section~\ref{sec:dscinfsup}, where we introduce a suitable multiplier space that satisfies the discrete inf-sup condition. Section~\ref{sec:numerics} presents numerical results that assess the effectiveness of the proposed method and demonstrate its optimal approximation properties. Finally, Section~\ref{sec:conclusions} draws the main conclusions { and provides an outlook towards possible extensions of the theory}.

	\section{Model problem  and mortar formulation}\label{sec:C1_mortar_formulation}

{ In the following we introduce the model problem that we consider throughout the paper and present the multi-patch setting over which we will then define the mortar formulation of the model problem. 

\subsection{Model problem}}

{
Let $\Omega \subset \mathbb{R}^2$ be a bounded connected domain with a piecewise polynomial boundary. We consider the following fourth-order elliptic problem with homogeneous boundary conditions:
\begin{equation}\label{eq:bilaplace_strong}
\begin{cases}
	\begin{aligned}
		\Delta^2 u &= f & &\text{in } \Omega, \\
		u &= 0 & &\text{on } \partial\Omega, \\
		\dn u &= 0 & &\text{on } \partial\Omega,
	\end{aligned}
\end{cases}
\end{equation}
where $f \in H^{-2}(\Omega)$. { Here $\dn$ denotes the normal derivative.}

\begin{remark}\label{rem:homogeneous_reduction}
This homogeneous framework can directly account for non-homogeneous boundary data via a standard lifting technique. Let $\hat{u}$ be the solution to the non-homogeneous problem:
\begin{equation*}
\begin{cases}
	\begin{aligned}
		\Delta^2 \hat{u} &= \hat{f} & &\text{in } \Omega, \\
		\hat{u} &= g_0 & &\text{on } \partial\Omega, \\
		\dn \hat{u} &= g_1 & &\text{on } \partial\Omega,
	\end{aligned}
\end{cases}
\end{equation*}
where $\hat{f} \in H^{-2}(\Omega)$, and $g_0, g_1$ are given boundary data. Assuming there exists a lifting function $\bar{u} \in H^2(\Omega)$ that satisfies the boundary conditions $\bar{u} = g_0$ and $\dn \bar{u} = g_1$ on $\partial\Omega$, we can define the homogeneous variable as $u := \hat{u} - \bar{u}$. 

By linearity, $u$ satisfies the reference problem~\eqref{eq:bilaplace_strong} with the modified source term $f = \hat{f} - \Delta^2 \bar{u} \in H^{-2}(\Omega)$.
\end{remark}
}
A possible variational formulation of~\eqref{eq:bilaplace_strong} is: \textit{Find $u \in H^2_0(\Omega)$ such that}
\begin{equation}\label{eq:bilaplace_var}
	\int_{\Omega} \nabla (\nabla u ) : \nabla (\nabla v ) \, \dx = \langle f, v \rangle_{H^{-1} \times H^1_0}, \qquad \forall v \in H^2_0(\Omega),
\end{equation}
{ where $\nabla (\nabla(\cdot))$ denotes the Hessian and} 
$H^2_0(\Omega)$ denotes the Sobolev space of functions in $L^2(\Omega)$ with second derivatives in $L^2( \Omega)$ and null traces and normal derivatives on the boundary $\partial \Omega$.
In~\eqref{eq:bilaplace_var}, we assume the datum  $f$ is more regular than what is necessary, that is we assume  $f \in H^{-1}(\Omega) = (H^1_0(\Omega))'$. The bilinear form in \eqref{eq:bilaplace_var} results from integration by parts and from the equivalence {
\[
\int_{\Omega} (\partial_{x}\partial_{x} w) \, (\partial_{y}\partial_{y} v) \, \dx = \int_{\Omega} (\partial_{x}\partial_{y} w)  \, (\partial_{x}\partial_{y} v) \, \dx
\]
which holds for any $w,v \in H^2_0(\Omega)$. Here $\partial_{x}$ and $\partial_{y}$ denote the derivatives in directions $x$ and $y$, respectively, where $\mathbf{x} = (x,y)$.}

{ \subsection{Multi-patch segmentation}}

The mortar method applies to a decomposition of the domain $\Omega$ into non-overlapping subdomains, which are commonly referred to as patches in the context of isogeometric analysis. Thus, we consider a family of open sets $\left\lbrace \Omega_k\right\rbrace _{k \in \io}$, edges (open line segments) $\left\lbrace \Gamma_\ell\right\rbrace _{\ell \in \mathcal{I}_\Gamma}$, and vertices $\left\lbrace \mathcal{V}_i\right\rbrace _{i \in \iv}$, such that
\begin{equation*}
	\overline{\Omega} = \bigcup_{k \in \io} {\Omega}_k \cup \bigcup_{\ell \in \mathcal{I}_\Gamma} {\Gamma}_\ell \cup \bigcup_{i \in \iv} {\mathcal{V}}_i, 
\end{equation*}
where $\io$, $\mathcal{I}_\Gamma$, and $\iv$ denote finite sets of patch, edge, and vertex indices, respectively.
\begin{assumption}\label{ass:mp_conf}
We assume that 
\begin{equation*}
	\Omega_k \cap \Omega_{k'} = \emptyset, \quad \forall \, k,k' \in \io, \, k \neq k',
\end{equation*}
and the patch segmentation has no hanging nodes; thus, for $k\neq k'$, the set $\overline{\Omega}_k \cap \overline{\Omega}_{k'}$ is either empty, a vertex, or exactly one edge and two vertices.
\end{assumption}
Let $\is \subset \mathcal{I}_\Gamma$ denote the index set of all interior edges, i.e., $\ell \in \is$ if and only if ${\Gamma}_\ell \subset \Omega$. The skeleton $\Sigma$, representing the union of all interior edges, is defined as
\begin{equation*}
	\Sigma := \bigcup_{\ell \in \is} \Gamma_\ell,
\end{equation*}
where
\begin{equation*}
	\overline\Sigma = \overline{\bigcup_{k \in \io} \partial \Omega_k \setminus \partial \Omega}.
\end{equation*}
For each interface $\Gamma_\ell$, with $\ell\in\is$, there exists a pair of indices $k, k' \in \io$ such that
\begin{equation*}
	\overline{\Gamma}_{\ell} = \partial \Omega_k \cap \partial \Omega_{k'}.
\end{equation*}
Moreover, we assign arbitrarily an ordering to the patch indices with $\{k,k'\}=\{\ml,\sell\}$, referring to $\Omega_{\ml}$ as the \textit{primary patch} of $\Gamma_\ell$ and to $\Omega_{\sell}$ as the \textit{secondary patch}. { Analogous to the normal at the domain boundary, we denote by $n$ the outward unit normal vector to $\Gamma_\ell$ from the primary side, and by $\dnl (\cdot) := n \cdot \nabla (\cdot)$ the outer normal derivative on $\Gamma_\ell$ along $n$.}

{ \subsection{Mortar formulation}

We define the continuous mortar formulation for a space that is globally only $H^1$, patch-wise $H^2$ and where all normal jumps vanish at all vertices. Thus, the space is formally defined as}
\begin{equation*}\label{eq:spaceX}
	X:=\Set{v \in H^1(\Omega) \left| \ \begin{aligned}
		&v \lvert_{\Omega_k} \in H^2(\Omega_k), \forall k \in \io, \\
		&v = \dn v = 0, \text{ on }\partial \Omega, \\
		&\jump{\dnl v}_{\gl} \in H^{\nicefrac{1}{2}}_{00}(\gl), \forall \ell \in \is
		\end{aligned} \right. },
\end{equation*}
where $\jump{\cdot}_{\gl}$ denotes the jump, defined as the trace from the secondary side $\sell$ minus the trace from the primary side $\ml$. The space $X$ is equipped with the norm
\begin{equation}\label{eq:normX}
	\norm{v}_X^2 := \sum_{k \in \io} \norm{v}_{H^2(\Omega_k)}^2 + \sum_{\ell \in \is} \norm{\jump{\dnl v}_{\gl}}_{H^{\nicefrac{1}{2}}_{00}(\gl)}^2.
\end{equation}
The continuous Lagrange multiplier space is defined as
\begin{equation*}
	M := \prod_{\ell \in \is} H^{\nicefrac{-1}{2}}(\Gamma_\ell).
\end{equation*}
We introduce the continuous bilinear form on $X \times X$
\begin{equation*}
	a(u,v) := \sum_{k \in \io} \int_{\Omega_k} \nabla (\nabla u ) : \nabla (\nabla v ) \, \dx,
\end{equation*}
and the continuous bilinear form on $X \times M$
\begin{equation*}
	b(v,\mu) := \sum_{\ell \in \is} \int_{\Gamma_\ell} \mu \jump{\dnl v}_{\gl} \, \ds,
\end{equation*}
{ which induces} the linear continuous operator $B:X \rightarrow M'$, defined as
\begin{align*}
    \left\langle Bv,\mu \right\rangle_{M' \times M}= b(v,\mu), \quad \forall v \in X,\quad \forall \mu \in M.
\end{align*}
The following characterization holds:
\begin{equation}\label{eq:char_ker}
	H^2_0(\Omega) = \Set{v \in X : b(v,\mu) = 0, \forall \mu \in M}=\ker B.
\end{equation}
\begin{lemma}\label{lemma:normeq2}
The seminorm $\abs{v}_{H^2(\Omega)} = \norm{\nabla (\nabla v )}_{L^2(\Omega)}$ is a norm on $H^2_0(\Omega)$. In particular, there exists a constant $C > 0$ such that 
for all $v \in H^2_0(\Omega)$ { we have}
\begin{equation*}
	\norm{v}_{H^2(\Omega)} \leq C \abs{v}_{H^2(\Omega)}.
\end{equation*}
\end{lemma}
\begin{proof}
For the proof  we refer to \citeRefDetail{Chapitre III, Lemme 11.1}{Ladyzhenskaya1968}.
\end{proof}

We can now introduce the mortar variational formulation, that is  the saddle point problem: 
\textit{Find $(u,\lambda)\in X \times M$ such that}
\begin{equation}\label{eq:mortar-variational-formulation}
\begin{cases}
	\begin{aligned}
		a(u, v) + b(v, \lambda) & = \langle f, v \rangle_{H^{-1} \times H^1_0} & & \forall v \in X \\
		b(u, \mu)          & = 0    & & \forall \mu \in M.
	\end{aligned}
\end{cases}
\end{equation}
The well posedness of \eqref{eq:mortar-variational-formulation} is guaranteed by \citeRefDetail{Theorem 4.2.1}{boffi} and Lemma~\ref{lemma:normeq2}. Furthermore, if  $(u, \lambda) \in X \times M$ is the solution of problem~\eqref{eq:mortar-variational-formulation}, then $u \in H^2_0(\Omega)$ thanks to \eqref{eq:char_ker} and then is the solution of problem~\eqref{eq:bilaplace_var}.
	\section{Preliminaries}
\label{sec:preliminaries}
In this paper, we denote by $C$ a positive constant which may vary with each occurrence but remains independent of the mesh size $h$. 
\subsection{B-splines}\label{sec:bsplines}
In the following, we provide an overview of B-splines, introducing the necessary notation and fundamental concepts. For a comprehensive exposition, we refer the reader to~\citeRefs{cottrell,bazilevs,schumaker}.\\
Given a positive integer $p$, a $p$-open knot vector $\Xi$ is defined as a non-decreasing sequence of real numbers:
\begin{equation*}
	0 = \xi_{1} = \ldots = \xi_{p+1} < \xi_{p+2} \le \ldots \le \xi_n < \xi_{n+1} = \ldots = \xi_{n+p+1} = 1,
\end{equation*}
The breakpoint vector associated to $\Xi$, defined as
\begin{equation*}
	Z = (\zeta_1,\ldots,\zeta_{\nZ}), \quad \text{with } 0=\zeta_1<\cdots<\zeta_{\nZ} = 1,
\end{equation*}
collects the knots without repetition. For any breakpoint $\zeta_j \in Z$, we denote by $m_j$ its multiplicity in $\Xi$.
The B-spline basis functions of degree $p$ on $\Xi$ are defined recursively using the \emph{Cox--de Boor formula}. For $1 \leq i \leq n$, the basis functions are given by:
\begin{equation*}
	\begin{split}
		B_{i,0}({ \H{x}}) &=
		\begin{cases}
			1 & \text{if } \xi_i \leq { \H{x}} < \xi_{i+1}, \\
			0 & \text{otherwise,}
		\end{cases}\\
		B_{i,p}({ \H{x}}) &= \frac{{ \H{x}} - \xi_i}{\xi_{i+p} - \xi_i} B_{i,p-1}({ \H{x}}) + \frac{\xi_{i+p+1} - { \H{x}}}{\xi_{i+p+1} - \xi_{i+1}} B_{i+1,p-1}({ \H{x}}), \quad p \geq 1,
	\end{split}
\end{equation*}
where any division by zero is defined to yield a coefficient of zero. These functions span the space of piecewise polynomials of degree $p$, with $p - m_j$ continuous derivatives at the breakpoints $\zeta_j$:
\begin{equation*}
	\Sp^p(\Xi) := \operatorname{span}\{ B_{i,p} \}_{i=1}^n.
\end{equation*}
\begin{assumption}\label{ass:reg}  
	We assume that the spline space $\Sp^p(\Xi)$ is globally at least $C^1$-continuous. Specifically, we require $p \geq 2$ and $m_j \leq p - 1$ for $j = 2, \dots, \nZ - 1$.  
\end{assumption}  
{ Given a bi-degree ${\mathbf p} = (p, p)$ and two open knot vectors $\Xi_1$ and $\Xi_2$, the tensor-product B-spline basis is defined as
\begin{equation*}  
	{\cal B}_{\mathbf p}(\boldsymbol{\Xi}) := \{ \H{B}_{{\mathbf i}, {\mathbf p}}(\H{\mathbf{x}}) = \H{B}_{i_1, p}(\H{x}_1) \H{B}_{i_2, p}(\H{x}_2), \quad \H{\mathbf{x}} \in (0,1)^2 \},  
\end{equation*}  
where the corresponding spline space they span is given by  
\begin{equation*}  
	\Sp^{\mathbf p} = \Sp^{\mathbf p}(\Xi_1, \Xi_2) := \Sp^{p}(\Xi_1) \otimes \Sp^{p}(\Xi_2) = \operatorname{span}({\cal B}_{\mathbf p}(\boldsymbol{\Xi})).  
\end{equation*}  
While it is possible to define spline spaces of different degree in different directions, this straight-forward extension would further complicate the presentation of our results. Thus we restrict ourselves to the case of equal degrees in all directions and at all patches.
}
The breakpoints $(Z_1, Z_2)$ associated with the knot vectors $(\Xi_1, \Xi_2)$ define a mesh $\Q_h$ on $\H{\Omega}$, which partitions the parametric domain into two-dimensional elements:  
\begin{equation}\label{eq:mesh}    
		\Q_h = \Q_h(Z_1, Z_2) := \{ Q = (\zeta_i^1, \zeta_{i+1}^1) \times (\zeta_j^2, \zeta_{j+1}^2) : 1 \leq i < {\nZ}_1,  1 \leq j < {\nZ}_2 \}.   
\end{equation}  
The global mesh size $h$ is defined as  
\begin{equation*}  
	h := \max_{\substack{1 \leq d \leq 2 \\ 1 \leq i < {\nZ}_d }} h^d_{i},  
\end{equation*}  
where $h^d_i := \zeta^d_{i+1} - \zeta^d_i$ and $\zeta^d_i$ represents the $i$-th breakpoint in $Z_d$.  
\begin{assumption}\label{ass:quasiunif}  
	We assume that the knot vectors are quasi-uniform. Specifically, there exists a constant $\alpha > 0$, independent of $h$, such that  
	\begin{equation*}  
		\alpha h \leq h^d_i, \quad \text{for } 1 \leq i < {\nZ}_d \text{ and } d = 1,2.  
	\end{equation*}  
\end{assumption}  
	\subsection{Multi-patch isogeometric discretization}\label{sec:iga_disc}
We consider a standard $C^0$-conforming isogeometric multi-patch discretization. Each subdomain $\Omega_k$ is defined as the image of the parametric domain $\H{\Omega} = (0,1)^2$ under a spline geometry mapping $\bm{F}_k \colon \H{\Omega} \to \Omega_k$. 
\begin{assumption}\label{ass:unif_reg}
We assume that $\bm{F}_k \in \Sp^{\mathbf{p}}_k \times \Sp^{\mathbf{p}}_k$, satisfies the uniform regularity condition: there exists a constant $\gamma >0$ such that
\begin{equation*}
	\det \left(\nabla \bm{F}_k({ \H{\mathbf{x}}})\right) \geq \gamma, \quad \forall {\H{\mathbf{x}}} \in \H{\Omega}. 
\end{equation*}
\end{assumption}
Here, $\Sp^{\mathbf{p}}_k = \Sp^{\mathbf{p}}(\Xi^k_1, \Xi^k_2)$ is a spline space of fixed degree $\mathbf{p}$, with (potentially distinct) quasi-uniform knot vectors in each parametric direction.\\
The isogeometric discretization of the patch-local primal space on $\Omega_k$ is defined as
\begin{equation}\label{eq:defXkh}
	X_{k,h} := \Set{ v \in H^2(\Omega_k) : v \circ \bm{F}_{k} \in \Sp^{\mathbf{p}}_k \text{ and } v = \dn v = 0 \text{ on } \partial \Omega \cap \partial \Omega_k }.
\end{equation}
To simplify the notation, we assume that the spline degree $p$ is the same in every direction and across all patches.
To construct a global space, we assume that the spaces are $C^0$-conforming, as outlined below:
\begin{assumption}\label{ass:confmesh}
	We assume that, for all $\ell \in \is$, the trace spaces on $\Gamma_\ell$ are equal, i.e., $X_{\ml,h}|_{\Gamma_\ell} = X_{\sell,h}|_{\Gamma_\ell}$.
\end{assumption}
{ This assumption implies that the knot vectors along the interface for the spline spaces of the two adjacent patches must be matching.}
We denote by $\mathcal{I}_{\mathcal{V},k}$ the set of indices corresponding to the four vertices of the patch $\Omega_k$. The vertices in the parametric domain $\H{\Omega}$ are denoted by $\H{\mathcal{V}}_{i,k} := \bm{F}_k^{-1}(\mathcal{V}_i)$, or simply $\H{\mathcal{V}}_{i} = \H{\mathcal{V}}_{i,k}$ when $k$ is fixed.\\
On the global domain $\Omega$, we define the space of globally $C^0$-continuous isogeometric functions that are $C^2$-continuous at the vertices of the decomposition as
\begin{equation}\label{eq:defXh}
	X_h := \left\{ v \in H^1(\Omega) : v|_{\Omega_k} \in X_{k,h}, \, \forall k \in \mathcal{I}_{\Omega}, \text{ and } v \in C^2(\mathcal{V}_i), \forall i \in \mathcal{I}_{\mathcal{V}} \right\}.
\end{equation}
\begin{prop}
	The inclusion $X_h \subset X$ holds.
\end{prop}
\begin{proof}
	Let $v \in X_h$. By Assumption~\ref{ass:reg}, we have $\jump{\dn v}_{\gl} \in C^0(\gl)$. Additionally, since $v \in C^2(\mathcal{V}_i)$ for all $i \in \iv$, it follows that $\jump{\dn v}_{\gl}$ vanishes at the endpoints of $\gl$. Consequently, $\jump{\dn v}_{\gl} \in H^{1/2}_{00}(\gl)$, which implies $v \in X$.
\end{proof}
Let $M_h := \prod_{\ell \in  \is} M_\ell$ denote a finite-dimensional subspace of $M$ supported on $\Gamma_{\ell}$, that will be specified in Section \ref{sec:dscinfsup}. We define the mortar kernel as  
\begin{equation}\label{eq:mortarkernel}
	V_h := \Set{v_h \in X_h : b(v_h, \mu_h) = 0, \forall \mu_h \in M_h}.
\end{equation}  
Here, $V_h$ represents the isogeometric space with weak $C^1$-continuity, which depends on $M_h$.  
With this framework, we can formulate the weak form of the model problem as a saddle-point problem incorporating $C^1$-mortar constraints along $\Sigma$:  
\textit{Find $(u_h, \lambda_h) \in X_h \times M_h$ such that}  
\begin{equation}\label{eq:qh}
	\begin{cases}
		\begin{aligned}
			a(u_h, v_h) + b(v_h, \lambda_h)	 & = \langle f, v_h \rangle_{H^{-1} \times H^1_0},  && \forall   v_h \in X_h, \\
			b(u_h, \mu_h)	                 & = 0,        && \forall \mu_h \in M_h.
		\end{aligned}
	\end{cases}
\end{equation}
	\section{Approximation properties of $X_h$}\label{sec:Xhprop}
{ In this section we provide approximation error bounds for the globally $C^0$ multi-patch isogeometric space $X_h$. We first define patchwise projectors on the parametric and physical domain and show their approximation properties, which then yield the approximation properties of the multi-patch projector.

\subsection{Approximation properties on the parametric domain}}

For $s > 3$, the following Sobolev compact embedding holds (see~\citeRefDetail{Theorem 6.2}{adams}):
	\begin{equation*}
		H^s(\H{\Omega}) \subset C^2\left(\overline{\H{\Omega}}\right).
	\end{equation*}
Given an integer $s \geq 0$, the bent Sobolev space of order $s$, denoted by $\hh^s(\H{\Omega})$, consists of functions that belong to $H^s(Q)$ for every element $Q$ in the mesh $\Q_h$ introduced in \eqref{eq:mesh}, and  exhibit matching traces at interelement interfaces, up to $\min\{p-m, s-1\}$-order derivatives, where $m$ represents the multiplicity of the knot separating two neighboring elements. The seminorms and norms of $\hh^s(\H{\Omega})$ are defined as follows:
\begin{equation*}
	\begin{split}
		\abs{v}^2_{\hh^i(\H{\Omega})}  &:= \sum_{Q  \in \Q_h} \abs{v}^2_{H^i(Q)}, \quad 0 \leq i \leq s,\\
		\norm{v}^2_{\hh^s(\H{\Omega})} &:= \sum_{i = 0}^{s} \abs{v}^2_{\hh^i(\H{\Omega})}.
	\end{split}
\end{equation*}
Bent Sobolev spaces on { parts of the parameter domain} $\H{\Omega}$ { as well as on the physical domain} $\Omega$ are constructed analogously. For further details, see \citeOneRef{bazilevs}. 
{
\begin{prop}\label{prop:dual-basis}
Let ${\cal B}_{\mathbf p}(\boldsymbol{\Xi}^k) = \{ \H{B}_{{\mathbf i}, {\mathbf p}}, (1,1)\leq \mathbf i \leq (n_1,n_2) \}$ be the tensor-product B-spline basis of the spline space $\Sp_{k}^{\mathbf{p}}$. Then there exists a dual basis $\{\LambdX_{\mathbf i}, (1,1)\leq \mathbf i \leq (n_1,n_2)\}$, which satisfies
\begin{align}\label{eq:lambdaprop}
	\abs{\LambdX_{(i_1,i_2)} (\H{v})} &\leq Ch^{-1} \norm{\H{v}}_{L^2\left( S^1_{i_1} \times S^2_{i_2} \right) }, & (2,2) \leq (i_1,i_2) \leq (n_1,n_2)-(1,1), \notag \\
	\abs{\LambdX_{(i_1,1)} (\H{v})} &\leq Ch^{-\nicefrac{1}{2}} \norm{\H{v}}_{L^2\left( S^1_{i_1}\times\{0\} \right) }, & 2 \leq i_1 \leq n_1-1, \notag \\
	\abs{\LambdX_{(i_1,n_2)} (\H{v})} &\leq Ch^{-\nicefrac{1}{2}} \norm{\H{v}}_{L^2\left( S^1_{i_1}\times\{1\} \right) }, & 2 \leq i_1 \leq n_1-1, \\
	\abs{\LambdX_{(1,i_2)} (\H{v})} &\leq Ch^{-\nicefrac{1}{2}} \norm{\H{v}}_{L^2\left( \{0\}\times S^2_{i_2} \right) }, & 2 \leq i_2 \leq n_2-1, \notag \\
	\abs{\LambdX_{(n_1,i_2)} (\H{v})} &\leq Ch^{-\nicefrac{1}{2}} \norm{\H{v}}_{L^2\left( \{1\}\times S^2_{i_2} \right) }, & 2 \leq i_2 \leq n_2-1, \notag
\end{align}
and
\begin{align}
    &\abs{\LambdX_{(1,1)} (\H{v})} = \abs{\H{v}(0,0)}, \quad 
    &&\abs{\LambdX_{(1,n_2)} (\H{v})} = \abs{\H{v}(0,1)}, \notag \\ 
    &\abs{\LambdX_{(n_1,1)} (\H{v})} = \abs{\H{v}(1,0)}, \quad 
    &&\abs{\LambdX_{(n_1,n_2)} (\H{v})} = \abs{\H{v}(1,1)}, \notag
\end{align}
where $S^\delta_i = \supp{\H{B}^\delta_{i, p}}$ are the supports of the univariate B-splines in direction $\delta\in\{1,2\}$.
\end{prop}
\begin{proof}
We define a dual basis $\{\LambdX^\delta_i\}_{1\leq i \leq n_\delta}$ for the univariate B-splines in direction $\delta$, where $\LambdX^\delta_1(\H{v}) = \H{v}(0)$, $\LambdX^\delta_{n_\delta}(\H{v}) = \H{v}(1)$ and $\LambdX^\delta_i(\H{v})$ for $2 \leq i \leq n_\delta-1$ are defined as in~\citeRefDetail{Theorem 4.41}{schumaker}, which states
\[
    \abs{\LambdX^\delta_i (f)} \leq C h^{-\nicefrac{1}{2}} \norm{f}_{L^2\left( S^\delta_i\right)}.
\]
By defining $\LambdX_{(i_1,i_2)} = \LambdX^1_{i_1} \otimes \LambdX^2_{i_2}$, the bounds follow directly from~\citeRefDetail{Theorem 4.41}{schumaker}, and the fact that $\abs{\LambdX^\delta_{1} (f)} = \abs{f(0)}$ and $\abs{\LambdX^\delta_{n_\delta} (f)} = \abs{f(1)}$. This completes the proof.
\end{proof}
See Figure~\ref{fig:lambdabasis} for a visualization of the dual basis defined in Proposition~\ref{prop:dual-basis}.
\begin{figure}[h]
\centering

\subfloat[ Index sets for $\LambdX_{(i_1,i_2)}$\label{fig:lambdabasis-a}]{
\begin{tikzpicture}
\begin{axis}[hide axis, axis equal,
			 xmin=-0.2, xmax=1.2, ymin=-0.2, ymax=1.2]
\addplot[black]coordinates{(0, 0) ( 1, 0) ( 1, 1) ( 0, 1) ( 0, 0)};
\addplot[gray, dashed]coordinates{(1/8, 0) (1/8, 1)};
\addplot[gray, dashed]coordinates{(2/8, 0) (2/8, 1)};
\addplot[gray, dashed]coordinates{(3/8, 0) (3/8, 1)};
\addplot[gray, dashed]coordinates{(4/8, 0) (4/8, 1)};
\addplot[gray, dashed]coordinates{(5/8, 0) (5/8, 1)};
\addplot[gray, dashed]coordinates{(6/8, 0) (6/8, 1)};
\addplot[gray, dashed]coordinates{(7/8, 0) (7/8, 1)};
\addplot[gray, dashed]coordinates{(0, 1/8) (1, 1/8)};
\addplot[gray, dashed]coordinates{(0, 2/8) (1, 2/8)};
\addplot[gray, dashed]coordinates{(0, 3/8) (1, 3/8)};
\addplot[gray, dashed]coordinates{(0, 4/8) (1, 4/8)};
\addplot[gray, dashed]coordinates{(0, 5/8) (1, 5/8)};
\addplot[gray, dashed]coordinates{(0, 6/8) (1, 6/8)};
\addplot[gray, dashed]coordinates{(0, 7/8) (1, 7/8)};
\addplot+[only marks, mark = ball, ball color = green, draw opacity=0, mark size = 3pt]
	coordinates{(0,0) (1,0) (1,1) (0,1)};
				
\addplot+[only marks, mark = ball, ball color = blue, draw opacity=0, mark size = 3pt]
	coordinates{(1/8,0) (2/8,0) (3/8,0) (4/8,0) (5/8,0) (6/8,0) (7/8,0)
	(0/8,1/8) (8/8,1/8)
	(0/8,2/8) (8/8,2/8)
	(0/8,3/8) (8/8,3/8)
	(0/8,4/8) (8/8,4/8)
	(0/8,5/8) (8/8,5/8)
	(0/8,6/8) (8/8,6/8)
	(0/8,7/8) (8/8,7/8)
	(1/8,8/8) (2/8,8/8) (3/8,8/8) (4/8,8/8) (5/8,8/8) (6/8,8/8) (7/8,8/8)};				
				
\addplot+[only marks, mark = ball, ball color = black,  draw opacity=0, mark size = 3pt]
	coordinates{(2/8,2/8) (3/8,2/8) (4/8,2/8) (5/8,2/8) (6/8,2/8) (2/8,3/8) (3/8,3/8) (4/8,3/8) (5/8,3/8) (6/8,3/8) (2/8,4/8) (3/8,4/8) (4/8,4/8) (5/8,4/8) (6/8,4/8) (2/8,5/8) (3/8,5/8) (4/8,5/8) (5/8,5/8) (6/8,5/8) (2/8,6/8) (3/8,6/8) (4/8,6/8) (5/8,6/8) (6/8,6/8)};
	
\addplot+[only marks, mark = ball, ball color = red,  draw opacity=0, mark size = 3pt]
coordinates{(1/8,1/8) (2/8,1/8) (3/8,1/8) (4/8,1/8) (5/8,1/8) (6/8,1/8) (7/8,1/8)
	(1/8,2/8) (7/8,2/8) (1/8,3/8) (7/8,3/8) (1/8,4/8) (7/8,4/8) (1/8,5/8) (7/8,5/8) (1/8,6/8) (7/8,6/8) (1/8,7/8) (2/8,7/8) (3/8,7/8) (4/8,7/8) (5/8,7/8) (6/8,7/8) (7/8,7/8)};

\node[below]  at (axis cs: 0,0) {$(1,1)$};
\node[above]  at (axis cs: 0,1) {$(1,n_2)$};
\node[below] at (axis cs: 1,0) {$(n_1,1)$};
\node[above] at (axis cs: 1,1) {$(n_1,n_2)$};

\end{axis}
\end{tikzpicture}}
\subfloat[Index sets for ${\LambdY}_{(i_1,i_2)}$\label{fig:lambdabasis-b}]{
\begin{tikzpicture}
\begin{axis}[hide axis, axis equal,
			 xmin=-0.2, xmax=1.2, ymin=-0.2, ymax=1.2]
\addplot[black]coordinates{(0, 0) ( 1, 0) ( 1, 1) ( 0, 1) ( 0, 0)};
\addplot[gray, dashed]coordinates{(1/8, 0) (1/8, 1)};
\addplot[gray, dashed]coordinates{(2/8, 0) (2/8, 1)};
\addplot[gray, dashed]coordinates{(3/8, 0) (3/8, 1)};
\addplot[gray, dashed]coordinates{(4/8, 0) (4/8, 1)};
\addplot[gray, dashed]coordinates{(5/8, 0) (5/8, 1)};
\addplot[gray, dashed]coordinates{(6/8, 0) (6/8, 1)};
\addplot[gray, dashed]coordinates{(7/8, 0) (7/8, 1)};
\addplot[gray, dashed]coordinates{(0, 1/8) (1, 1/8)};
\addplot[gray, dashed]coordinates{(0, 2/8) (1, 2/8)};
\addplot[gray, dashed]coordinates{(0, 3/8) (1, 3/8)};
\addplot[gray, dashed]coordinates{(0, 4/8) (1, 4/8)};
\addplot[gray, dashed]coordinates{(0, 5/8) (1, 5/8)};
\addplot[gray, dashed]coordinates{(0, 6/8) (1, 6/8)};
\addplot[gray, dashed]coordinates{(0, 7/8) (1, 7/8)};
\addplot+[only marks, mark = ball, ball color = cyan, draw opacity=0, mark size = 3pt]
	coordinates{(0,0) (1/8,0) (2/8,0) (6/8,0) (7/8,0) (1,0) (0,1) (1/8,1) (2/8,1) (6/8,1) (7/8,1) (1,1) (0,1/8) (0,2/8) (0,6/8) (0,7/8)
	  (1,1/8) (1,2/8) (1,6/8) (1,7/8) (1/8,1/8) (7/8,1/8) (1/8,7/8) (7/8,7/8)};
				
\addplot+[only marks, mark = ball, ball color = blue, draw opacity=0, mark size = 3pt]
	coordinates{(3/8,0) (4/8,0) (5/8,0) 
	(0/8,3/8) (8/8,3/8)
	(0/8,4/8) (8/8,4/8)
	(0/8,5/8) (8/8,5/8)
	 (3/8,8/8) (4/8,8/8) (5/8,8/8) };				
				
\addplot+[only marks, mark = ball, ball color = black,  draw opacity=0, mark size = 3pt]
	coordinates{(2/8,2/8) (3/8,2/8) (4/8,2/8) (5/8,2/8) (6/8,2/8) (2/8,3/8) (3/8,3/8) (4/8,3/8) (5/8,3/8) (6/8,3/8) (2/8,4/8) (3/8,4/8) (4/8,4/8) (5/8,4/8) (6/8,4/8) (2/8,5/8) (3/8,5/8) (4/8,5/8) (5/8,5/8) (6/8,5/8) (2/8,6/8) (3/8,6/8) (4/8,6/8) (5/8,6/8) (6/8,6/8)};
	
\addplot+[only marks, mark = ball, ball color = red,  draw opacity=0, mark size = 3pt]
coordinates{ (2/8,1/8) (3/8,1/8) (4/8,1/8) (5/8,1/8) (6/8,1/8) 
	(1/8,2/8) (7/8,2/8) (1/8,3/8) (7/8,3/8) (1/8,4/8) (7/8,4/8) (1/8,5/8) (7/8,5/8) (1/8,6/8) (7/8,6/8) (2/8,7/8) (3/8,7/8) (4/8,7/8) (5/8,7/8) (6/8,7/8)};

\node[below]  at (axis cs: 0,0) {$(1,1)$};
\node[above]  at (axis cs: 0,1) {$(1,n_2)$};
\node[below] at (axis cs: 1,0) {$(n_1,1)$};
\node[above] at (axis cs: 1,1) {$(n_1,n_2)$};

\end{axis}
\end{tikzpicture}}

\caption{ The dual functionals ${\LambdX}_{(i_1,i_2)}$ appearing in Proposition~\ref{prop:dual-basis} and ${\LambdY}_{(i_1,i_2)}$ appearing in the proof of Proposition~\ref{prop:projparamdom} have a different support depending on $i_1,i_2$. For the indices in black and red, the supports are those of the corresponding basis functions. For the blue indices the supports are proper intervals of the boundary edges $\partial \widehat{\Omega}$. For the indices depicted in green and cyan, the supports of the dual functionals are the corresponding vertices, where values (green) and up to second order derivatives (cyan), respectively, are interpolated.}
\label{fig:lambdabasis}

\end{figure}
}
{
To prove error bounds we require suitable trace estimates.
\begin{lemma}
Let $B\subset \mathbb{R}^2$ be a box with dimensions that are $O(h)$. We then have for all edges $\Gamma\subset \partial B$ and for all $f \in H^1(B)$
\begin{equation}\label{eq:trace_estimate}
	\norm{f}_{L^2\left(  \Gamma \right)} \leq C \left( h^{-\nicefrac{1}{2}} \norm{f}_{L^2(B)} + h^{\nicefrac{1}{2}} \abs{f}_{H^1(B)} \right).
\end{equation}
Moreover, for all vertices $\mathbf{x}$ of $B$ and for all $f \in H^2(B)$ we have
\begin{equation}\label{eq:corner_estimate}
\begin{aligned}
	\abs{f(\mathbf{x})} \leq C \left( h^{-1}\norm{f}_{L^2(B)} + \abs{f}_{H^1(B)} + h \abs{f}_{H^2(B)} \right).
\end{aligned}
\end{equation}
\end{lemma}
\begin{proof}
Let us consider a linear parametrization $T:[0,1]^2 \rightarrow B$ of $B$, i.e. $T([0,1]^2)=B$. Thanks to a scaling argument, we have
\begin{equation*}
\begin{aligned}
	\norm{f}_{L^2\left( \Gamma \right)} &\leq C h^{\nicefrac{1}{2}} \norm{f \circ T}_{L^2(T^{-1}(\Gamma))} \leq C h^{\nicefrac{1}{2}} \norm{f \circ T}_{H^1([0,1]^2)} \\
	&= C h^{\nicefrac{1}{2}} \left( \norm{f \circ T}_{L^2([0,1]^2)} + \abs{f \circ T}_{H^1([0,1]^2)} \right)\\
	&\leq C h^{\nicefrac{1}{2}} \left( h^{-1} \norm{f}_{L^2(B)} +  \abs{f}_{H^1(B)} \right) \\
	&= C \left( h^{-\nicefrac{1}{2}} \norm{f}_{L^2(B)} + h^{\nicefrac{1}{2}} \abs{f}_{H^1(B)} \right)
\end{aligned}
\end{equation*}
and, for $\mathbf{x}=T(\hat{\mathbf{x}}) \in \partial\Gamma$, 
\begin{equation*}
\begin{aligned}
	\abs{f(\mathbf{x})} &= \abs{f(T (\hat{\mathbf{x}}))} \leq C \left( \norm{f \circ T}_{L^2(T^{-1}({\Gamma}))} + \abs{f \circ T}_{H^1(T^{-1}({\Gamma}))} \right) \\
	&\leq C \left( \norm{f \circ T}_{L^2(T^{-1}(B))} + \abs{f \circ T}_{H^1(T^{-1}(B))} + \abs{f \circ T}_{H^2(T^{-1}(B))} \right) \\
	&\leq C \left( h^{-1}\norm{f}_{L^2(B)} + \abs{f}_{H^1(B)} + h \abs{f}_{H^2(B)} \right).
\end{aligned}
\end{equation*}
This completes the proof.
\end{proof}
We can now define a projector onto the tensor-product spline space which interpolates up to second order derivatives at all vertices and satisfies suitable local approximation error bounds.}
\begin{prop}\label{prop:projparamdom}
Let $m$ and $t$ be integers such that $0 \leq t \leq m \leq p+1$. There exists an operator $\H{\Pi}_{k,h} \colon \hh^4(\H{\Omega}) \to \Sp_{k}^{\mathbf{p}}$ such that
\begin{equation*}\label{eq:defproj}
	\H{\Pi}_{k,h} \H{v}_h = \H{v}_h, \quad \forall \, \H{v}_h \in \Sp_{k}^{\mathbf{p}},
\end{equation*}
and for all $\iota \in \mathcal{I}_{\mathcal{V}, k}$ { we have}
\begin{equation}\label{eq:projc2}
\begin{split}
	(\H{\Pi}_{k,h} \H{v})(\H{\mathcal{V}}_\iota) &= \H{v}(\H{\mathcal{V}}_\iota), \\
	\nabla (\H{\Pi}_{k,h} \H{v})(\H{\mathcal{V}}_\iota) &= \nabla \H{v}(\H{\mathcal{V}}_\iota), \\
	\nabla (\nabla (\H{\Pi}_{k,h} \H{v}))(\H{\mathcal{V}}_\iota) &= \nabla(\nabla \H{v})(\H{\mathcal{V}}_\iota),
\end{split}
\end{equation}
{ for all $\H{v}\in \hh^4(\H{\Omega})$}. Furthermore,  for  $s = \max(4, m)$ and  for all $\H{v} \in \hh^s(\H{\Omega})$, we have
\begin{equation}\label{eq:appproj}
	\norm*{\H{v} - \H{\Pi}_{k,h} \H{v}}_{\hh^t(\H{\Omega})} \leq C \, h^{m-t} \norm*{\H{v}}_{\hh^s(\H{\Omega})}.
\end{equation}
\end{prop}
\begin{proof} {
The operator $\H{\Pi}_{k,h}$ is defined as
\begin{equation}
	\H{\Pi}_{k,h}\H{v} = \sum_{i_1=1}^{n_1}\sum_{i_2=1}^{n_2} {\LambdY}_{\mathbf{i}}(\H{v}) \H{B}_{\mathbf{i},\mathbf{p}},
\end{equation}
using a modified version ${\LambdY}_{\mathbf{i}}$ of the dual basis $\LambdX_{\mathbf{i}}$ introduced in Proposition~\ref{prop:dual-basis}. Let $\mathcal{I}_v := \Set{(1,1),(n_1,1),(1,n_2),(n_1,n_2)}$ and let 
\begin{equation*}
	\ic := \Set{ \textbf{i}=(i_1,i_2) \in \mathcal{I}: \abs{\iota_1-i_1}+\abs{\iota_2-i_2} \leq 2, \mbox{ for some }(\iota_1,\iota_2)\in \mathcal{I}_v }.
\end{equation*}
This index set corresponds to exactly those functions where at least one of the up to second order derivatives at the corners is non-zero.
For all $\mathbf{i} \in \mathcal{I} \setminus \ic$, we define $\LambdY_{\mathbf{i}}$ as
\begin{equation*}
	\LambdY_{\mathbf{i}} := \LambdX_{\mathbf{i}}.
\end{equation*}
For all $\mathbf{i} \in \ic$, the functionals $\LambdY_{\mathbf{i}}$ are defined to fit the following constraints
\begin{align}
	\H{\Pi}_{k,h} \H{v}(\H{\mathcal{V}}_\iota) &= \H{v}(\H{\mathcal{V}}_\iota), & \partial_{\H{x}}\partial_{\H{x}} \H{\Pi}_{k,h} \H{v}(\H{\mathcal{V}}_\iota) &= \partial_{\H{x}}\partial_{\H{x}} \H{v}(\H{\mathcal{V}}_\iota) \notag\\
	\partial_{\H{x}} \H{\Pi}_{k,h} \H{v}(\H{\mathcal{V}}_\iota) &= \partial_{\H{x}} \H{v}(\H{\mathcal{V}}_\iota) & \partial_{\H{x}}\partial_{\H{y}} \H{\Pi}_{k,h} \H{v}(\H{\mathcal{V}}_\iota) &= \partial_{\H{x}}\partial_{\H{y}} \H{v}(\H{\mathcal{V}}_\iota) \label{eq:const}\\
	\partial_{\H{y}} \H{\Pi}_{k,h} \H{v}(\H{\mathcal{V}}_\iota) &= \partial_{\H{y}} \H{v}(\H{\mathcal{V}}_\iota) & \partial_{\H{y}}\partial_{\H{y}} \H{\Pi}_{k,h} \H{v}(\H{\mathcal{V}}_\iota) &= \partial_{\H{y}}\partial_{\H{y}} \H{v}(\H{\mathcal{V}}_\iota), \notag
\end{align}
for all $\iota \in \mathcal{I}_{\mathcal{V},k}$.}

It is easy to see that~\eqref{eq:const} is a well-posed interpolation problem that determines { $\LambdY_{\mathbf{i}}$ for all $\mathbf{i} \in \ic$}. In particular, for $\H{\mathcal{V}}_\iota = (0,0)$ and { $\bigl\{ \LambdY_{\mathbf{i}} \bigr\}_{i_1+i_2 \leq 4}$}, we obtain
%
{
\begin{equation}\label{eq:interp_system}
	\begin{bmatrix}
	1 & 0 & 0 & 0 & 0 & 0 \\
	-1 & 1 & 0 & 0 & 0 & 0 \\
	-1 & 0 & 1 & 0 & 0 & 0 \\
	c_1 & -1-c_1 & 0 & 1 & 0 & 0 \\
	1 & -1 & -1 & 0 & 1 & 0 \\
	c_2 & 0 & -1-c_2 & 0 & 0 & 1
	\end{bmatrix}
	\begin{bmatrix}
	\LambdY_{(1,1)}\\
	\LambdY_{(2,1)}\\
	\LambdY_{(1,2)}\\
	\LambdY_{(3,1)}\\
	\LambdY_{(2,2)}\\
	\LambdY_{(1,3)}
	\end{bmatrix}
	=
	\begin{bmatrix}
	\H{v}(0,0)\\
	\frac{h_{1,1}}{p_1} \, \partial_{\H{x}} \H{v}(0,0)\\
	\frac{h_{2,1}}{p_2} \, \partial_{\H{y}} \H{v}(0,0)\\
	\frac{(h_{1,1})^2c_1}{p_1(p_1-1)} \, \partial_{\H{x}}\partial_{\H{x}} \H{v}(0,0)\\
	\frac{h_{1,1} h_{2,1}}{p_1 p_2} \, \partial_{\H{x}}\partial_{\H{y}} \H{v}(0,0)\\
	\frac{(h_{2,1})^2c_2}{p_2(p_2-1)} \, \partial_{\H{y}}\partial_{\H{y}} \H{v}(0,0)
	\end{bmatrix},
\end{equation}
with $c_k:=\frac{h_{k,1}+h_{k,2}}{h_{k,1}}$ if the first knot in direction $k$ is a single knot, otherwise $c_k:=1$. In the above system the} matrix is non-singular, since it is lower-triangular. Thanks to Assumption~\ref{ass:quasiunif}, if $\H{v} \in \hh^4(\H{\Omega})$, { we obtain
\begin{equation}\label{eq:1}
	\sum_{i_1+i_2 \leq 4} \abs{\LambdY_{\mathbf{i}}\H{v}} \leq C \left( \abs{\H{v}(0,0)} + h \norm{\nabla \H{v}(0,0)} +  h^2 \norm{\nabla(\nabla \H{v})(0,0)}\right),
\end{equation}
where $\norm{\cdot}$ denotes the vector or matrix norm.
A similar bound holds for all other $\LambdY_{\mathbf{i}}$ with $\mathbf{i} \in \ic$. Assume now that $Q \in \Q_h$ satisfies  $Q \subseteq \bigcup_{i_1 + i_2 \leq 4} S^1_{i_1}\times S^2_{i_2}$, then we have
\begin{equation}\label{eq:11}
	\begin{split}
	\norm{\H{\Pi}_{k,h}(\H{v})}_{L^2(Q)} &\leq \sum_{(i_1,i_2):\ Q \subseteq S^1_{i_1}\times S^2_{i_2}} \abs{\LambdY_{\mathbf{i}}\H{v}} \norm{\H{B}_{\mathbf{i},\mathbf{p}}}_{L^2(Q)} \\
	&\leq C \sum_{\substack{(i_1,i_2):\ Q \subseteq S^1_{i_1}\times S^2_{i_2} \\ i_1+i_2 > 4}} h \abs{\LambdY_{\mathbf{i}}\H{v}} + C \sum_{ i_1+i_2 \leq 4} h \abs{\LambdY_{\mathbf{i}}\H{v}}.
	\end{split}
\end{equation}}
Then, from~\eqref{eq:lambdaprop}, \eqref{eq:1} and~\eqref{eq:11}, we conclude that
\begin{equation}\label{eq:12}
	\begin{aligned}
		\norm{\H{\Pi}_{k,h} \H{v}}_{L^2(Q)} 
		&\leq C \left(
		\norm{\H{v}}_{L^2(\widetilde{Q})} 
		+ h^{\nicefrac{1}{2}} \norm{\H{v}}_{L^2(\partial \widetilde{Q} \cap \partial \H{\Omega})} \right. \\
		&\left.\qquad \quad 
		+ h \abs{\H{v}(0,0)} 
		+ h^2 \norm{\nabla \H{v}(0,0)} 
		+ h^3 \norm{\nabla(\nabla \H{v})(0,0)}
		\right) \\
		&= C \opnorm{\H{v}}_{\widetilde{Q}},
	\end{aligned}
\end{equation}
where $\widetilde{Q}$ is the {support extension} of $Q$, that is the union of the supports of basis functions whose support intersects $Q$, and $\opnorm{\cdot}_{\widetilde{Q}}$ is a short notation for the sum of norms appearing on the right-hand-side of~\eqref{eq:12}. { The estimate~\eqref{eq:trace_estimate} immediately implies
\begin{equation*}
	\norm{\H{v}}_{L^2\left(  \partial \widetilde{Q} \cap \partial \H{\Omega} \right)} \leq C \left( h^{-\nicefrac{1}{2}} \norm{f}_{L^2(\widetilde{Q})} + h^{\nicefrac{1}{2}} \abs{f}_{\hh^1(\widetilde{Q})} \right).
\end{equation*}
Moreover, \eqref{eq:corner_estimate} applied to $\H{v}$, or components of $\nabla \H{v}$ and $\nabla(\nabla \H{v})$ implies
\begin{equation*}
	\abs{\H{v}(0,0)} \leq C \left( h^{-1}\norm{\H{v}}_{L^2(\widetilde{Q})} + \abs{\H{v}}_{\hh^1(\widetilde{Q})} + h \abs{\H{v}}_{\hh^2(\widetilde{Q})} \right),
\end{equation*}
\begin{equation*}
	\abs{\nabla \H{v}(0,0)} \leq C \left( h^{-1}\abs{\H{v}}_{\hh^1(\widetilde{Q})} + \abs{\H{v}}_{\hh^2(\widetilde{Q})} + h \abs{\H{v}}_{\hh^3(\widetilde{Q})} \right)
\end{equation*}
and
\begin{equation*}
	\abs{\nabla(\nabla \H{v})(0,0)} \leq C \left( h^{-1}\abs{\H{v}}_{\hh^2(\widetilde{Q})} + \abs{\H{v}}_{\hh^3(\widetilde{Q})} + h \abs{\H{v}}_{\hh^4(\widetilde{Q})} \right).
\end{equation*}
}
In summary, we obtain
\begin{equation*}\label{eq:trinorm}
	\opnorm{\H{v}}_{\widetilde{Q}} \leq C \left( \norm{\H{v}}_{L^2(\widetilde{Q})} + h \norm{\H{v}}_{\hh^1(\widetilde{Q})} + h^2 \norm{\H{v}}_{\hh^2(\widetilde{Q})} + h^3 \norm{\H{v}}_{\hh^3(\widetilde{Q})} + h^4 \norm{\H{v}}_{\hh^4(\widetilde{Q})}  \right).
\end{equation*}
{ Using the same proving strategy as in~\citeRefDetail{Lemma 3.1}{bazilevs}, based on a scaling argument, and replacing the Sobolev seminorm on the left by a sum of correctly scaled norms,}
we can prove that for any $\H{v} \in \hh^s (\widetilde{Q})$, there exists a spline $\H{w} \in \Sp_{k}^{\mathbf{p}}$ such that
\begin{equation*}
	\opnorm{\H{v} - \H{w}}_{\widetilde{Q}} + h^t \abs{\H{v} - \H{w}}_{\hh^t (\widetilde{Q})} \leq C\,h^m \norm{\H{v}}_{\hh^s (\widetilde{Q})}. 
\end{equation*}
Combining the previous results and a standard inverse inequality, when $Q \subseteq \cup_{i_1 + i_2 \leq 4} S^1_{i_1}\times S^2_{i_2}$
\begin{equation}\label{eq:hnorm}
\begin{split}
	\abs{\H{v} - \H{\Pi}_{k,h} \H{v}}_{\hh^t(Q)} &= \abs{\H{v} - \H{w} - \H{\Pi}_{k,h} (\H{v} - \H{w})}_{\hh^t(Q)} \\
	&\leq \abs{\H{v} - \H{w}}_{\hh^t(Q)} + C\,h^{-t} \norm{\H{\Pi}_{k,h} (\H{v} - \H{w})}_{L^2(Q)} \\
	&\leq \abs{\H{v} - \H{w}}_{\hh^t(Q)} + C\,h^{-t} \opnorm{\H{v} - \H{w}}_{\widetilde{Q}} \\
	&\leq C\,h^{m-t} \norm{\H{v}}_{\hh^s (\widetilde{Q})}.
\end{split}
\end{equation}
We remark that the same construction can be repeated for all $\H{\mathcal{V}}_\iota$, with $\iota \in \mathcal{I}_{\mathcal{V},k}$, and~\eqref{eq:hnorm} holds for $Q \nsubseteq \cup_{i_1 + i_2 \leq 4} S^1_{i_1}\times S^2_{i_2}$, as proved in~\citeOneRef{bazilevs}. Finally~\eqref{eq:appproj} follows by summation over all mesh elements.
\end{proof}
{
We moreover have the following error bound for all traces.
\begin{prop}\label{prop:projparamdom-trace}
Let $m$ be an integer such that $0\leq m\leq p+1$, let $\H{\Pi}_{k,h}$ be the operator from Proposition~\ref{prop:projparamdom}
and consider an edge $\H{\Gamma} \subset \partial \H{\Omega}$. Then for  $s = \max(4, m)$ and  for all $\H{v} \in \hh^s(\H{\Omega})$, we have
\begin{equation*}
	\norm{\nabla(\H{v} - \H{\Pi}_{k,h} \H{v})}_{L^2(\H{\Gamma})} + h\, \abs{\nabla(\H{v} - \H{\Pi}_{k,h} \H{v})}_{\hh^1(\H{\Gamma})} \leq C\, h^{m - \nicefrac{3}{2}} \norm{\H{v}}_{\hh^s(\H{\Omega})}.
\end{equation*}
\end{prop}
\begin{proof}
We proceed similar to the proof of Proposition~\ref{prop:projparamdom}. Given a mesh element $Q$ adjacent to the boundary $\partial \H{\Omega}$, let $E$ be one edge of $\partial Q$ on $\partial \H{\Omega}$. Setting $\eta := \H{v} - \H{\Pi}_{k,h} \H{v}$ and using~\eqref{eq:trace_estimate} and~\eqref{eq:hnorm}, we obtain for $j \in\{0,1\}$ 
\begin{equation*}
\begin{aligned}
	\abs{\nabla \eta}_{H^j(E)} &\leq C \left( h^{-\nicefrac{1}{2}} \abs{\nabla \eta}_{H^{j}(Q)} + h^{\nicefrac{1}{2}} \abs{\nabla \eta}_{H^{j+1}(Q)} \right) \\
	&\leq C \left( h^{-\nicefrac{1}{2}} \abs{\eta}_{H^{j+1}(Q)} + h^{\nicefrac{1}{2}} \abs{\eta}_{H^{j+2}(Q)} \right) \\
	&\leq C h^{m-\nicefrac{5}{2}} \norm{\H{v}}_{\hh^s(\widetilde{Q})}.
\end{aligned}
\end{equation*}
Then the result follows follows.
\end{proof}
Note that $\nabla \H{v} -  \nabla \H{\Pi}_{k,h} \H{v} \in H^1_0(\H{\Gamma})$ for each edge $\H{\Gamma}$ of $\partial \H{\Omega}$. Indeed $\nabla \H{v}$ is interpolated at the vertices, that is $\nabla \H{v} -  \nabla \H{\Pi}_{k,h} \H{v}$ vanishes at the endpoints of $\H{\Gamma}$.

With homogeneous boundary conditions, we have a result similar to Proposition~\ref{prop:projparamdom-trace}
under weaker regularity conditions}, as stated below.
\begin{prop}\label{prop:projparamdom_low-H2_regularity}
Let $\Gamma$ be an edge of $\Omega_k$ and $\H{\Gamma}$ its pre-image (an edge of $ \H{\Omega}$). There exists an operator 
\begin{equation*}
	\H{\Pi}^0_{k,h} \colon
	\left \{\H{ v} \in \hh^2( \H{\Omega}) \cap H^1_0( \H{\Omega}) : \nabla \H{v}|_{{\H{\Gamma}}} = 0 \right \} \to \left \{\H{v}\in \Sp^{p}(\Xi) \cap H^1_0( \H{\Omega}) : \nabla \H{v}|_{{\H{\Gamma}}} = 0 \right \}
\end{equation*}
such that
\begin{equation*}
	\norm{\nabla(\H{v} - \H{\Pi}^0_{k,h} \H{v})}_{L^2(\H{\Gamma})} + h\, \abs{\nabla(\H{v} - \H{\Pi}^0_{k,h} \H{v})}_{\hh^1(\H{\Gamma})} \leq C\,h^{\nicefrac{1}{2}} \norm{\H{v}}_{\hh^2(\H{\Omega})}.
\end{equation*}
Here $\Sp^{p}(\Xi)$ is the trace space of $\Sp^{\mathbf{p}}_k$ corresponding to the edge $\H{\Gamma}$.
\end{prop}
\begin{proof}
The operator $\H{\Pi}^0_{k,h} $ is as $\H{\Pi}_{k,h} $ in Proposition \ref{prop:projparamdom} where
only interior dual basis functions are active, see Figure~\ref{fig:lambdabasis}. { This is analogous to \citeRefDetail{Lemma 3.6}{bazilevs}}. 
\end{proof}
{ \subsection{Approximation properties on the physical domain}

We now define an isogeometric projector on a single patch, which inherits the approximation properties from the spline projector on $\H{\Omega}$ defined in Proposition~\ref{prop:projparamdom}.}

\begin{prop}\label{prop:projsplinec2}
Let $m$ and $t$ be integers such that $0 \leq t \leq m \leq p+1$.
Then, there exists an operator $\Pi_{k,h} \colon \left \{v \in H^4(\Omega_k) : v = \partial_n v = 0, \text{on } \partial \Omega \cap \partial \Omega_k \right \} \to X_{k,h}$ such that
\begin{equation*}
	\Pi_{k,h} v = v, \quad \forall v \in X_{k,h}, \\
\end{equation*}
with
\begin{equation}\label{eq:c2phys}
\begin{split}
	(\Pi_{k,h} v)(\mathcal{V}_\iota) 		&=  v(\mathcal{V}_\iota) \\
	\nabla (\Pi_{k,h} v)(\mathcal{V}_\iota)	&=\nabla v(\mathcal{V}_\iota)\\
	\nabla(\nabla (\Pi_{k,h} v))(\mathcal{V}_\iota)	&= \nabla(\nabla v)(\mathcal{V}_\iota),
\end{split}
\end{equation}
for all $\iota \in \mathcal{I}_{\mathcal{V},k}$, and,  for  $s = \max (4,m)$,
\begin{equation}\label{eq:physapp}
	\norm*{v - \Pi_{k,h} v}_{\hh^t(\Omega_k)} \leq C \, h^{m-t} \norm*{v}_{\hh^s(\Omega_k)}.
\end{equation}
For each edge $\Gamma$ of $\partial \Omega_k \setminus \partial \Omega$, the following inequality holds:
\begin{equation}\label{eq:physappbis}
\norm{\nabla(v - \Pi_{k,h} v)}_{L^2(\Gamma)} + h\, \abs{\nabla(v - \Pi_{k,h} v)}_{\hh^1(\Gamma)} \leq C\, h^{m - \nicefrac{3}{2}} \norm{v}_{\hh^s(\Omega_k)}.
\end{equation}
Additionally, for every boundary edge $\Gamma \subset \partial \Omega$, we have
\begin{equation*}
	\Pi_{k,h} v|_\Gamma = \dn \Pi_{k,h} v|_\Gamma = 0.
\end{equation*}
\end{prop}
\begin{proof}
Let $\Pi_{k,h} \colon \left \{v \in H^4(\Omega_k) : v = \partial_n v = 0, \text{on } \partial \Omega \cap \partial \Omega_k \right \} \to X_{k,h}$ such that  $\Pi_{k,h} v := (\H{\Pi}_{k,h}(v \circ \bm{F}_k)) \circ \bm{F}_k^{-1}$, where $\H{\Pi}_{k,h}$ is defined as in Proposition~\ref{prop:projparamdom}, where for each edge $\Gamma \subset \partial
\Omega$ the dual basis functions that correspond to inner basis functions that have non-vanishing gradient at the global boundary are removed. In Figure 1 those are the indices marked in red. If, e.g., the bottom edge is a boundary edge, the corresponding indices are $(i,2)$, with $i\in\{3,\ldots,n_1-2\}$. For $\mathcal{V}_\iota \in \partial \Omega$, the interpolation problem~\eqref{eq:const} is well-posed by construction, since~\eqref{eq:interp_system}.
Then, \eqref{eq:c2phys} follows from \eqref{eq:const} and the chain rule. Similarly, \eqref{eq:physapp}-\eqref{eq:physappbis} also follow from the chain rule, see \citeRefDetail{Section 3}{bazilevs} for further details.
\end{proof}
We can now define a projector onto $X_h$ by a patchwise construction.
\begin{prop}\label{prop:globalproj} 
Let $m$ be an integer such that $0 \leq m \leq p+1$, there exists a multi-patch operator $\Pi_h \colon H^4(\Omega) \cap H^2_0(\Omega) \to X_h$ such that for $s = \max (4,m)$
\begin{equation*}
	\norm{v - \Pi_h v}_X \leq C\,h^{m-2}\,\norm{v}_{H^s(\Omega)}.
\end{equation*}
\end{prop}
\begin{proof}
Let $v$ be in $H^4(\Omega) \cap H^2_0(\Omega)$. We define $\Pi_h$ as
\begin{equation*}
	(\Pi_h v)|_{\Omega_k} = \Pi_{k,h} (v|_{\Omega_k}), \quad \forall k \in \io.
\end{equation*}
The norm $\norm{\cdot}_X$ consists of patchwise $H^2$-norm contributions and norms of jumps of normal derivatives across all interfaces. The bounds for the $H^2$-norms follow from Proposition~\ref{prop:projsplinec2} with $t=2$. What remains is to prove that, for all $\ell \in \is$,
\begin{equation*}
	\norm{\jump{\dn v - \dn \Pi_h v}}_{H^{\nicefrac{1}{2}}(\gl)} \leq C\,h^{m-2}\,\norm{v}_{H^s\left( \Omega_{\ml} \cup \Omega_{\sell}\right) }.
\end{equation*}
From~\eqref{eq:physappbis}, it holds
\begin{align*}
	\norm{(\nabla v - \nabla \Pi_h v)|_{\Omega_{\sell}}}_{L^2(\gl)} &\leq C\,h^{m-\nicefrac{3}{2}}  \norm{v}_{H^s\left( \Omega_{\sell}\right) }, \\
	\norm{(\nabla v - \nabla \Pi_h v)|_{\Omega_{\sell}}}_{H^1(\gl)} &\leq C\,h^{m-\nicefrac{5}{2}} \norm{v}_{H^s\left( \Omega_{\sell}\right)}.
\end{align*}
Interpolating between $L^2(\gl)$ and $H^1(\gl)$, we get
\begin{equation*}
	\norm{(\nabla v - \nabla \Pi_h v)|_{\Omega_{\sell}}}_{H^{\nicefrac{1}{2}}(\gl)} \leq C\,h^{m-2} \abs{v}_{H^s\left( \Omega_{\sell}\right) }.
\end{equation*}
Similarly, we also obtain
\begin{equation*}
	\norm{(\nabla v - \nabla \Pi_h v)|_{\Omega_{\ml}}}_{H^{\nicefrac{1}{2}}(\gl)} \leq C\,h^{m-2} \abs{v}_{H^s\left( \Omega_{\ml}\right) }.
\end{equation*}
Finally, it holds
\begin{align*}
	\norm{\jump{\dn v - \dn \Pi_h v}}_{H^{\nicefrac{1}{2}}(\gl)} &\leq \norm{(\nabla v - \nabla \Pi_h v)|_{\Omega_{\sell}}}_{H^{\nicefrac{1}{2}}(\gl)}\\
	& \quad + \norm{(\nabla v - \nabla \Pi_h v)|_{\Omega_{\ml}}}_{H^{\nicefrac{1}{2}}(\gl)}\\
	& \leq  C\,h^{m-2}\,\norm{v}_{H^s\left( \Omega_{\ml} \cup \Omega_{\sell}\right), }
\end{align*}
which concludes the proof.
\end{proof}
	\section{Approximation properties of $V_h$}\label{sec:Vhprop}
The goal of this section is to establish the estimate
\begin{equation*}
	\inf_{v_h \in V_h} \norm{u - v_h}_X \leq C h^{p - 1} \norm{u}_{H^{s}(\Omega)},
\end{equation*}
where $s = \max(4, p+1)$.
We introduce the discrete trace derivative space $W_h$, defined as $W_h := \prod_{\ell \in \is} W_{\ell}$, where each component $W_{\ell}$ is given by
\begin{equation}\label{eq:defWl}
	W_{\ell} := \left\{w_h \in L^2(\gl) : \; w_h = \dn \varphi_h |_{\gl}, \text{ with }\quad
	\begin{aligned} 
		& \varphi_h \in X_{\sell,h} \cap H^1_0\left( \Omega_{\sell}\right), \\
		& \dn \varphi_h|_{\partial \Omega_{\sell} \setminus {\Gamma}_{\ell}} = 0
	\end{aligned}\right\}.
\end{equation}
Functions in $W_\ell$ and their derivatives vanish at the endpoints of the interface $\gl$. Furthermore, the degrees of freedom (control points) associated with $W_\ell$ are independent of those associated with $W_{\ell'}$ for $\ell' \neq \ell$. { Figure~\ref{fig:lambdabasis-b}} illustrates this distinction: degrees of freedom corresponding to the $C^2$ constraints at the vertices are marked in { light blue}, those related to the $C^0$ constraints along the interfaces are marked in { dark blue}, and the degrees of freedom of $\varphi_h$ spanning $W_\ell$, including $\dn \varphi_h$, are marked in { red. All remaining degrees of freedom are marked in black}.

We emphasize that the definition and properties of $V_h$ depend  on the choice of the multiplier space $M_h$. The construction of $M_h  = \prod_{\ell \in  \is} M_\ell$ must ensure that the following general assumptions are satisfied, which are essential for the well-posedness and approximation properties of the method.
\begin{assumption}\label{ass:oderM}
For every $\ell \in \is$ and any $\tau \in H^{p - \nicefrac{3}{2}}(\gl)$, the space $M_\ell$ satisfies
\begin{equation*}
	\inf_{\mu_h \in M_\ell} \norm{\tau - \mu_h}_{L^2(\gl)} \leq C h^{p - \nicefrac{3}{2}} \norm{\tau}_{H^{p - \nicefrac{3}{2}}(\gl)}.
\end{equation*}
\end{assumption}
\begin{assumption}\label{ass:dim_infsup}
For every $\ell \in \is$, the discrete trace space $W_\ell$ and the discrete multiplier space $M_\ell$ have the same dimension, i.e., $\dim(W_\ell) = \dim(M_\ell)$. Moreover, there exists a constant $\beta_\ell > 0$, independent of $h$, such that
\begin{equation*}
		\sup_{w_h \in W_\ell} \frac{\int_{\gl} w_h \, \mu_h \, \ds}{\norm{w_h}_{L^2(\gl)}} \geq \beta_\ell \norm{\mu_h}_{L^2(\gl)}, \quad \forall \mu_h \in M_\ell.
 \end{equation*}
\end{assumption}
\begin{prop}\label{prop:orderW-H12}
For each $w \in H_{00}^{\nicefrac{1}{2}}(\gl)$, the following inequality holds:
\begin{equation*}
	\inf_{w_h \in W_{\ell}} \norm{w - w_h}_{L^2(\gl)} \leq C \, h^{\nicefrac{1}{2}} \norm{w}_{H^{\nicefrac{1}{2}}(\gl)}.
\end{equation*}
\end{prop}

\begin{proof}
Let $w \in H_{00}^{\nicefrac{1}{2}}(\gl)$. Consider the solution $v$ to the boundary value problem
\begin{equation*}
\begin{cases}
	\Delta^2 v = 0 & \text{in } \Omega_{\sell}, \\
	v = 0 & \text{on } \partial \Omega_{\sell}, \\
	\dn v = w & \text{on } \Gamma_{\ell}, \\
	\dn v = 0 & \text{on } \partial \Omega_{\sell} \setminus \Gamma_{\ell}.
\end{cases}
\end{equation*}
Here, $\partial_{n}$ denotes the derivative in the direction of the { inward} unit normal vector to $\partial \Omega_{\sell}$.
According to \citeOneRef{grisvard2011elliptic}, this problem has a unique solution $v \in H^2(\Omega_{\sell}) \cap H^1_0(\Omega_{\sell})$, and the following estimate holds:
\begin{equation*}
	\norm{v}_{H^2(\Omega_{\sell})} \leq C \norm{w}_{H_{00}^{\nicefrac{1}{2}}(\gl)}.
\end{equation*}
Using Proposition \ref{prop:projparamdom_low-H2_regularity}, we can define a projector
\begin{align*}
	\Pi^{0}_{\sell,h}: & \left\{ v_h \in H^2(\Omega_{\sell}) \cap H^1_0(\Omega_{\sell}) \mid 
	\dn v_h|_{\partial \Omega_{\sell} \setminus \mathring{\Gamma}_{\ell}} = 0 \right\} \\
	& \to \left\{ v_h \in X_{\sell,h} \cap H^1_0(\Omega_{\sell}) \mid 
	\dn v_h|_{\partial \Omega_{\sell} \setminus \mathring{\Gamma}_{\ell}} = 0 \right\},
\end{align*}
such that
\begin{equation*}
	\norm{\nabla(v - \Pi^{0}_{\sell,h} v)}_{L^2(\gl)} \leq C h^{\nicefrac{1}{2}} \norm{v}_{H^2(\Omega_{\sell})}.
\end{equation*}
Next, we define the projector $P_{h,\ell}: H_{00}^{\nicefrac{1}{2}}(\gl) \to W_{\ell}$ by
\begin{equation*}
	P_{h,\ell}(w) := (\dn \Pi^{0}_{\sell,h} v)|_{\gl}.
\end{equation*}
Finally, we estimate the error between $w$ and $P_{h,\ell}w$:
\begin{align*}
	\norm{w - P_{h,\ell} w}_{L^2(\gl)} &\leq \norm{\nabla (v - \Pi^{0}_{\sell,h} v)}_{L^2(\gl)} \\
	&\leq C h^{\nicefrac{1}{2}} \norm{v}_{H^2(\Omega_{\sell})} \\
	&\leq C h^{\nicefrac{1}{2}} \norm{w}_{H^{\nicefrac{1}{2}}(\gl)},
\end{align*}
which completes the proof of the desired inequality.
\end{proof}
We define the mortar projection $\M_h$ by
\begin{equation*}
	(\M_h v)|_{\gl} := \M_{\ell}(v|_{\gl}), \quad \forall v \in L^2(\Sigma),
\end{equation*}
where, for each $\ell \in \is$, the local projection $\M_\ell \colon L^2(\gl) \to W_\ell$ is defined as
\begin{equation*}
	\int_{\gl} \M_\ell v \, \mu \, \ds = \int_{\gl} v \, \mu \, \ds, \quad \forall \mu \in M_{\ell}.
\end{equation*}
\begin{prop}\label{rem:bothstab}
Given Assumption~\ref{ass:dim_infsup},  for every $\ell \in \is$, it holds
\begin{equation*}
	\sup_{\mu_h \in M_\ell} \frac{\int_{\gl} w_h \, \mu_h \, \ds}{\norm{\mu_h}_{L^2(\gl)}} \geq \beta_\ell \norm{w_h}_{L^2(\gl)}, \quad \forall w_h \in W_\ell,
\end{equation*}
\end{prop}
\begin{proof}
It follows from \citeRefDetail{Proposition 3.4.3}{boffi}.
\end{proof}
\begin{prop}\label{cor:pistablel2}
For every $\ell \in \is$, the mortar projection $\M_\ell \colon L^2(\gl) \to W_\ell$ is well-defined and $L^2$-stable:
\begin{equation*}
	\norm{\M_\ell w}_{L^2(\gl)} \leq C \norm{w}_{L^2(\gl)}, \quad \forall w \in L^2(\gl).
\end{equation*}
\end{prop}
\begin{proof}
From Assumptions~\ref{ass:dim_infsup} and \citeRefDetail{Proposition 3.4.3}{boffi}, the matrix associated with the scalar product $(w, \mu)_{L^2}$ on $W_{\ell} \times M_{\ell}$ is square and invertible. Therefore, $\M_\ell$ is well-defined. Furthermore, from \citeRefDetail{Theorem 4.2.3}{boffi}, it follows that $\M_\ell$ is $L^2$-stable.
\end{proof}
\begin{lemma}\label{lemma:projproperties} 
For every $\ell \in \is$  the following inequalities hold:
\begin{equation*}
	\norm{w - \M_\ell \, w}_{L^2(\gl)} \leq C \,h^{\nicefrac{1}{2}}\norm{w}_{H^{\nicefrac{1}{2}}(\gl)}, \, \forall w \in H^{\nicefrac{1}{2}}_{00}(\gl)
\end{equation*}
and
\begin{equation*}
	\norm{\M_\ell w}_{H^{\nicefrac{1}{2}}(\gl)} \leq C \norm{w}_{H^{\nicefrac{1}{2}}(\gl)}, \quad \forall w \in H^{\nicefrac{1}{2}}_{00}(\gl).
\end{equation*}
\end{lemma}
\begin{proof}
Let $w \in H^1_{00}(\Gamma_\ell)$ and let $P_{h,\ell}$ be the spline approximation operator defined in the proof of Proposition~\ref{prop:orderW-H12}. It follows that
\begin{equation}\label{eq:mortar_est}
\begin{aligned}
	\norm{w - \M_\ell w}_{L^2(\Gamma_\ell)} &= \norm{w - P_{h,\ell} w - \M_\ell (w - P_{h,\ell} w)}_{L^2(\Gamma_\ell)} \\
	&\leq C \norm{w - P_{h,\ell} w}_{L^2(\Gamma_\ell)} \\
	&\leq C h^{1/2} \norm{w}_{H^{1/2}(\Gamma_\ell)}.
\end{aligned}
\end{equation}
Using  the inverse inequality for splines and~\eqref{eq:mortar_est} we obtain
\begin{align*}
	\norm{\M_\ell w}_{H^{1/2}(\Gamma_\ell)} &\leq \norm{\M_\ell (w - P_{h,\ell} w)}_{H^{1/2}(\Gamma_\ell)} + \norm{ P_{h,\ell} w}_{H^{1/2}(\Gamma_\ell)}\\
	&\leq C h^{-1/2} \norm{\M_\ell (w - P_{h,\ell} w)}_{L^{2}(\Gamma_\ell)} + \norm{ w}_{H^{1/2}(\Gamma_\ell)} \\
	& \leq C \norm{ w}_{H^{1/2}(\Gamma_\ell)}.
\end{align*}
\end{proof}
In addition to the local projector $\M_\ell$ for each $\ell \in \is$, we define the dual operator $\M_{\ell}^{*} \colon L^2(\gl) \to M_\ell$ by the following condition:
\begin{equation}\label{eq:mortarprojdual}
	\int_{\gl} \M_{\ell}^{*} \mu \, v \, \ds = \int_{\gl} \mu \, v \, \ds, \quad \forall v \in W_{\ell}.
\end{equation}
\begin{lemma}\label{lemma:dualproperties}
For every $\ell \in \is$, the dual operator $\M_{\ell}^{*} \colon L^2(\gl) \to M_\ell$ is well-defined and $L^2$-stable:
\begin{equation*}
	\norm{\M_{\ell}^{*} \mu}_{L^2(\gl)} \leq C \norm{\mu}_{L^2(\gl)}, \quad \forall \mu \in L^2(\gl).
\end{equation*}
Furthermore, for all $\mu \in H^{p - \nicefrac{3}{2}}(\gl)$, the following inequality holds:
\begin{equation*}
	\norm{\mu - \M_{\ell}^{*} \mu}_{L^2(\gl)} \leq C\,h^{p - \nicefrac{3}{2}}\, \norm{\mu}_{H^{p - \nicefrac{3}{2}}(\gl)}.
\end{equation*}
\end{lemma}

\begin{proof}
As in the proof of Proposition~\ref{cor:pistablel2}, Assumption~\ref{ass:dim_infsup} guarantees that $\M_{\ell}^{*}$ is well-defined. By combining Assumption~\ref{ass:dim_infsup} with~\eqref{eq:mortarprojdual}, we obtain the following estimate:
\begin{equation*}
\begin{split}
	\norm{\M_{\ell}^{*} \mu}_{L^2(\gl)} 
	&\leq \frac{1}{\beta_\ell} \sup_{w_h \in W_{\ell}} \frac{\int_{\gl} \M_{\ell}^{*} \mu \, w_h \, \ds}{\norm{w_h}_{L^2(\gl)}} \\
	&= \frac{1}{\beta_\ell} \sup_{w_h \in W_{\ell}} \frac{\int_{\gl} \mu \, w_h \, \ds}{\norm{w_h}_{L^2(\gl)}} \\
	&\leq \frac{1}{\beta_\ell} \sup_{w_h \in W_{\ell}} \frac{\norm{\mu}_{L^2(\gl)} \norm{w_h}_{L^2(\gl)}}{\norm{w_h}_{L^2(\gl)}} \\
	&= \frac{1}{\beta_\ell} \norm{\mu}_{L^2(\gl)}.
\end{split}
\end{equation*}
The approximation property of $\M_{\ell}^{*}$ follows from Assumption~\ref{ass:oderM}. Indeed, for $\mu \in H^{p - \nicefrac{3}{2}}(\gl)$, let $\mu_h \in M_\ell$ be such that
\begin{equation*}
	\norm{\mu - \mu_h}_{L^2(\gl)} \leq C h^{p - \nicefrac{3}{2}} \norm{\mu}_{H^{p - \nicefrac{3}{2}}(\gl)},
\end{equation*}
then, we have the following estimate:
\begin{equation*}
\begin{split}
	\norm{\mu - \M_{\ell}^{*} \mu}_{L^2(\gl)} 
	&= \norm{\mu - \mu_h - \M_{\ell}^{*}(\mu - \mu_h)}_{L^2(\gl)} \\
	&\leq C' \norm{\mu - \mu_h}_{L^2(\gl)} \\
	&\leq C h^{p - \nicefrac{3}{2}} \norm{\mu}_{H^{p - \nicefrac{3}{2}}(\gl)}.
\end{split}
\end{equation*}
This completes the proof.
\end{proof}
We now proceed with the study of the properties of $V_h$.
\begin{lemma}\label{lemma:eqbrokennorm2}
For every $v_h \in V_h$, the following inequality holds:
\begin{equation*}
	\sum_{k \in \io} \norm{v_h}^2_{H^2(\Omega_k)} \leq C \sum_{k \in \io} \abs{v_h}^2_{H^2(\Omega_k)}.
\end{equation*}
\end{lemma}
\begin{proof}
The proof follows the approach outlined in~\citeRefDetail{Theorem 1}{marcinkowski}. We first show that for all $v_h \in V_h$, it holds:
\begin{equation}\label{eq:est1}
	\sum_{k \in \io} \norm{\partial_x v_h}^2_{L^2(\Omega_k)} \leq C \sum_{k \in \io} \abs{v_h}^2_{H^2(\Omega_k)}.
\end{equation}
Let ${ \mathbf{x} =} (x,y) \in \Omega$ be an arbitrary point. Consider the endpoints $(a,y)$ and $(b,y) \in \partial \Omega$ of a segment containing $(x,y)$, and define the set
\begin{equation*}
	I[a,x] := \left\{ t \in [a,x] \mid (t,y) \in \Gamma \right\},
\end{equation*}
where $N_I = |I[a,b]|$. Then, we have
\begin{equation}\label{eq:est2}
	\partial_x v_h(x,y) = \int_a^x { \partial_{x}\partial_{x}} v_h(t,y) \, \mathrm{d} t + \sum_{t_{\ell} \in I[a,x]} \jump{\partial_x v_h}_{\gl}(t_{\ell},y).
\end{equation}
The first term on the right-hand side of~\eqref{eq:est2} can be bounded as follows:
\begin{equation*}
	\left| \int_a^x { \partial_{x}\partial_{x}} v_h(t,y) \, \mathrm{d} t \right| \leq \int_a^b \left| { \partial_{x}\partial_{x}} v_h(t,y) \right| \, \mathrm{d} t \leq M^{\frac{1}{2}} \left( \int_a^b \left| { \partial_{x}\partial_{x}} v_h(t,y) \right|^2 \, \mathrm{d} t \right)^{\frac{1}{2}},
\end{equation*}
where $M = \text{diam}(\Omega)$. The second term is bounded by
\begin{equation*}
\begin{split}
	\left| \sum_{t_{\ell} \in I[a,x]} \jump{\partial_x v_h}_{\gl}(t_{\ell},y) \right| &\leq \sum_{t_{\ell} \in I[a,b]} \left| \jump{\partial_x v_h}_{\gl}(t_{\ell},y) \right| \\
	&\leq \left( N_I \right)^{\frac{1}{2}} \left( \sum_{t_{\ell} \in I[a,b]} \left| \jump{\partial_x v_h}(t_{\ell},y) \right|^2 \right)^{\frac{1}{2}}.
\end{split}
\end{equation*}
Therefore, we obtain the following estimate:
\begin{equation}\label{eq:aux01}
	\left| \partial_x v_h(x,y) \right|^2 \leq C \int_a^b \left| { \partial_{x}\partial_{x}} v_h(t,y) \right|^2 \, \mathrm{d} t + C \sum_{t_{\ell} \in I[a,b]} \left| \jump{\partial_x v_h}_{\gl}(t_{\ell},y) \right|^2.
\end{equation}
Without loss of generality, assume that $\Omega \subset [0, M]^2$, and extend $\left| { \partial_{x}\partial_{x}} v_h(t,y) \right|^2$ to zero outside of $\Omega$. Integrating the first term on the right-hand side of~\eqref{eq:aux01} over $\Omega$, we get
\begin{equation*}
\begin{split}
	\int_\Omega \left( \int_a^b \left| { \partial_{x}\partial_{x}} v_h(t,y) \right|^2 \, \mathrm{d} t \right) \, {\mathrm{d} \mathbf{x}}
	&\leq \int_0^M \int_0^M \int_0^M \left| { \partial_{x}\partial_{x}} v_h(t,y) \right|^2 \, \mathrm{d} t \, \mathrm{d} x \, \mathrm{d} y \\
	&= M \int_0^M \int_0^M \left| { \partial_{x}\partial_{x}} v_h(t,y) \right|^2 \, \mathrm{d} t \, \mathrm{d} y \\
	&= M \int_\Omega \left| { \partial_{x}\partial_{x}} v_h(t,y) \right|^2 \, {\mathrm{d} \mathbf{x}}.
\end{split}
\end{equation*}
Similarly, integrating the second term of~\eqref{eq:aux01} over $\Omega$, we obtain
\begin{equation*}
	\int_\Omega \sum_{t_{\ell} \in I[a,b]} \left| \jump{\partial_x v_h}_{\gl}(t_{\ell},y) \right|^2 \, {\mathrm{d}\mathbf{x}} \leq C \sum_{\ell \in \is} \int_{\gl} \left| \jump{\dn v_h}_{\gl} \right|^2 \, {\ds}.
\end{equation*}
Therefore, we have the following inequality:
\begin{equation*}
	\int_\Omega \left| \partial_x v_h(x,y) \right|^2 \, {\mathrm{d}\mathbf{x}} \leq C \int_\Omega \left| { \partial_{x}\partial_{x}} v_h \right|^2 \, \mathrm{d}\mathbf{x} + C \sum_{\ell \in \is} \int_{\gl} \left| \jump{\dn v_h}_{\gl} \right|^2 \, {\ds}.
\end{equation*}
To estimate the last integral, note that by the mortar condition~\eqref{eq:mortarkernel}, the average value of $\jump{\dn v_h}_{\gl}$ over $\gl$ is zero. Applying the standard trace theorem and the Poincaré inequality, we obtain
\begin{equation*}
	\int_{\gl} \left| \jump{\dn v_h}_{\gl} \right|^2 \, {\ds} \leq C_{\ml} \abs{v_h}^2_{H^2(\Omega_{\ml})} + C_{\sell} \abs{v_h}^2_{H^2(\Omega_{\sell})}.
\end{equation*}
Summing over all interfaces, we obtain~\eqref{eq:est1}.\\
Using similar arguments for the $y$-derivative, we can find a constant $C > 0$ such that for every $v_h \in V_h$, the following inequality holds:
\begin{equation}\label{eq:est_h1}
	\sum_{k \in \io} \abs{v_h}^2_{H^1(\Omega_k)} \leq C \sum_{k \in \io} \abs{v_h}^2_{H^2(\Omega_k)}.
\end{equation}
Moreover, since $V_h \subset H^1_0(\Omega)$, applying the standard Poincaré inequality yields
\begin{equation}\label{eq:est_l2}
	\sum_{k \in \io} \norm{v_h}^2_{L^2(\Omega_k)} \leq C \sum_{k \in \io} \abs{v_h}^2_{H^1(\Omega_k)}.
\end{equation}
Finally, the thesis follows by combining~\eqref{eq:est_h1} and~\eqref{eq:est_l2}.
\end{proof}
\begin{lemma}\label{lemma:es1}
Let $s = \max(4, p+1)$. For all $v \in H^s(\Omega)$, 
\begin{equation*}
	\inf_{v_h \in V_h} \|v - v_h\|_X \leq C h^{p - 1} \norm{v}_{H^s(\Omega)}.
\end{equation*}
\end{lemma}
\begin{proof}
We define the extension operator $\hh_{h,\ell}: W_\ell \to X_{\sell,h}$ as follows.  
Given $w_h \in W_\ell$, let ${\varphi_h \in X_{\sell,h} \cap H^1_0(\Omega_{\sell})}$ be as in~\eqref{eq:defWl}.  
We then set
\begin{equation*}
\hh_{h,\ell} w_h := \varphi_h + \psi_h,
\end{equation*}
where $\psi_h \in X_{\sell,h} \cap H^2_0(\Omega_{\sell})$ is the unique function satisfying
\begin{equation*}
	\int_{\Omega_{\sell}} \Delta (\psi_h + \varphi_h)\, \Delta z_h \, { \mathrm{d}\mathbf{x}} = 0 
	\quad \forall z_h \in X_{\sell,h} \cap H^2_0(\Omega_{\sell}).
\end{equation*}
It follows that $\partial_{n} \hh_{h,\ell} w_h |_{\Gamma_{\ell}} = w_h$, { $\partial_{n} \hh_{h,\ell} w_h |_{\partial \Omega_{\sell} \setminus {\Gamma}_{\ell}} = 0$}, and
\begin{equation*}
	\|\hh_{h,\ell} w_h\|_{H^2(\Omega_{\sell})} \leq C \norm{w_h}_{H^{1/2}(\partial \Omega_{\sell})} \leq C \norm{w_h}_{H^{1/2}(\Gamma_\ell)}.
\end{equation*}
We further extend $\hh_{h,\ell} {w_h}$ to be zero outside $\Omega_{\sell}$.\\
Now, let $v \in H^{p+1}(\Omega)$, and let $\Pi_h v \in X_h$ be the global spline projection of $v$ in $X_h$, as defined in Proposition~\ref{prop:globalproj}. We define an element $v_h \in V_h$ by
\begin{equation*}
	v_h := \Pi_h v - \sum_{\ell \in \is} \hh_{h,\ell} \M_\ell \jump{\dnl \Pi_h v}_{\Gamma_\ell}.
\end{equation*}
From the definition of $\M_\ell$, we have that for all $\tau_h \in M_h$,
\begin{equation*}
	\int_{\Gamma_\ell} \jump{\dnl v_h}_{\Gamma_\ell} \, \tau_h \, \ds = \int_{\Gamma_\ell} \left( \jump{\dnl \Pi_h v}_{\Gamma_\ell} - \M_\ell \jump{\dnl \Pi_h v} \right) \tau_h \, \ds = 0.
\end{equation*}
Thus, $v_h \in V_h$. Consequently, we can write
\begin{equation*}
	\|v - v_h\|_X \leq \|v - \Pi_h v\|_X + \left\|\sum_{\ell \in \is} \hh_{h,\ell} \M_\ell \jump{\dnl \Pi_h v}_{\Gamma_\ell}\right\|_X.
\end{equation*}
The first term is estimated in Proposition~\ref{prop:globalproj}.\\
On the other hand, using the definition of $\hh_{h,\ell}$, we obtain the following estimate:
\begin{equation*}
\begin{split}
	\left\|\sum_{\ell \in \is} \hh_{h,\ell} \M_\ell \jump{\dnl \Pi_h v}_{\Gamma_\ell}\right\|_X^2 &\leq \sum_{\ell \in \is} \|\hh_{h,\ell} \M_\ell \jump{\dnl \Pi_h v}_{\Gamma_\ell}\|_{H^2\left(\Omega_{\sell}\right)}^2 \\
	& \quad+ \sum_{\ell \in \is} \|\M_\ell \jump{\dnl \Pi_h v}_{\Gamma_\ell}\|_{H^{1/2}(\Gamma_\ell)}^2 \\
	&\leq C \sum_{\ell \in \is} \|\M_\ell \jump{\dnl \Pi_h v}_{\Gamma_\ell}\|_{H^{1/2}(\Gamma_\ell)}^2.
\end{split}
\end{equation*}
The $H^{1/2}_{00}$-stability of the mortar projector $\M_\ell$, stated in Lemma~\ref{lemma:projproperties}, and Proposition~\ref{prop:globalproj} yield
\begin{equation*}
\begin{split}
	\|\M_\ell \jump{\dnl \Pi_h v}_{\Gamma_\ell}\|_{H^{1/2}(\Gamma_\ell)}^2 &\leq \|\jump{\dnl \Pi_h v}_{\Gamma_\ell}\|_{H^{1/2}(\Gamma_\ell)}^2 \\
	&\leq C h^{2(p-1)} \left(\norm{v}_{H^s\left(\Omega_{\ml}\right)}^2 + \norm{v}_{H^s\left(\Omega_{\sell}\right)}^2 \right).
\end{split}
\end{equation*}
Finally, by Assumption~\ref{ass:quasiunif}, summing over the interfaces, we obtain the desired result.
\end{proof}
	\section{Error analysis for the model problem}\label{sec:error_analysis}
Since $V_h \not\subseteq H^2(\Omega)$ (as the space $V_h$ is only weakly $C^1$), we must employ Strang's lemma to bound the numerical error in our formulation.
\begin{lemma}\label{lemma:strang}
Let $u$ be the solution of Problem~\eqref{eq:mortar-variational-formulation}, and $u_h$ the solution of the discrete Problem~\eqref{eq:qh}. It holds
\begin{equation}\label{eq:strang}
	\norm{u - u_h}_X \leq C \left( \inf_{v_h \in V_h} \norm{u - v_h}_X + \sup_{v_h \in V_h} \frac{\abs*{\int_{\Sigma} { \dn \dn }u  \jump{\dn v_h} \, \ds}}{\norm{v_h}_X}\right).
\end{equation}
\end{lemma}
\begin{proof}
From Lemma~\ref{lemma:eqbrokennorm2}, the bilinear form $a(\cdot, \cdot)$ is coercive on $V_h$ with respect to the broken norm $\sum_{k \in \io}\norm{\cdot}^2_{H^2(\Omega_k)}$. Our goal is to establish coercivity  with respect to $\norm{\cdot}_X$. 

We recall the inverse inequality:
\begin{equation*}
	\norm{\jump{\dn v_h}_{\gl}}_{H^{\nicefrac{1}{2}}(\gl)}^2 \leq C h^{-1} \norm{\jump{\dn v_h}_{\gl}}_{L^2(\gl)}^2  ,
\end{equation*}
where this equality follows from the fact that the jump $\jump{\dn v_h}_{\gl}$ is orthogonal to $M_\ell$ for all $v_h \in V_h$ and for every $\ell \in \is$. Using the approximation properties of $\M_{\ell}^{*}$ stated in Lemma~\ref{lemma:dualproperties}, we derive:
\begin{equation*}
\begin{split}
	\norm{\jump{\dn v_h}_{\gl}}_{H^{\nicefrac{1}{2}}(\gl)}^2 &\leq C h^{-1} \norm{\jump{\dn v_h}_{\gl}}_{L^2(\gl)}^2 \leq C h^{-1} \norm{(\id - \M_{\ell}^{*})\jump{\dn v_h}_{\gl}}_{L^2(\gl)}^2 \\
	&\leq C h^{-1} \left( \norm{(\id - \M_{\ell}^{*}) \dn v_h|_{\Omega_{\ml}}}_{L^2(\gl)}^2 \right. \\
	&\left. \quad+ \norm{(\id - \M_{\ell}^{*})\dn v_h|_{\Omega_{\sell}}}_{L^2(\gl)}^2 \right) \\
	&\leq C h^{-1} \left( h \norm{\dn v_h|_{\Omega_{\ml}}}_{H^{\nicefrac{1}{2}}(\gl)}^2 + h \norm{\dn v_h|_{\Omega_{\sell}}}_{H^{\nicefrac{1}{2}}(\gl)}^2 \right) \\
	&\leq C \left(\norm{v_h}_{H^2\left( \Omega_{\ml}\right) }^2 + \norm{v_h}_{H^2\left( \Omega_{\sell}\right) }^2 \right).
\end{split}
\end{equation*}
Thus, $a(\cdot, \cdot)$ is coercive with respect to $\norm{\cdot}_X$, allowing us to apply the second Strang Lemma (see, e.g., \cite[Theorem 4.2.2]{ciarlet}), which gives~\eqref{eq:strang}.
\end{proof}
The first term in~\eqref{eq:strang} represents the best approximation error in $V_h$, which was analyzed in Lemma~\ref{lemma:es1}. The second term corresponds to the consistency error. The following result establishes an optimal estimate for the consistency error.
\begin{lemma}\label{lemma:es2}
The following inequality holds
\begin{equation*}
	\sup_{v_h \in V_h} \frac{\abs*{\int_{\Sigma} { \dn\dn} u \jump{\partial_n v_h} \, \ds}}{\norm{v_h}_X} \leq C \sum_{k \in \io} h^{p - 1} \norm{u}_{H^{p+1}(\Omega_k)}.
\end{equation*}
\end{lemma}
\begin{proof}
Let $v_h \in V_h$. By applying Lemma~\ref{lemma:dualproperties}, we have
\begin{equation*}
\begin{split}
	\int_{\Sigma} { \dn\dn} u \jump{\partial_n v_h} \, \ds &= \int_{\Sigma} ({ \dn\dn} u - \M_{\ell}^{*}{ \dn\dn} u) \jump{\partial_n v_h} \, \ds \\
	&\leq \sum_{\ell \in \is} \norm{{ \dn\dn} u - \M_{\ell}^{*}{ \dn\dn} u}_{L^2(\gl)} \norm{\jump{\dn v_h}_{\gl}}_{L^2(\gl)} \\
	&\leq C h^{p-\nicefrac{3}{2}} \sum_{\ell \in \is} \norm{{ \dn\dn} u}_{H^{p-\nicefrac{3}{2}}(\gl)} \norm{\jump{\dn v_h}_{\gl}}_{L^2(\gl)}.
\end{split}
\end{equation*}
By Lemma~\ref{lemma:projproperties}, we further obtain
\begin{equation*}
\begin{split}
	\norm{\jump{\dn v_h}_{\gl}}_{L^2(\gl)} &\leq \norm{(\id - \M_{\ell})\jump{\dn v_h}_{\gl}}_{L^2(\gl)} \\
	&\leq C h^{\nicefrac{1}{2}} \norm{\jump{\dn v_h}_{\gl}}_{H^{\nicefrac{1}{2}}(\gl)} \leq C h^{\nicefrac{1}{2}} \norm{v_h}_X.
\end{split}
\end{equation*}
Using the standard trace theorem, we deduce that
\begin{equation*}
	\norm{{ \dn\dn} u}_{H^{p-\nicefrac{3}{2}}(\gl)} \leq C \norm{u}_{H^{p+\nicefrac{1}{2}}(\gl)} \leq C \norm{u}_{H^{p+1}(\Omega_{\sell})}.
\end{equation*}
Consequently, we arrive at the following bound:
\begin{equation*}
	\int_{\Sigma} \partial_{n}u  \jump{\partial_n v_h} \,\ds \leq C h^{p-1} \sum_{k \in \io} \norm{u}_{H^{p+1}(\Omega_{k})} \norm{v_h}_X.
\end{equation*}
This completes the proof.
\end{proof}
We combine the results from Lemmas~\ref{lemma:es1},~\ref{lemma:strang}, and~\ref{lemma:es2} to derive the optimal a priori error estimate.
\begin{theorem}\label{thm:orderofconv}
Let $s = \max(4, p+1)$. For $u \in H^{s}(\Omega)$, the following inequality holds:
\begin{equation*}
	\| u - u_h \|_X \leq C\, h^{p-1} \| u \|_{H^{s}(\Omega)}.
\end{equation*}
\end{theorem}
	\section{Discrete {inf-sup} condition}\label{sec:dscinfsup}
This section focuses on defining the discrete multiplier spaces $M_\ell$, for $\ell \in \is$, that satisfy Assumption~\ref{ass:dim_infsup}. For each interface, we construct a spline space of degree reduced by two, with increased smoothness near the vertices. We begin by considering a single interface $\Gamma$ and denote the knot vector on the reference interface $\H{\Gamma} := (0,1)$ as $\Xi$, where
\begin{equation*}
	\Xi = (\xi_1, \xi_2, \ldots, \xi_{n+p+1}),
\end{equation*}
with $\xi_1 = \dots = \xi_{p+1} = 0$ and $\xi_{n+1} = \dots = \xi_{n+p+1} = 1$.\\
The key step is to establish the inf-sup condition on the parameter domain $\H{\Gamma}$. On this domain, the trace space of the normal derivative is denoted by $\mathcal{S}_{0,\partial}^p$, which consists of splines of degree $p$ with vanishing value and derivative at the interface boundaries. We define the corresponding Lagrange multiplier space  as  $\mathcal{S}_{\textsc{m}}^{p-2}$, that is the spline space of degree $p-2$ with increased smoothness between the first two and last two elements, respectively (see~\eqref{eq:spaces} below). 

As a reminder, $\Sp^p(\Xi)$ is the spline space on $\H{\Gamma}$ generated by the knot vector $\Xi$. We introduce the modified knot vector
\begin{equation*}
 \Xi^p_{\textsc{m}} = (\xi_1,\ldots,\xi_{p+1},\xi_{p+3},\ldots,\xi_{n-1},\xi_{n+1},\ldots,\xi_{n+p+1})
\end{equation*}
by removing the first and last inner knots, $\xi_{p+2}$ and $\xi_{n}$, respectively. We define the degree-reduced knot vectors $\Xi^{p-1}$ and $\Xi^{p-2}$, where the boundary knots have multiplicity  $p$ and $p-1$, respectively, instead of $p+1$. Similarly, we define the knot vectors $\Xi_{\textsc{m}}^{p-1}$ and $\Xi_{\textsc{m}}^{p-2}$.
If the first and last inner knots have multiplicity $1$, then $\Sp^p(\Xi^p_{\textsc{m}})$ is the spline space of degree $p$ on $\H{\Gamma}$ with the first two and last two elements merged, see Figure~\ref{fig:basisSp0SpM} for an example of $\Sp^p_{0}$ and the corresponding multiplier space $\Sp^{p-2}_{\textsc{m}}$.
\begin{figure}
\centering
\begin{tikzpicture}
\begin{axis}[xtick={0.001,0.125,0.25,0.375,0.5,0.625,0.75,0.875,1}, 
			 xticklabels={$0$, $h$, $2h$, $3h$, $4h$, $5h$, $6h$, $7h$, $1$}, 
			 axis x line=center, axis y line=center,
			 xmin = 0, xmax = 1, enlargelimits=0.05, 
			 width=0.75\linewidth,
			 height=0.25\linewidth]
			 
\addplot[domain=0:0.125, 
mesh, patch type=cubic spline, patch type sampling,  thick, black]
{-469.3333*x^3 +96.0000*x^2};
\addplot[domain=0.125:0.25,
mesh, patch type=cubic spline, patch type sampling,  thick, black]
{298.6667*(x-0.125)^3 -80.0000*(x-0.125)^2 +2.0000*(x-0.125) +0.5833};
\addplot[domain=0.25:0.375,
mesh, patch type=cubic spline, patch type sampling,  thick, black]
{-85.3333*(x-0.25)^3 +32.0000*(x-0.25)^2 -4.0000*(x-0.25) +0.1667};		 
			
\addplot[domain=0:0.125, 
mesh, patch type=cubic spline, patch type sampling,  thick, blue]
{85.3333*x^3};
\addplot[domain=0.125:0.25,
mesh, patch type=cubic spline, patch type sampling,  thick, blue]
{-256.0000*(x - 0.125)^3 +32.0000*(x - 0.125)^2 +4.0000*(x - 0.125) +0.1667};
\addplot[domain=0.25:0.375,
mesh, patch type=cubic spline, patch type sampling, thick, blue]
{256.0000*(x - 0.25)^3 -64.0000*(x - 0.25)^2 +0.6667};
\addplot[domain=0.375:0.5,
mesh, patch type=cubic spline, patch type sampling,  thick, blue]
{-85.3333*(x - 0.375)^3 +32.0000*(x - 0.375)^2 -4.0000*(x - 0.375) +0.1667};

\addplot[domain=0.125:0.25,
mesh, patch type=cubic spline, patch type sampling,  thick, red]
{85.3333*(x - 0.125)^3};
\addplot[domain=0.25:0.375,
mesh, patch type=cubic spline, patch type sampling,  thick, red]
{-256.0000*(x - 0.25)^3 +32.0000*(x - 0.25)^2 +4.0000*(x - 0.25) +0.1667};
\addplot[domain=0.375:0.5,
mesh, patch type=cubic spline, patch type sampling,  thick, red]
{256.0000*(x - 0.375)^3 -64.0000*(x - 0.375)^2 +0.6667};
\addplot[domain=0.5:0.625,
mesh, patch type=cubic spline, patch type sampling,  thick, red]
{-85.3333*(x - 0.5)^3 +32.0000*(x - 0.5)^2 -4.0000*(x - 0.5) +0.1667};

\addplot[domain=0.25:0.375,
mesh, patch type=cubic spline, patch type sampling,  thick, green]
{85.3333*(x - 0.25)^3};
\addplot[domain=0.375:0.5,
mesh, patch type=cubic spline, patch type sampling,  thick, green]
{-256.0000*(x - 0.375)^3 +32.0000*(x - 0.375)^2 +4.0000*(x - 0.375) +0.1667};
\addplot[domain=0.5:0.625,
mesh, patch type=cubic spline, patch type sampling,  thick, green]
{256.0000*(x - 0.5)^3 -64.0000*(x - 0.5)^2 +0.6667};
\addplot[domain=0.625:0.75,
mesh, patch type=cubic spline, patch type sampling,  thick, green]
{-85.3333*(x - 0.625)^3 +32.0000*(x - 0.625)^2 -4.0000*(x - 0.625) +0.1667};

\addplot[domain=0.375:0.5,
mesh, patch type=cubic spline, patch type sampling,  thick, orange]
{85.3333*(x - 0.375)^3};
\addplot[domain=0.5:0.625,
mesh, patch type=cubic spline, patch type sampling,  thick, orange]
{-256.0000*(x - 0.5)^3 +32.0000*(x - 0.5)^2 +4.0000*(x - 0.5) +0.1667};
\addplot[domain=0.625:0.75,
mesh, patch type=cubic spline, patch type sampling,  thick, orange]
{256.0000*(x - 0.625)^3 -64.0000*(x - 0.625)^2 +0.6667};
\addplot[domain=0.75:0.875,
mesh, patch type=cubic spline, patch type sampling,  thick, orange]
{-85.3333*(x - 0.75)^3 +32.0000*(x - 0.75)^2 -4.0000*(x - 0.75) +0.1667};

\addplot[domain=0.375:0.5,
mesh, patch type=cubic spline, patch type sampling,  thick, orange]
{85.3333*(x - 0.375)^3};
\addplot[domain=0.5:0.625,
mesh, patch type=cubic spline, patch type sampling,  thick, orange]
{-256.0000*(x - 0.5)^3 +32.0000*(x - 0.5)^2 +4.0000*(x - 0.5) +0.1667};
\addplot[domain=0.625:0.75,
mesh, patch type=cubic spline, patch type sampling,  thick, orange]
{256.0000*(x - 0.625)^3 -64.0000*(x - 0.625)^2 +0.6667};
\addplot[domain=0.75:0.875,
mesh, patch type=cubic spline, patch type sampling,  thick, orange]
{-85.3333*(x - 0.75)^3 +32.0000*(x - 0.75)^2 -4.0000*(x - 0.75) +0.1667};

\addplot[domain=0.5:0.625,
mesh, patch type=cubic spline, patch type sampling,  thick, magenta]
{85.3333*(x - 0.5)^3};
\addplot[domain=0.625:0.75,
mesh, patch type=cubic spline, patch type sampling,  thick, magenta]
{-256.0000*(x - 0.625)^3 +32.0000*(x - 0.625)^2 +4.0000*(x - 0.625) +0.1667};
\addplot[domain=0.75:0.875,
mesh, patch type=cubic spline, patch type sampling,  thick, magenta]
{256.0000*(x - 0.75)^3 -64.0000*(x - 0.75)^2 +0.6667};
\addplot[domain=0.875:1,
mesh, patch type=cubic spline, patch type sampling,  thick, magenta]
{-85.3333*(x - 0.875)^3 +32.0000*(x - 0.875)^2 -4.0000*(x - 0.875) +0.1667};

\addplot[domain=0.625:0.75,
mesh, patch type=cubic spline, patch type sampling,  thick, teal]
{85.3333*(x - 0.625)^3};
\addplot[domain=0.75:0.875,
mesh, patch type=cubic spline, patch type sampling,  thick, teal]
{-298.6667*(x - 0.75)^3 +32.0000*(x - 0.75)^2 +4.0000*(x - 0.75) +0.1667};
\addplot[domain=0.875:1,
mesh, patch type=cubic spline, patch type sampling,  thick, teal]
{469.3333*(x - 0.875)^3 -80.0000*(x - 0.875)^2 -2.0000*(x - 0.875) +0.5833};

\end{axis}
\end{tikzpicture}

\begin{tikzpicture}
	\begin{axis}[xtick={0.001,0.25,0.375,0.5,0.625,0.75,1}, 
		xticklabels={$0$, $2h$, $3h$, $4h$, $5h$, $6h$, $1$}, 
		axis x line=center, axis y line=center,
		xmin = 0, xmax = 1, enlargelimits=0.05, width=0.75\linewidth,
		height=0.25\linewidth]
		
		\addplot[domain=0:0.25, 
		mesh, patch type=cubic spline, patch type sampling,  thick, black]
		{-4*x + 1};	 
		
		\addplot[domain=0:0.25, 
		mesh, patch type=cubic spline, patch type sampling,  thick, blue]
		{4*x};
		\addplot[domain=0.25:0.375,
		mesh, patch type=cubic spline, patch type sampling,  thick, blue]
		{-8*(x - 0.25) + 1};
		
		\addplot[domain=0.25:0.375,
		mesh, patch type=cubic spline, patch type sampling,  thick, red]
		{8*(x - 0.25)};
		\addplot[domain=0.375:0.5,
		mesh, patch type=cubic spline, patch type sampling,  thick, red]
		{-8*(x - 0.375) + 1};
		
		\addplot[domain=0.375:0.5,
		mesh, patch type=cubic spline, patch type sampling,  thick, green]
		{8*(x - 0.375)};
		\addplot[domain=0.5:0.625,
		mesh, patch type=cubic spline, patch type sampling,  thick, green]
		{-8*(x - 0.5) + 1};
		
		\addplot[domain=0.5:0.625,
		mesh, patch type=cubic spline, patch type sampling,  thick, orange]
		{8*(x - 0.5)};
		\addplot[domain=0.625:0.75,
		mesh, patch type=cubic spline, patch type sampling,  thick, orange]
		{-8*(x - 0.625) + 1};
		
		\addplot[domain=0.625:0.75,
		mesh, patch type=cubic spline, patch type sampling,  thick, magenta]
		{8*(x - 0.625)};
		\addplot[domain=0.75:1,
		mesh, patch type=cubic spline, patch type sampling,  thick, magenta]
		{-4*(x - 0.75) + 1};
		
		\addplot[domain=0.75:1,
		mesh, patch type=cubic spline, patch type sampling,  thick, teal]
		{4*(x - 0.75)};
		
	\end{axis}
\end{tikzpicture}

\caption{Top: basis functions of the space $\Sp^3_{0,\partial}$ on the knot vector $\Xi = \{0, 0, h, \dots, 7h, 1, 1\}$. Bottom: basis functions of the space $\Sp^1_{\textsc{m}}$ on the knot vector $\Xi_{\textsc{m}} = \{0, 0, 2h, \dots, 6h, 1, 1\}$}
\label{fig:basisSp0SpM}

\end{figure}
We introduce the notation $\partial\H{\Gamma} := \{0,1\}$ to indicate the endpoints of $\H{\Gamma}$.
In this section we make use of the following spline spaces:
\begin{equation}\label{eq:spaces}
\begin{split}
	\Sp^p_{0,\partial}	 &:= \Set{v \in \Sp^p(\Xi) : v = v' = 0 \text{ on } \partial\H{\Gamma}} \\
	\Sp^{p-1}_{0,\int} &:= \Set{v \in \Sp^{p-1}(\Xi^{p-1}) : v = 0 \text{ on }  \partial\H{\Gamma} \text{ and } \int_{\H{\Gamma}} v \, \ds = 0} \\
	\Sp^{p-1}_{\textsc{m},\int} &:= \Set{v \in \Sp^{p-1}(\Xi^{p-1}_{\textsc{m}}) : \int_{\H{\Gamma}} v \, \ds = 0}\\
	\Sp^{p-2}_{\textsc{m}}     &:= \Sp^{p-2}(\Xi^{p-2}_{\textsc{m}}).
\end{split}
\end{equation}
Note that the derivative operator is a bijection between $\Sp^p_{0,\partial}$ and $\Sp^{p-1}_{0,\int}$ and between $\Sp^{p-1}_{\textsc{m},\int}$ and $\Sp^{p-2}_{\textsc{m}}$.

As previously stated in Section~\ref{sec:Xhprop}, $C$ is a positive, dimensionless constant that may vary with each occurrence but remains independent of $h$ and the choice of function within the involved spaces.
\begin{prop}\label{prop:bijder2}
	For any $v \in H^1_0(\H{\Gamma})$, the following inequality holds:
	\begin{equation*}
		\norm{v}_{L^2(\H{\Gamma})} \leq C \norm{v'}_{H^{1}(\H{\Gamma})'},
	\end{equation*}
	where $H^{1}(\H{\Gamma})'$ denotes the dual space of $H^1(\H{\Gamma})$.
\end{prop}

\begin{proof}
	For all $w \in L^2(\H{\Gamma})$, there exists a function $\bar{w} \in H^1_{\int}(\H{\Gamma})$ such that $\bar{w}' = w$. Here, the space $H^1_{\int}(\H{\Gamma})$ is defined as
	\begin{equation*}
		H^1_{\int}(\H{\Gamma}) := \Set{\bar{w} \in H^1(\H{\Gamma}) : \int_{\H{\Gamma}} \bar{w} \, \ds = 0}.
	\end{equation*}
	Now, for any $v \in H^1_0(\H{\Gamma})$, we have
	\begin{equation*}
		\begin{split}
			\norm{v}_{L^2(\H{\Gamma})} &= \sup_{w \in L^2(\H{\Gamma})} \frac{\int_{\H{\Gamma}} v \, w \, \ds}{\norm{w}_{L^2(\H{\Gamma})}} 
			= \sup_{\bar{w} \in H^1_{\int}(\H{\Gamma})} \frac{\int_{\H{\Gamma}} v \bar{w}' \, \ds}{\abs{\bar{w}}_{H^1(\H{\Gamma})}} \\
			&\leq \sqrt{2} \sup_{\bar{w} \in H^1_{\int}(\H{\Gamma})} \frac{\int_{\H{\Gamma}} \bar{w} \, v' \, \ds}{\norm{\bar{w}}_{H^1(\H{\Gamma})}} 
			\leq \sqrt{2} \norm{v'}_{H^1(\H{\Gamma})'},
		\end{split}
	\end{equation*}
	where the Poincaré inequality, $\norm{\bar{w}}_{L^2(\H{\Gamma})} \leq \abs{\bar{w}}_{H^1(\H{\Gamma})}$ for all $\bar{w} \in H^1_{\int}(\H{\Gamma})$, is used in the last step.
\end{proof}
{ 
The proof of Proposition~\ref{prop:infsupl2} relies on the following technical assumption.}
{
\begin{assumption}\label{ass:technical}
Let $2 \leq p\leq 14$ and $1\leq r \leq p-1$, $Z$ be the break points corresponding to the knot vector $\Xi$. Then the first as well as the last $\lceil \frac{p+1}{p-r} \rceil$ inner break points of $Z$ are equidistant and have multiplicity $p-r$.
\end{assumption}
\begin{remark}\label{rem:p_limitation}
The upper bound $p \leq 14$ in Assumption~\ref{ass:technical} arises from our specific proof strategy rather than a fundamental limitation of the inf-sup condition. As detailed in the proof of Proposition~\ref{prop:infsupl2}, our verification relies on the positivity of the minimum generalized eigenvalue $\mu_{\min}$ in~\eqref{eq:aux05}. For polynomial degrees $p \geq 15$, this eigenvalue becomes negative for certain choices of the regularity $r$, rendering our localized algebraic approach inconclusive for verifying the estimate~\eqref{eq:aux04}.
\end{remark}}
{
For, e.g., $p=5$ and $r=3$ the knot vector looks like
\[
(0,0,0,0,0,0,h_0,h_0,2h_0,2h_0,3h_0,3h_0,\xi_{13},\ldots),
\]
with $\xi_{13} > 3h_0$. Any knot vector that is obtained from a non-uniform knot vector by sufficiently many uniform refinements (always introducing equally spaced knots of fixed multiplicity $m$) satisfies this assumption.}
\begin{prop}\label{prop:infsupl2}
{ For any $\phi \in \Sp_{0}^{p-1}$ we have}
\begin{equation*}
	\sup_{\tau \in \Sp^{p-1}_{\textsc{m}}} \frac{\int_{\H{\Gamma}} \phi \tau \, \ds}{\norm{\tau}_{L^2(\H{\Gamma})}} \geq C \norm{\phi}_{L^2(\H{\Gamma})}.
\end{equation*}
\end{prop}
\begin{proof}
 For all knot vectors $\Xi$ satisfying Assumption~\ref{ass:technical}, with $|\Xi|$ large enough, we provide a constructive proof of the statement. 
We prove that, for every $\phi \in \Sp^{p-1}_{0}$, there exists a $\tau_\phi \in \Sp^{p-1}_{\textsc{m}}$ such that 
\begin{subequations}
\begin{align}
	\norm{\tau_\phi}_{L^2(\H{\Gamma})} &\leq C_1 \norm{\phi}_{L^2(\H{\Gamma})}, \label{eq:norm1} \\
	\norm{\phi}_{L^2(\H{\Gamma})}^2 &\leq C_2 \int_{\H{\Gamma}} \phi \, \tau_\phi \, \ds. \label{eq:norm2}
\end{align}
\end{subequations}
This implies the desired result, since
\[
 \sup_{\tau \in \Sp^{p-1}_{\textsc{m}}} \frac{\int_{\H{\Gamma}} \phi \tau \, \ds}{\norm{\tau}_{L^2(\H{\Gamma})}} \geq \frac{\int_{\H{\Gamma}} \phi \tau_\phi \, \ds}{\norm{\tau_\phi}_{L^2(\H{\Gamma})}} \geq \frac{\int_{\H{\Gamma}} \phi \tau_\phi \, \ds}{C_1 \norm{\phi}_{L^2(\H{\Gamma})}} \geq \frac{1}{C_1C_2}\norm{\phi}_{L^2(\H{\Gamma})}.
\]

\textbf{Definition of $\tau_\phi$:} Let $\{\H{B}_i\}_{i=1}^n$ and $\{\widetilde{B}_i\}_{i=1}^n$ be the B-spline basis functions for $\Sp^{p-1}_{0}$ and $\Sp^{p-1}_{\textsc{m}}$, respectively. We denote by $\H{I}_i$ and $\widetilde{I}_i$ the supports of $\H{B}_i$ and $\widetilde{B}_i$, respectively, and by $\abs{\H{I}_i}$ and $\abs{\widetilde{I}_i}$ their respective lengths. Note that the set $\{\widetilde{B}_{i}\}_{i=2}^{n-1}$, forms a basis for $\Sp^{p-1}_{0}(\H{\Gamma}) \cap \Sp^{p-1}_{\textsc{m}}(\H{\Gamma})$. Consequently, 
$\{\widetilde{B}_{i}\}_{i=2}^{n-1}\cup\{ b_L,b_R \}$, with $b_L:= \H{B}_{1+p-r}$ and $b_R:= \H{B}_{n-p+r}$, forms a basis for $\Sp^{p-1}_{0}(\H{\Gamma})$.

Any $\phi \in \Sp^{p-1}_{0}(\H{\Gamma})$ can thus be expressed as
\begin{equation*}
	\phi = \alpha_L b_L + \sum_{i=2}^{n-1} c_i \widetilde{B}_{i} + \alpha_R b_R
\end{equation*}
and we define $\tau_\phi \in \Sp^{p-1}_{\textsc{m}}(\H{\Gamma})$ as
\begin{equation*}
	\tau_\phi = \alpha_L \Pi_L + \sum_{i=2}^{n-1} c_i \widetilde{B}_{i} + \alpha_R \Pi_R,
\end{equation*}
where $\Pi_L$ and $\Pi_R$ are the $L^2$-projections of $b_L$ and $b_R$ onto the spaces spanned by $\{\widetilde{B}_i\}_{i=1}^{2p-r}$ and $\{\widetilde{B}_i\}_{i=n-2p+r+1}^n$, respectively. Let $n_0$ be the smallest integer such that $\supp(\widetilde{B}_{2p-r}) \cap \supp(\widetilde{B}_{n_0-2p+r+1}) = \emptyset$. Thus, for all $n\geq n_0$, the supports of $\Pi_L$ and $\Pi_R$ are disjoint.


\textbf{Proof of \eqref{eq:norm1}:} 
Let $I \subset (0,1)$ be an element of the mesh. Assume first that $I \subseteq \supp(\Pi_L)$. 
We then have the estimate
\begin{align*}
	\norm{\tau_\phi}_{L^2(I)}^2 &\leq \int_I \left( \alpha_L \Pi_L + \sum_{i: I \subset \widetilde{I}_i} c_i \widetilde{B}_i \right)^2 \ds \\ &\leq \int_I \left( \alpha_L^2 + \sum_{i: I \subset \widetilde{I}_i} c_i^2 \right) \left( \Pi_L^2 + \sum_{i: I \subset \widetilde{I}_i} \widetilde{B}_i^2 \right) \ds.
\end{align*}
%
%
Since each basis function $\{\widetilde{B}_i\}_{i=1}^n$ satisfies $0 \leq \widetilde{B}_i \leq 1$ and the basis functions sum to one, we have the bound
\begin{equation*}
	\sum_{i : I \subset \widetilde{I}_i} \widetilde{B}_i^2 \leq \sum_{i : I \subset \widetilde{I}_i} \widetilde{B}_i \leq 1.
\end{equation*}
Furthermore, since the mesh is quasi-uniform, there exists a constant $C > 0$ independent of $h$, such that $\| \Pi_L \|^2_{L^2(I)} \leq \| b_L \|^2_{L^2(\H{\Gamma})} \leq \abs{\supp(b_L)} \leq C \abs{I}$. Thus, we obtain
\begin{equation*}
	\norm{\tau_\phi}_{L^2(I)}^2 \leq C \left( \alpha_L^2 + \sum_{i : I \subset \widetilde{I}_i} c_i^2 \right) \abs{I}.
\end{equation*}
%
Applying~\cite[Theorem 4.41]{schumaker} to $\phi$, we obtain the bound
\begin{equation*}
	\norm{\tau_\phi}_{L^2(I)}^2 \leq C \left( \alpha_L^2 + \sum_{i: I \subset \widetilde{I}_i} c_i^2 \right) \abs{I} \leq C' \sum_{i: I \subset \widetilde{I}_i} \norm{\phi}^2_{L^2(\H{I}_i)}.
\end{equation*}
Analogously, the same bound can be shown for all $I \subseteq \supp(\Pi_R)$. For all elements $I \nsubseteq \supp(b_L)\cup \supp(b_R)$ we trivially have 
$\norm{\tau_\phi}_{L^2(I)} = \norm{\phi}_{L^2(I)}$, since $\tau_\phi|_I = \phi |_I$ in that case. Summing over all $I \subset (0,1)$, we obtain the inequality in~\eqref{eq:norm1}.

\textbf{Proof of \eqref{eq:norm2}:} Since for all $2p-r+1 \leq i \leq n - 2p + r$, we have $\H{B}_i = \widetilde{B}_i$, we can write 
\begin{equation}\label{eq:aux03}
	\int_{\H{\Gamma}} \phi \, \tau_\phi \, \ds = \int_0^{\xi_{3p-r+2}} \phi \, \tau_\phi \, \ds + \|\phi\|^2_{L^2 \left( \xi_{3p-r+2}, \xi_{n-2p+r+1}\right)} + \int_{\xi_{n-2p+r+1}}^1 \phi \, \tau_\phi \, \ds.
\end{equation}
We now prove that there exists a constant $C>0$, independent of $h$, such that the first term of~\eqref{eq:aux03}satisfies
\begin{equation}\label{eq:aux04}
	\int_0^{\xi_{3p-r+2}} \phi \, \tau \, \ds \geq C \, \|\phi\|^2_{L^2 \left(0, \xi_{3p-r+2}\right)},
\end{equation}
for all $\phi \in \Sp^{p-1}_{0}(\H{\Gamma})$. Let $\boldsymbol{\phi} := (\alpha_L,c_2,\ldots,c_{n-1},\alpha_R) \in \mathbb{R}^n$ be the vector of coefficients. The generalized eigenvalue problem
\begin{equation}\label{eq:aux05}
	\mu_{\min} := \inf_{\boldsymbol{\phi} \neq 0} \frac{\boldsymbol{\phi}^t M_1 \boldsymbol{\phi}}{\boldsymbol{\phi}^t M_2 \boldsymbol{\phi}},
\end{equation}
where $M_1$ and $M_2$ are the matrices corresponding to the quadratic forms $\int_0^{\xi_{3p-r+2}} \phi \, \tau_\phi \, \ds$ and $\|\phi\|^2_{L^2 \left(0, \xi_{3p-r+2}\right)}$, respectively, yields~\eqref{eq:aux04} with $C= \mu_{\min}$. 

%

By Assumption~\ref{ass:technical}, all involved functions are defined on a uniform mesh, and all integrals can be computed explicitly for fixed values of $p$, $r$, and mesh size $h$. A scaling argument shows that the value of $\mu_{\min}$ does not depend on $h$. The dependence of $\mu_{\min}$ on both $p$ and $r$ is computed numerically and presented in Table~\ref{tab:mu_min}. 
Due to symmetry, we also have 
\begin{equation}\label{eq:aux06}
	\int_{\xi_{n-2p+r+1}}^1 \phi \, \tau_\phi \, \ds \geq C \, \|\phi\|^2_{L^2 \left(\xi_{n-2p+r+1}, 1\right)},
\end{equation}
for all $\phi \in \Sp^{p-1}_{0}$. 
\begin{table}[htbp]
\centering
\small
\rotatebox{90}{
\begin{tabular}{rccccccccccccc}
	\toprule
	r & $p=2$  & $p=3$ & $p=4$ & $p=5$ & $p=6$ & $p=7$ & $p=8$ & $p=9$ & $p=10$ & $p=11$ & $p=12$ & $p=13$ & $p=14$\\
	\midrule
	\begin{filecontents*}{tabMuMin.csv}
Row,2,3,4,5,6,7,8,9,10,11,12,13,14,15
1,0.475,0.179,0.095,0.059,0.04,0.029,0.022,0.017,0.014,0.011,0.010,0.008,0.007,0.00605
2,,0.448,0.186,0.104,0.066,0.046,0.033,0.025,0.020,0.016,0.013,0.011,0.009,0.0079
3,,,0.427,0.189,0.104,0.066,0.045,0.032,0.024,0.018,0.014,0.011,0.009,0.00749
4,,,,0.409,0.184,0.109,0.063,0.042,0.028,0.020,0.014,0.010,0.007,0.00515
5,,,,,0.395,0.183,0.104,0.070,0.037,0.023,0.014,0.008,0.004,0.00159
6,,,,,,0.383,0.176,0.098,0.064,0.045,0.018,0.008,0.002,-0.0022
7,,,,,,,0.372,0.176,0.102,0.057,0.039,0.028,0.002,-0.00483
8,,,,,,,,0.363,0.169,0.096,0.050,0.032,0.021,0.015
9,,,,,,,,,0.355,0.168,0.089,0.059,0.024,0.014
10,,,,,,,,,,0.348,0.162,0.093,0.052,0.017
11,,,,,,,,,,,0.341,0.162,0.087,0.045
12,,,,,,,,,,,,0.335,0.156,0.081
13,,,,,,,,,,,,,0.330,0.156
\end{filecontents*}
\csvreader[late after line=\\,
	head to column names]
	{tabMuMin.csv}{2=\ppa, 3=\ppb, 4=\ppc, 5=\ppd, 6=\ppe, 7=\ppf, 8=\ppg, 9=\pph, 10=\ppi, 11=\ppj, 12=\ppk, 13=\ppl, 14=\ppm}{\Row & {\ppa} & {\ppb} & {\ppc} & {\ppd} & {\ppe} & {\ppf} & {\ppg} & {\pph} & {\ppi} & {\ppj} & {\ppk} & {\ppl} & {\ppm}}
	\bottomrule
\end{tabular}
}
\caption{Values of $\mu_{\min}$ obtained from~\eqref{eq:aux05} using a uniform mesh.}
 \label{tab:mu_min}
\end{table}
Finally, by combining equations~\eqref{eq:aux03}, \eqref{eq:aux04}, and~\eqref{eq:aux06}, we obtain equation~\eqref{eq:norm2}.

\textbf{Proof for small $n$:} For all $n< n_0$, the estimate is trivially satisfied for some $C(n_0)>0$. Since $n_0$ depends only on $p$ and $r$, the desired result follows.
\end{proof}
We observe that higher regularity leads to better stability (higher values of $\mu_{\min}$). Numerical experiments (not reported) also suggest that $\mu_{\min}$ remains stable even for  quasi-uniform meshes.
It is important to observe that meshes initially generated from non-uniform knot configurations become uniform near the boundary after a sufficient number of uniform refinement steps. As a result, Assumption~\ref{ass:technical} is typically satisfied in practical applications. Therefore, the stability result of Proposition~\ref{prop:infsupl2} extends to general meshes, provided the mesh size $h$ is sufficiently small.

\begin{lemma}\label{lem:bijection}
The operator $\widecheck{\M}^{*}_h  \colon L^2(\H{\Gamma}) \to \Sp^{p-1}_{\textsc{m}}$ defined by
\begin{equation*}
	\int_{\H{\Gamma}} \widecheck{\M}^{*}_h \phi \, \psi \, \ds = \int_{\H{\Gamma}} \phi \, \psi \, \ds, \quad \forall \psi \in \Sp^{p-1}_{0},
\end{equation*}
is well-defined, fulfills
\begin{equation*}
	\| \widecheck{\M}^{*}_h \phi \|_{L^2(\H{\Gamma})} \leq C \| \phi \|_{L^2(\H{\Gamma})}, \, \forall \phi \in {L^2(\H{\Gamma})} 
\end{equation*}
and is bijective from $\Sp^{p-1}_{0}$ into $ \Sp^{p-1}_{\textsc{m}}$.
\end{lemma}
\begin{proof}
To demonstrate that the operator is well-defined, assume there exist two distinct elements $\tau_1, \tau_2 \in \Sp^{p-1}_{\textsc{m}}$ such that
\begin{equation*}
	\int_{\H{\Gamma}} \tau_k \, \psi \, \ds = \int_{\H{\Gamma}} \phi \, \psi \, \ds, \quad \forall \psi \in \Sp^{p-1}_{0}.
\end{equation*}
By Proposition~\ref{prop:infsupl2}, we have
\begin{equation*}
	\|\tau_1 - \tau_2\|_{L^2(\H{\Gamma})} \leq C \sup_{\psi \in \Sp^{p-1}_{0},\, \|\psi\|_{L^2}=1} \int_{\H{\Gamma}} (\tau_1 - \tau_2) \, \psi \, \ds = 0.
\end{equation*}
This implies $\tau_1 = \tau_2$, confirming that the operator is well-defined.
The $L^2$-stability follows similarly to the proof of Lemma~\ref{lemma:dualproperties}.
To establish bijectivity as an operator between $\Sp^{p-1}_{0}$ and $\Sp^{p-1}_{\textsc{m}}$, note that
\begin{equation*}
	\|\psi\|^2_{L^2(\H{\Gamma})} = \int_{\H{\Gamma}} \widecheck{\M}^{*}_h \psi \, \psi \, \ds \leq \|\widecheck{\M}^{*}_h \psi\|_{L^2(\H{\Gamma})} \|\psi\|_{L^2(\H{\Gamma})}.
\end{equation*}
Consequently, we obtain
\begin{equation*}
	\|\psi\|_{L^2(\H{\Gamma})} \leq \|\widecheck{\M}^{*}_h \psi\|_{L^2(\H{\Gamma})}, \quad \forall \psi \in \Sp^{p-1}_{0}.
\end{equation*}
If there exist two distinct functions $\psi_1, \psi_2 \in \Sp^{p-1}_{0}$ such that $\widecheck{\M}^{*}_h \psi_1 = \widecheck{\M}^{*}_h \psi_2$, it follows that $\|\psi_1 - \psi_2\|_{L^2(\H{\Gamma})} = 0$, implying $\psi_1 = \psi_2$. Since $\dim(\Sp^{p-1}_{0}) = \dim(\Sp^{p-1}_{\textsc{m}})$, surjectivity also follows.
\end{proof}
\begin{prop}\label{prop:infsuph1}
For any $\phi \in \Sp^{p-1}_{0}$, the following inequality holds:
\begin{equation}\label{eq:infsuph1}
	\sup_{\tau \in \Sp^{p-1}_{\textsc{m}}} \frac{\int_{\H{\Gamma}} \phi \tau \, \ds}{\norm{\tau}_{H^1(\H{\Gamma})}} \geq C \norm{\phi}_{H^1(\H{\Gamma})'}.
\end{equation}
\end{prop}
\begin{proof}
From Lemma~\ref{lem:bijection}, we know that the operator $\widecheck{\M}^{*}_h$ is $L^2$-stable. We now verify that it is also $H^1$-stable. Let $\tau \in H^1(\H{\Gamma})$, and let $\widetilde{\Pi}_h: H^1(\H{\Gamma}) \rightarrow \Sp^{p-1}_{\textsc{m}}$ denote the quasi-interpolant operator defined { using the univariate dual basis from~\citeRefDetail{Theorem 4.41}{schumaker}, cf. Proposition~\ref{prop:dual-basis}}. 
Recall that for $t = 0, 1$, it holds that
\begin{equation}\label{eq:qiest}
	\norm{\tau - \widetilde{\Pi}_h \tau}_{H^t(\H{\Gamma})} \leq C h^{1-t} \norm{\tau}_{H^1(\H{\Gamma})}.
\end{equation}
Using the inverse inequality for splines, we derive:
\begin{equation*}
\begin{split}
	\norm{\widecheck{\M}^{*}_h \tau}_{H^1(\H{\Gamma})} &\leq \norm{\widecheck{\M}^{*}_h (\tau - \widetilde{\Pi}_h \tau)}_{H^1(\H{\Gamma})} + \norm{\widecheck{\M}^{*}_h \widetilde{\Pi}_h \tau}_{H^1(\H{\Gamma})} \\
	&\leq C h^{-1} \norm{\tau - \widetilde{\Pi}_h \tau}_{L^2(\H{\Gamma})} + \norm{\widetilde{\Pi}_h \tau}_{H^1(\H{\Gamma})}.
\end{split}
\end{equation*}
Applying $L^2$-stability and inequality~\eqref{eq:qiest}, we have:
\begin{equation*}
\begin{split}
	\norm{\widecheck{\M}^{*}_h \tau}_{H^1(\H{\Gamma})} &\leq C h^{-1} \norm{\tau - \widetilde{\Pi}_h \tau}_{L^2(\H{\Gamma})} + \norm{\widetilde{\Pi}_h \tau}_{H^1(\H{\Gamma})} \\
	&\leq C \norm{\tau}_{H^1(\H{\Gamma})}.
\end{split}
\end{equation*}
For any $\phi \in H^1(\H{\Gamma})'$, it holds that:
\begin{equation}\label{eq:aux02}
	\norm{\phi}_{H^1(\H{\Gamma})'} = \sup_{\eta \in H^1(\H{\Gamma})} \frac{\int_{\H{\Gamma}} \eta \phi \, \ds}{\norm{\eta}_{H^1(\H{\Gamma})}},
\end{equation}
and, in particular, this holds for any $\phi \in \Sp^{p-1}_{0}$. By Lemma~\ref{lem:bijection}, for any $\psi \in \Sp^{p-1}_{\textsc{m}}$, there exists $\eta \in \Sp^{p-1}_{0}$ such that $\psi = \widecheck{\M}^{*}_h \eta$, which implies:
\begin{equation*}
	\int_{\H{\Gamma}} \psi \phi \, \ds = \int_{\H{\Gamma}} \eta \phi \, \ds, \quad \forall \phi \in \Sp^{p-1}_{0},
\end{equation*}
and $\norm{\psi}_{H^1(\H{\Gamma})} \leq C^{-1} \norm{\eta}_{H^1(\H{\Gamma})}$. 
Thus, for all $\phi \in \Sp^{p-1}_{0}$, there exists $\psi = \widecheck{\M}^{*}_h \eta \in \Sp^{p-1}_{\textsc{m}}$ such that:
\begin{equation*}
	\frac{\int_{\H{\Gamma}} \psi \phi \, \ds}{\norm{\psi}_{H^1(\H{\Gamma})}} \geq C \frac{\int_{\H{\Gamma}} \eta \phi \, \ds}{\norm{\eta}_{H^1(\H{\Gamma})}}.
\end{equation*}
Taking the supremum over $\eta \in H^1(\H{\Gamma})$ and applying~\eqref{eq:aux02}, we establish~\eqref{eq:infsuph1}.
\end{proof}
\begin{prop}\label{prop:infsupconst}
	For any $\tau \in \Sp^{p-1}_{\textsc{m},\int}$, the following inequality holds:
	\begin{equation*}
		\sup_{\phi \in \Sp^{p-1}_{0,\int}} \frac{\int_{\H{\Gamma}} \phi \tau \, \ds}{\norm{\phi}_{H^1(\H{\Gamma})'}} \geq C \norm{\tau}_{H^1(\H{\Gamma})}.
	\end{equation*}
\end{prop}

\begin{proof}
	From Proposition~\ref{prop:infsuph1}, we know that for all $\phi \in \Sp^{p-1}_{0}$, there exists $\tau \in \Sp^{p-1}_{\textsc{m}}$ such that
	\begin{equation*}
		\int_{\H{\Gamma}} \phi \tau \, \ds \geq C \norm{\phi}_{H^1(\H{\Gamma})'} \norm{\tau}_{H^1(\H{\Gamma})}.
	\end{equation*}
	This inequality holds in particular for all $\phi \in \Sp^{p-1}_{0,\int}$. Now, consider $\tau^* := \tau - \bar{\tau}$, where $\bar{\tau}$ denotes the mean value of $\tau$ on $\H{\Gamma}$. Substituting $\tau^*$, we have
	\begin{equation*}
		\int_{\H{\Gamma}} \phi \tau^* \, \ds = \int_{\H{\Gamma}} \phi \tau \, \ds \geq C \norm{\phi}_{H^1(\H{\Gamma})'} \norm{\tau}_{H^1(\H{\Gamma})}.
	\end{equation*}
	It is straightforward to verify that $\norm{\tau^*}_{H^1(\H{\Gamma})} \leq \norm{\tau}_{H^1(\H{\Gamma})}$. Thus, for all $\phi \in \Sp^{p-1}_{0,\int}$, we obtain
	\begin{equation*}
		\sup_{\tau^* \in \Sp^{p-1}_{\textsc{m},\int}} \frac{\int_{\H{\Gamma}} \phi \tau^* \, \ds}{\norm{\tau^*}_{H^1(\H{\Gamma})}} \geq C \norm{\phi}_{H^1(\H{\Gamma})'}.
	\end{equation*}
	Finally, as discussed in Remark~\ref{rem:bothstab}, the interchange of the inf-sup spaces leads to the desired result.
\end{proof}
\begin{theorem}\label{theo:infsupparam}
For any $\mu \in \Sp^{p-2}_{\textsc{m}}$, the following inequality holds:
\begin{equation}\label{eq:sup2}
	\sup_{\psi \in \Sp^p_{0,\partial}} \frac{\int_{\H{\Gamma}} \psi \mu \, \ds}{\norm{\psi}_{L^2(\H{\Gamma})}} \geq C \norm{\mu}_{L^2(\H{\Gamma})}.
\end{equation}
\end{theorem}
\begin{proof}
We begin by recalling Proposition~\ref{prop:infsupconst}, which states that for any $\tau \in \Sp^{p-1}_{\textsc{m},\int}$, the following inequality holds:
\begin{equation}\label{eq:sup1}
	\sup_{\psi \in \Sp^p_{0,\partial}} \frac{\int_{\H{\Gamma}} \psi' \tau}{\norm{\psi'}_{H^1(\H{\Gamma})'}} \geq C \norm{\tau}_{H^1(\H{\Gamma})},
\end{equation}
where we use the bijectivity between $\Sp^p_{0,\partial}$ and $\Sp^{p-1}_{0,\int}$ via the derivative operator, denoting $\psi' = \phi$.\\
By integration by parts and applying Proposition~\ref{prop:bijder2}, the quotient in~\eqref{eq:sup1} can be estimated as:
\begin{equation*}
	\frac{\int_{\H{\Gamma}} \psi' \tau}{\norm{\psi'}_{H^1(\H{\Gamma})'}} \leq C \frac{\int_{\H{\Gamma}} \psi \tau'}{\norm{\psi}_{L^2(\H{\Gamma})}} = C \frac{\int_{\H{\Gamma}} \psi \mu}{\norm{\psi}_{L^2(\H{\Gamma})}},
\end{equation*}
where we set $\mu = \tau'$. Here, $\mu$ is an arbitrary function in $\Sp^{p-2}_{\textsc{m}}$, and it satisfies
\begin{equation*}
	\norm{\mu}_{L^2(\H{\Gamma})} = \abs{\tau}_{H^1(\H{\Gamma})} \leq \norm{\tau}_{H^1(\H{\Gamma})}.
\end{equation*}
Thus, inequality~\eqref{eq:sup2} follows directly.
\end{proof}

We finally introduce a specific  multiplier space $M_\ell$ through the following assumption.
\begin{assumption} \label{ass:multipliers-defn} Let $\mathcal{F}_{\sell}|_{\H{\Gamma}}$ be the mapping of $\H{\Gamma}$ into the edge $\Gamma_\ell = \partial\Omega_{\ml} \cap \partial\Omega_{\sell}$, we assume
	\begin{equation*}
	M_\ell = \left\{ \mu \in L^2(\Gamma_\ell) : \mu \circ \mathcal{F}_{\sell}|_{\H{\Gamma}} = \H{\mu} \in \mathcal{S}^{p-2}_{\textsc{m}}(\H{\Gamma}) \right\}.
\end{equation*}
\end{assumption}
{
To proceed, we require a super-approximation property for local projections onto the spline space $\Sp^{p-2}_{\textsc{m}}(\H{\Gamma})$. Inspired by the approach in~\citeOneRef{brivadis}, the following lemma extends classical finite element approximation bounds to our isogeometric setting.

\begin{lemma}\label{lem:superapp}
Let $\H{\pi}^{*} \colon L^2(\H{\Gamma}) \to \Sp^{p-2}_{\textsc{m}}(\H{\Gamma})$ be a local projection operator with optimal approximation properties. Assume the weight function $w$ is smooth in the interior of each mesh element and continuous across the mesh lines. Under these conditions, the following super-approximation estimate holds:
\begin{equation*}
	\norm{\H{\mu}w - \H{\pi}^{*}(\H{\mu} w)}_{L^2(\H{\Gamma})} \leq C^* h \norm{\H{\mu}}_{L^2(\H{\Gamma})}, \quad \forall \H{\mu} \in \Sp^{p-2}_{\textsc{m}}(\H{\Gamma}).
\end{equation*}
\end{lemma}
\begin{proof}
This result follows by adapting the classical super-approximation arguments of~\cite[Theorem 2.3.1]{wahlbin} to the isogeometric setting, leveraging standard spline approximation theory (see, e.g.,~\citeOneRef{bazilevs}).
\end{proof}
The discrete inf-sup stability property of the pairing $\mathcal{S}^p_{0,\partial}/\mathcal{S}^{p-2}_{\textsc{m}}$ on $\H{\Gamma}$ passes to $W_{\ell}/M_{\ell}$, thanks to the following results.
\begin{lemma}\label{lem:geom_reduction}
Let $\mu \in M_\ell$, and let $\H{\mu}$ be its parametric counterpart defined on the reference boundary $\H{\Gamma} = \{0\} \times [0,1]$ via the mapping $\F_{\sell}$. Under Assumption~\ref{ass:unif_reg}, the continuous inf-sup condition on the physical domain satisfies the following lower bound in the parametric domain:
\begin{equation*}
	\sup_{\psi \in W_\ell} \frac{\int_{\Gamma_\ell} \psi \, \mu \, \ds}{\norm{\psi}_{L^2(\Gamma_\ell)}} \geq C^{-1} \sup_{\H{\eta} \in \Sp^p_{0,\partial}(\H{\Gamma})} \frac{\int_{\H{\Gamma}} (\H{\mu} \, \alpha_{\sell} \H{\rho}) \, \H{\eta} \, \mathrm{d}\H{y}}{\norm{\H{\eta}}_{L^2(\H{\Gamma})}},
\end{equation*}
where $\H{\rho}(\H{y}) := \norm{\partial_{\H{y}} \F_{\sell}(0,\H{y})}$, and the geometric scaling coefficient is:
\begin{equation*}
	\alpha_{\sell}(\H{y}) := \det \left[ \partial_{\H{x}} \F_{\sell}(0,\H{y}) \quad \partial_{\H{y}} \F_{\sell}(0,\H{y}) \right]^{-1}.
\end{equation*}
\end{lemma}
\begin{proof}
For any $\psi \in W_\ell$, by definition~\eqref{eq:defWl}, there exists a function $\varphi \in \widetilde{X}_{\sell,h}$ such that $\psi = \dn \varphi|_{\Gamma_\ell}$, where 
\begin{equation*}
	\widetilde{X}_{\sell,h} := \left\{ \varphi \in X_{\sell,h} \cap H^1_0(\Omega_{\sell}) :\; \dn \varphi|_{\partial \Omega_{\sell} \setminus \Gamma_\ell} = 0 \right\}.
\end{equation*}
Assume, without loss of generality, that $\Gamma_\ell$ is parameterized by $\F_{\sell}(0,\H{y})$ for $\H{y} \in [0,1]$. We can write $\varphi$ in parametric coordinates as $\varphi = \H{\varphi} \circ \F_{\sell}^{-1}$ with $\H{\varphi} \in \H{Y}$, where
\begin{equation*}
	\H{Y} := \left\{ \H{\varphi} \in \Sp^{\mathbf{p}}_k \cap H^1_0(\H{\Omega}) :\; \partial_{\H{n}} \H{\varphi}|_{\partial \H{\Omega} \setminus \H{\Gamma}} = 0 \right\}.
\end{equation*}
Following the framework in~\citeOneRef{collin}, the normal derivative transforms using $\alpha_{\sell}$ and the cross-scaling coefficient $\beta_{\sell}$, defined as:
\begin{equation*}
	\beta_{\sell}(\H{y}) := - \frac{\partial_{\H{x}} \F_{\sell}(0,\H{y}) \cdot \partial_{\H{y}} \F_{\sell}(0,\H{y})}{{\H{\rho}}^2 \det \left[ \partial_{\H{x}} \F_{\sell}(0,\H{y}) \quad \partial_{\H{y}} \F_{\sell}(0,\H{y}) \right]}.
\end{equation*}
Using these mapping components, we rewrite the continuous supremum in the parametric domain $\H{\Gamma}$:
\begin{equation*}
	\sup_{\psi \in W_\ell} \frac{\int_{\Gamma_\ell} \psi \, \mu \, \ds}{\norm{\psi}_{L^2(\Gamma_\ell)}} = \sup_{\H{\varphi} \in \H{Y}} \frac{\int_{\H{\Gamma}} \left( \alpha_{\sell} \partial_{\H{x}} \H{\varphi} + \beta_{\sell} \partial_{\H{y}} \H{\varphi} \right) \H{\mu} \, \H{\rho} \, \mathrm{d}\H{y}}{\norm{\left( \alpha_{\sell} \partial_{\H{x}} \H{\varphi} + \beta_{\sell} \partial_{\H{y}} \H{\varphi} \right) \sqrt{\H{\rho}}}_{L^2(\H{\Gamma})}}.
\end{equation*}
Since $\partial_{\H{y}} \H{\varphi}|_{\H{\Gamma}} = 0$ for all $\H{\varphi} \in \H{Y}$, the term with $\beta_{\sell}$ is identically zero on the boundary. We set $\H{\eta} := \partial_{\H{x}} \H{\varphi}|_{\H{\Gamma}} \in \Sp^p_{0,\partial}(\H{\Gamma})$ and apply the uniform mesh regularity from Assumption~\ref{ass:unif_reg} to bound the denominator. This completes the proof.
\end{proof}

\begin{theorem}
Under Assumption~\ref{ass:multipliers-defn}, there exists $\bar{h} > 0$ such that for all $h \leq \bar{h}$ and for each $\mu \in M_{\ell}$,
\begin{equation*}
	\sup_{\psi \in W_{\ell}} \frac{\int_{\Gamma_\ell} \psi \mu \, \ds}{\|\psi\|_{L^2(\Gamma_\ell)}} \geq C \|\mu\|_{L^2(\Gamma_\ell)}.
\end{equation*}
\end{theorem}
\begin{proof}
We bound the continuous supremum using Lemma~\ref{lem:geom_reduction}:
\begin{equation}\label{eq:proof_start}
	\sup_{\psi \in W_\ell} \frac{\int_{\Gamma_\ell} \psi \, \mu \, \ds}{\norm{\psi}_{L^2(\Gamma_\ell)}} \geq C^{-1} \sup_{\H{\eta} \in \Sp^p_{0,\partial}(\H{\Gamma})} \frac{\int_{\H{\Gamma}} (\H{\mu} \, \alpha_{\sell} \H{\rho}) \, \H{\eta} \, \mathrm{d}\H{y}}{\norm{\H{\eta}}_{L^2(\H{\Gamma})}}.
\end{equation}
Let $w = \alpha_{\sell}\H{\rho}$. The domain parametrizations are smooth inside the elements and have $C^{r-1}$ regularity across the element boundaries. Therefore, $w$ satisfies the conditions of Lemma~\ref{lem:superapp}. The local projection operator $\H{\pi}^{*} \colon L^2(\H{\Gamma}) \to \Sp^{p-2}_{\textsc{m}}(\H{\Gamma})$ satisfies the super-approximation estimate:
\begin{equation}\label{eq:thm_superapp}
	\norm{\H{\mu} \alpha_{\sell}\H{\rho} - \H{\pi}^{*}(\H{\mu} \alpha_{\sell}\H{\rho})}_{L^2(\H{\Gamma})} \leq C^* h \norm{\H{\mu}}_{L^2(\H{\Gamma})}.
\end{equation}
We add and subtract the projection $\H{\pi}^{*}(\H{\mu} \alpha_{\sell}\H{\rho})$ in the numerator of~\eqref{eq:proof_start}. Then, we use the Cauchy--Schwarz inequality to obtain:
\begin{equation*}
	\frac{\int_{\H{\Gamma}} (\H{\mu} \, \alpha_{\sell} \H{\rho}) \, \H{\eta} \, \mathrm{d}\H{y}}{\norm{\H{\eta}}_{L^2(\H{\Gamma})}} \geq \frac{\int_{\H{\Gamma}} \H{\pi}^{*}(\H{\mu} \, \alpha_{\sell} \H{\rho}) \, \H{\eta} \, \mathrm{d}\H{y}}{\norm{\H{\eta}}_{L^2(\H{\Gamma})}} - \norm{\H{\mu} \, \alpha_{\sell} \H{\rho}-\H{\pi}^{*}(\H{\mu} \, \alpha_{\sell} \H{\rho})}_{L^2(\H{\Gamma})}.
\end{equation*}
Taking the supremum over $\H{\eta} \in \Sp^p_{0,\partial}(\H{\Gamma})$ and using the discrete inf-sup condition from Theorem~\ref{theo:infsupparam}, we get:
\begin{equation*}
	\sup_{\H{\eta} \in \Sp^p_{0,\partial}(\H{\Gamma})} \frac{\int_{\H{\Gamma}} (\H{\mu} \, \alpha_{\sell} \H{\rho}) \, \H{\eta} \, \mathrm{d}\H{y}}{\norm{\H{\eta}}_{L^2(\H{\Gamma})}} \geq
	C \norm{\H{\pi}^{*} (\H{\mu} \, \alpha_{\sell} \H{\rho})}_{L^2(\H{\Gamma})} - \norm{\H{\mu} \, \alpha_{\sell} \H{\rho}-\H{\pi}^{*}(\H{\mu} \, \alpha_{\sell} \H{\rho})}_{L^2(\H{\Gamma})}.
\end{equation*}
Using the super-approximation bound~\eqref{eq:thm_superapp}, for $h \leq \bar{h}$ small enough, we have:
\begin{equation*}
\begin{aligned}
	\sup_{\H{\eta} \in \Sp^p_{0,\partial}(\H{\Gamma})} \frac{\int_{\H{\Gamma}} (\H{\mu} \, \alpha_{\sell} \H{\rho}) \, \H{\eta} \, \mathrm{d}\H{y}}{\norm{\H{\eta}}_{L^2(\H{\Gamma})}} &\geq
	C \norm{\H{\mu} \, \alpha_{\sell} \H{\rho}}_{L^2(\H{\Gamma})} - (C+1) \norm{\H{\mu} \, \alpha_{\sell} \H{\rho}-\H{\pi}^{*}(\H{\mu} \, \alpha_{\sell} \H{\rho})}_{L^2(\H{\Gamma})} \\
	&\geq C \norm{\H{\mu}}_{L^2(\H{\Gamma})} - C^* h \norm{\H{\mu}}_{L^2(\H{\Gamma})} \\
	&\geq C \norm{\H{\mu}}_{L^2(\H{\Gamma})}.
\end{aligned}	
\end{equation*}
Finally, the result follows from the standard equivalence of norms between the reference domain $\H{\Gamma}$ and the physical boundary $\Gamma_\ell$.
\end{proof}
}
\begin{corollary}
Under Assumption \ref{ass:multipliers-defn}, there exist $\bar h > 0 $ such that for all $h\leq \bar h$ and for each $\mu \in M_h$, the following inequality holds:

\begin{equation*}
	\sup_{\psi \in W_{h}} \frac{\int_{\Sigma} \psi \mu \ds}{\norm{\psi}_{L^2(\Sigma)}} \geq C \norm{\mu}_{L^2(\Sigma)}.
\end{equation*}
\end{corollary}
	\section{Numerical results}\label{sec:numerics}
\begin{figure}[h]
	\centering
	\begin{subfigure}{0.45\linewidth}
		\includegraphics[width=\linewidth]{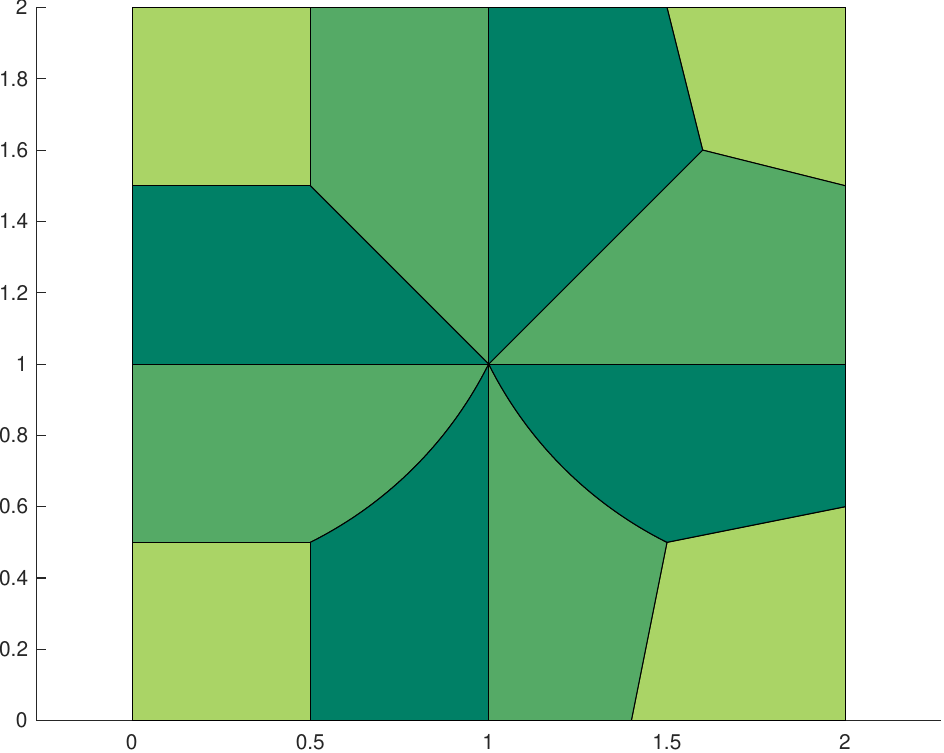}
		\caption{Square.}
		\label{fig:Square_mp}
	\end{subfigure}
	\hfill
	\begin{subfigure}{0.45\linewidth}
		\includegraphics[width=\linewidth]{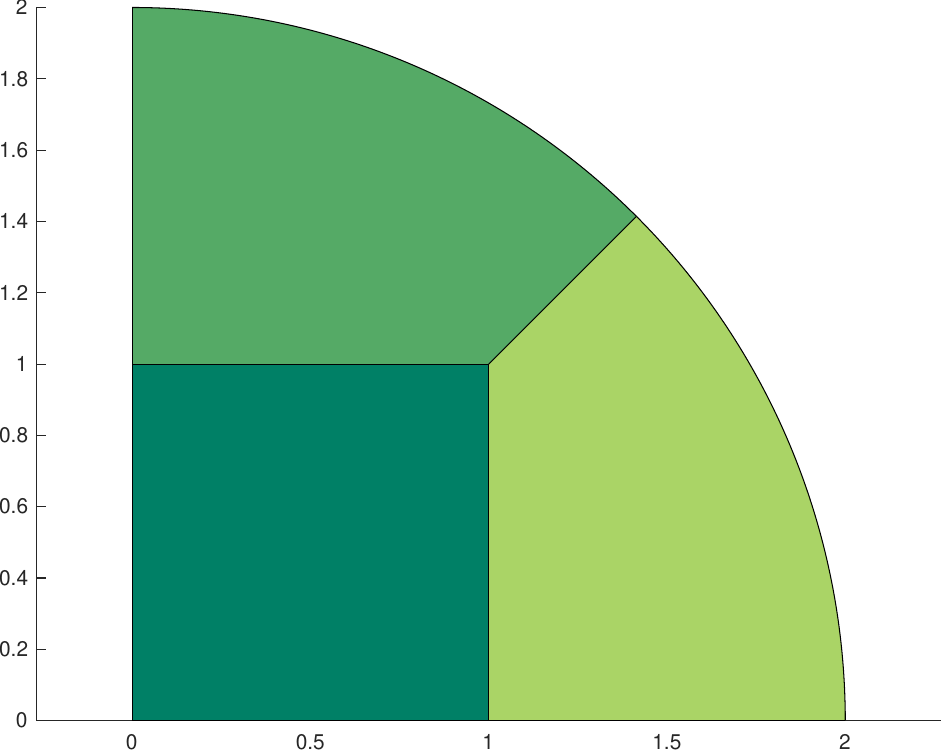}
		\caption{Quarter circle.}
		\label{fig:Quartercircle_mp}
	\end{subfigure}
	\\
	\vspace{5mm}
	\begin{subfigure}{0.45\linewidth}
		\includegraphics[width=\linewidth]{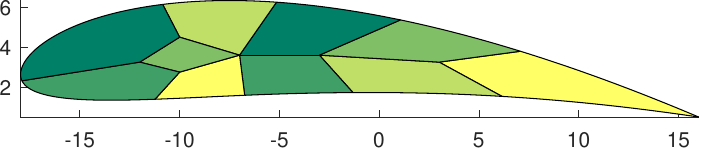}
		\caption{Airfoil.}
		\label{fig:Airfoil_mp}
	\end{subfigure}
	\caption{Physical domains.}
	\label{fig:Domains}
\end{figure}

In this section, we present numerical tests validating the optimal error estimate stated in Corollary~\ref{thm:orderofconv} for the pairing $\mathcal{S}^{p} / \mathcal{S}_{\textsc{m}}^{p-2}$, discussed in Section~\ref{sec:dscinfsup}. Additionally, we demonstrate that optimal convergence is achieved for pairings not covered by the theoretical framework.\\
All tests were conducted using Matlab R2024a and the GeoPDEs toolbox~\citeOneRef{defalco}. The domains considered are depicted in Figures~\ref{fig:Domains}: a 12-patch square (\subref{fig:Square_mp}), a 3-patch quarter-circle (\subref{fig:Quartercircle_mp}), and a 10-patch airfoil (\subref{fig:Airfoil_mp}). For all simulations, boundary conditions and source terms were chosen to produce the exact solution ${u_{\mathrm{ex}}(x,y) = \cos(x) \cos(y)}$.\\
For simplicity, previous sections assumed that the domain is described by a spline geometrical mapping (see Assumption~\ref{ass:unif_reg}). However, this assumption is not critical to the theoretical results. Numerical experiments, such as those on the quarter-circle domain, indicate that stability and accuracy are preserved on NURBS-parametrized domains.\\
For all tests, results obtained with higher-degree, refined meshes that are corrupted by round-off errors are excluded from the plots.
\subsection{{ Tests with} $C^2$-continuity at vertices}
We begin with a primal space $X_h$, defined as in~\eqref{eq:defXkh} and~\eqref{eq:defXh}, i.e., as the space of functions constructed from spline spaces of degree $p$ in each patch, with $C^0$ regularity across interfaces and strong $C^2$-continuity enforced at the vertices.\\ 
The multiplier space on the interfaces is constructed as described in Assumption~\ref{ass:multipliers-defn}.
From a theoretical perspective, this pairing satisfies the discrete inf-sup stability condition, as demonstrated in Section~\ref{sec:dscinfsup}.\\
Figures~\ref{fig:Square_C2},~\ref{fig:Airfoil_C2}, and~\ref{fig:Quartercircle_C2} illustrate the convergence plots in patch-wise $H^2$ broken norm (denoted by $\mathcal{H}^2$), as well as in $H^1$, $L^2$, and $L^{\infty}$ norms, for primal spaces of degrees $p = 2, \dots, 5$, with maximal global regularity $r = p - 1$. Figure~\ref{fig:Quartercircle_C2R1} shows results for spline spaces that are globally only $C^1$-continuous, i.e. $r=1$.\\
In all cases for $p = 2$, the rates of convergence are suboptimal in $L^2$, a phenomenon observed even in the single-patch case (see~\citeOneRef{strang}). In all other cases, the computed rates of convergence are optimal.
\tikzstyle{Linea1}=[thick,mark=*]
\tikzstyle{Linea2}=[thick,dashed]


\pgfplotscreateplotcyclelist{Lista1}{%
	{Linea1,Red},
	{Linea2,Red},
	{Linea1,Green},
	{Linea2,Green},
	{Linea1,Cyan},
	{Linea2,Cyan},
	{Linea1,Violet},
	{Linea2,Violet}}

\def \DATAFILE {ErrSquare.csv}

\begin{figure}[htbp]
\centering
\hspace*{\fill}
	\begin{subfigure}[t]{0.45\linewidth}
	\centering
		\begin{tikzpicture}[font=\small, trim axis left]
			\begin{loglogaxis}[
				cycle list name=Lista1,
				xmajorgrids=true,
				ymajorgrids=true,
				width=\linewidth,
				height=1.25\linewidth,
				xlabel={$h$},
				xminorticks=false,
				yminorticks=false,
				ymin={1e-14},
				ylabel={$\nicefrac{\|u-u_h\|_{\mathcal{H}^2}}{\|u\|_{\mathcal{H}^2}}$},
				legend columns=2,
				legend pos=south east,
				legend style={at={(0.99,0.01)},anchor=south east},
				legend entries={$\Sp^2/\Sp^0_{\textsc{m}}$,$O(h)$,$\Sp^3/\Sp^1_{\textsc{m}}$,$O(h^2)$,$\Sp^4/\Sp^2_{\textsc{m}}$,$O(h^3)$,$\Sp^5/\Sp^{3}_{\textsc{m}}$,$O(h^4)$}]
				\addplot table [x=h, y=H2P2, col sep=comma] {\DATAFILE};
				\addplot table [x=h, y=h1H2, col sep=comma] {\DATAFILE};
				\addplot table [x=h, y=H2P3, col sep=comma] {\DATAFILE};
				\addplot table [x=h, y=h2H2, col sep=comma] {\DATAFILE};
				\addplot table [x=h, y=H2P4, col sep=comma] {\DATAFILE};
				\addplot table [x=h, y=h3H2, col sep=comma] {\DATAFILE};
				\addplot table [x=h, y=H2P5, col sep=comma] {\DATAFILE};
				\addplot table [x=h, y expr=0.1*\thisrow{h4H2}, col sep=comma] {\DATAFILE};
			\end{loglogaxis}
		\end{tikzpicture}
		\label{fig:Square_H2_C0}
	\end{subfigure}
\hfill
	\begin{subfigure}[t]{0.45\linewidth}
		\centering
		\begin{tikzpicture}[font=\small, trim axis right]
			\begin{loglogaxis}[
				cycle list name=Lista1,
				xmajorgrids=true,
				ymajorgrids=true,
				width=\linewidth,
				height=1.25\linewidth,
				xlabel={$h$},
				xminorticks=false,
				yminorticks=false,
				ymin={1e-16},
				ylabel={$\nicefrac{\|u-u_h\|_{H^1}}{\|u\|_{H^1}}$},
				legend columns=2,
				legend pos=south east,
				legend style={at={(0.99,0.01)},anchor=south east},
				legend entries={$\Sp^2/\Sp^0_{\textsc{m}}$,$O(h^2)$,$\Sp^3/\Sp^1_{\textsc{m}}$,$O(h^3)$,$\Sp^4/\Sp^2_{\textsc{m}}$,$O(h^4)$,$\Sp^5/\Sp^{3}_{\textsc{m}}$,$O(h^5)$}]
				\addplot table [x=h, y=H1P2, col sep=comma] {\DATAFILE};
				\addplot table [x=h, y=h2H1, col sep=comma] {\DATAFILE};
				\addplot table [x=h, y=H1P3, col sep=comma] {\DATAFILE};
				\addplot table [x=h, y=h3H1, col sep=comma] {\DATAFILE};
				\addplot table [x=h, y=H1P4, col sep=comma] {\DATAFILE};
				\addplot table [x=h, y=h4H1, col sep=comma] {\DATAFILE};
				\addplot table [x=h, y=H1P5, col sep=comma, restrict expr to domain={\coordindex}{0:2}] {\DATAFILE};
				\addplot table [x=h, y expr=0.001*\thisrow{h5H1}, col sep=comma] {\DATAFILE};
			\end{loglogaxis}
		\end{tikzpicture}
		\label{fig:Square_H1_C2}
	\end{subfigure}
\\
\hspace*{\fill}
	\begin{subfigure}[t]{0.45\linewidth}
		\centering
		\begin{tikzpicture}[font=\small, trim axis left]
			\begin{loglogaxis}[
				cycle list name=Lista1,
				xmajorgrids=true,
				ymajorgrids=true,
				width=\linewidth,
				height=1.25\linewidth,
				xlabel={$h$},
				xminorticks=false,
				yminorticks=false,
				ymin={1e-18},
				ylabel={$\nicefrac{\|u-u_h\|_{L^2}}{\|u\|_{L^2}}$},
				legend columns=2,
				legend pos=south east,
				legend style={at={(0.99,0.01)},anchor=south east},
				legend entries={$\Sp^2/\Sp^0_{\textsc{m}}$,$O(h^3)$,$\Sp^3/\Sp^1_{\textsc{m}}$,$O(h^4)$,$\Sp^4/\Sp^2_{\textsc{m}}$,$O(h^5)$,$\Sp^5/\Sp^{3}_{\textsc{m}}$,$O(h^6)$}]
				\addplot table [x=h, y=L2P2, col sep=comma] {\DATAFILE};
				\addplot table [x=h, y=h3L2, col sep=comma] {\DATAFILE};
				\addplot table [x=h, y=L2P3, col sep=comma] {\DATAFILE};
				\addplot table [x=h, y=h4L2, col sep=comma] {\DATAFILE};
				\addplot table [x=h, y=L2P4, col sep=comma, restrict expr to domain={\coordindex}{0:3}] {\DATAFILE};
				\addplot table [x=h, y expr=0.01*\thisrow{h5L2}, col sep=comma] {\DATAFILE};
				\addplot table [x=h, y=L2P5, col sep=comma, restrict expr to domain={\coordindex}{0:2}] {\DATAFILE};
				\addplot table [x=h, y expr=0.00001*\thisrow{h6L2}, col sep=comma] {\DATAFILE};
			\end{loglogaxis}
		\end{tikzpicture}
		\label{fig:Square_L2_C2}
	\end{subfigure}
\hfill
	\begin{subfigure}[t]{0.45\linewidth}
		\centering
		\begin{tikzpicture}[font=\small, trim axis right]
			\begin{loglogaxis}[
				cycle list name=Lista1,
				xmajorgrids=true,
				ymajorgrids=true,
				width=\linewidth,
				height=1.25\linewidth,
				xlabel={$h$},
				xminorticks=false,
				yminorticks=false,
				ymin={1e-18},
				ylabel={$\nicefrac{\|u-u_h\|_{L^{\infty}}}{\|u\|_{L^{\infty}}}$},
				legend columns=2,
				legend pos=south east,
				legend style={at={(0.99,0.01)},anchor=south east},
				legend entries={$\Sp^2/\Sp^0_{\textsc{m}}$,$O(h^3)$,$\Sp^3/\Sp^1_{\textsc{m}}$,$O(h^4)$,$\Sp^4/\Sp^2_{\textsc{m}}$,$O(h^5)$,$\Sp^5/\Sp^{3}_{\textsc{m}}$,$O(h^6)$}]
				\addplot table [x=h, y=LinfP2, col sep=comma] {\DATAFILE};
				\addplot table [x=h, y=h3Linf, col sep=comma] {\DATAFILE};
				\addplot table [x=h, y=LinfP3, col sep=comma] {\DATAFILE};
				\addplot table [x=h, y=h4Linf, col sep=comma] {\DATAFILE};
				\addplot table [x=h, y=LinfP4, col sep=comma, restrict expr to domain={\coordindex}{0:3}] {\DATAFILE};
				\addplot table [x=h, y expr=0.025*\thisrow{h5Linf}, col sep=comma] {\DATAFILE};
				\addplot table [x=h, y=LinfP5, col sep=comma, restrict expr to domain={\coordindex}{0:2}] {\DATAFILE};
				\addplot table [x=h, y expr=0.00001*\thisrow{h6Linf}, col sep=comma] {\DATAFILE};
			\end{loglogaxis}
		\end{tikzpicture}
		\label{fig:Square_LInf_C2}
	\end{subfigure}
	\caption{Relative error for the square domain with $C^2$-constraints at the vertices.}
	\label{fig:Square_C2}
\end{figure}
\tikzstyle{Linea1}=[thick,mark=*]
\tikzstyle{Linea2}=[thick,dashed]


\pgfplotscreateplotcyclelist{Lista1}{%
	{Linea1,Red},
	{Linea2,Red},
	{Linea1,Green},
	{Linea2,Green},
	{Linea1,Cyan},
	{Linea2,Cyan},
	{Linea1,Violet},
	{Linea2,Violet}}

\def \DATAFILE {ErrAirfoil.csv}

\begin{figure}[h]
	\centering
	\hspace*{\fill}
	\begin{subfigure}[t]{0.45\linewidth}
		\centering
		\begin{tikzpicture}[font=\small, trim axis left]
			\begin{loglogaxis}[
				cycle list name=Lista1,
				xmajorgrids=true,
				ymajorgrids=true,
				width=\linewidth,
				height=1.5\linewidth,
				xlabel={$h$},
				xminorticks=false,
				yminorticks=false,
				ymin=1e-7,
				ylabel={$\nicefrac{\|u-u_h\|_{\mathcal{H}^2}}{\|u\|_{\mathcal{H}^2}}$},
				legend columns=2,
				legend pos=south east,
				legend style={at={(0.99,0.01)},anchor=south east},
				legend entries={$\Sp^2/\Sp^0_{\textsc{m}}$,$O(h)$,$\Sp^3/\Sp^1_{\textsc{m}}$,$O(h^2)$,$\Sp^4/\Sp^2_{\textsc{m}}$,$O(h^3)$,$\Sp^5/\Sp^{3}_{\textsc{m}}$,$O(h^4)$}]
				\addplot table [x=h, y=H2P2, col sep=comma] {\DATAFILE};
				\addplot table [x=h, y=h1H2, col sep=comma] {\DATAFILE};
				\addplot table [x=h, y=H2P3, col sep=comma] {\DATAFILE};
				\addplot table [x=h, y=h2H2, col sep=comma] {\DATAFILE};
				\addplot table [x=h, y=H2P4, col sep=comma] {\DATAFILE};
				\addplot table [x=h, y=h3H2, col sep=comma] {\DATAFILE};
				\addplot table [x=h, y=H2P5, col sep=comma] {\DATAFILE};
				\addplot table [x=h, y=h4H2, col sep=comma] {\DATAFILE};
			\end{loglogaxis}
		\end{tikzpicture}
		\label{fig:Airfoil_H2_C2}
	\end{subfigure}
	\hfill
	\begin{subfigure}[t]{0.45\linewidth}
		\centering
		\begin{tikzpicture}[font=\footnotesize, trim axis right]
			\begin{loglogaxis}[
				cycle list name=Lista1,
				xmajorgrids=true,
				ymajorgrids=true,
				width=\linewidth,
				height=1.5\linewidth,
				xlabel={$h$},
				xminorticks=false,
				yminorticks=false,
				ymin=1e-9,
				ylabel={$\nicefrac{\|u-u_h\|_{H^1}}{\|u\|_{H^1}}$},
				legend columns=2,
				legend pos=south east,
				legend style={at={(0.99,0.01)},anchor=south east},
				legend entries={$\Sp^2/\Sp^0_{\textsc{m}}$,$O(h^2)$,$\Sp^3/\Sp^1_{\textsc{m}}$,$O(h^3)$,$\Sp^4/\Sp^2_{\textsc{m}}$,$O(h^4)$,$\Sp^5/\Sp^{3}_{\textsc{m}}$,$O(h^5)$}]
				\addplot table [x=h, y=H1P2, col sep=comma] {\DATAFILE};
				\addplot table [x=h, y=h2H1, col sep=comma] {\DATAFILE};
				\addplot table [x=h, y=H1P3, col sep=comma] {\DATAFILE};
				\addplot table [x=h, y=h3H1, col sep=comma] {\DATAFILE};
				\addplot table [x=h, y=H1P4, col sep=comma] {\DATAFILE};
				\addplot table [x=h, y=h4H1, col sep=comma] {\DATAFILE};
				\addplot table [x=h, y=H1P5, col sep=comma] {\DATAFILE};
				\addplot table [x=h, y=h5H1, col sep=comma] {\DATAFILE};
			\end{loglogaxis}
		\end{tikzpicture}
		\label{fig:Airfoil_H1_C2}
	\end{subfigure}
	\\
	\hspace*{\fill}
	\begin{subfigure}[t]{0.45\linewidth}
		\centering
		\begin{tikzpicture}[font=\small, trim axis left]
			\begin{loglogaxis}[
				cycle list name=Lista1,
				xmajorgrids=true,
				ymajorgrids=true,
				width=\linewidth,
				height=1.5\linewidth,
				xlabel={$h$},
				xminorticks=false,
				yminorticks=false,
				ymin=1e-11,
				ylabel={$\nicefrac{\|u-u_h\|_{L^2}}{\|u\|_{L^2}}$},
				legend columns=2,
				legend pos=south east,
				legend style={at={(0.99,0.01)},anchor=south east},
				legend entries={$\Sp^2/\Sp^0_{\textsc{m}}$,$O(h^3)$,$\Sp^3/\Sp^1_{\textsc{m}}$,$O(h^4)$,$\Sp^4/\Sp^2_{\textsc{m}}$,$O(h^5)$,$\Sp^5/\Sp^{3}_{\textsc{m}}$,$O(h^6)$}]
				\addplot table [x=h, y=L2P2, col sep=comma] {\DATAFILE};
				\addplot table [x=h, y=h3L2, col sep=comma] {\DATAFILE};
				\addplot table [x=h, y=L2P3, col sep=comma] {\DATAFILE};
				\addplot table [x=h, y=h4L2, col sep=comma] {\DATAFILE};
				\addplot table [x=h, y=L2P4, col sep=comma] {\DATAFILE};
				\addplot table [x=h, y=h5L2, col sep=comma] {\DATAFILE};
				\addplot table [x=h, y=L2P5, col sep=comma] {\DATAFILE};
				\addplot table [x=h, y=h6L2, col sep=comma] {\DATAFILE};
			\end{loglogaxis}
		\end{tikzpicture}
		\label{fig:Airfoil_L2_C2}
	\end{subfigure}
	\hfill
	\begin{subfigure}[t]{0.45\linewidth}
		\centering
		\begin{tikzpicture}[font=\small, trim axis right]
			\begin{loglogaxis}[
				cycle list name=Lista1,
				xmajorgrids=true,
				ymajorgrids=true,
				width=\linewidth,
				height=1.5\linewidth,
				xlabel={$h$},
				xminorticks=false,
				yminorticks=false,
				ymin=1e-10,
				ylabel={$\nicefrac{\|u-u_h\|_{L^{\infty}}}{\|u\|_{L^{\infty}}}$},
				legend columns=2,
				legend pos=south east,
				legend style={at={(0.99,0.01)},anchor=south east},
				legend entries={$\Sp^2/\Sp^0_{\textsc{m}}$,$O(h^3)$,$\Sp^3/\Sp^1_{\textsc{m}}$,$O(h^4)$,$\Sp^4/\Sp^2_{\textsc{m}}$,$O(h^5)$,$\Sp^5/\Sp^{3}_{\textsc{m}}$,$O(h^6)$}]
				\addplot table [x=h, y=LinfP2, col sep=comma] {\DATAFILE};
				\addplot table [x=h, y=h3Linf, col sep=comma] {\DATAFILE};
				\addplot table [x=h, y=LinfP3, col sep=comma] {\DATAFILE};
				\addplot table [x=h, y=h4Linf, col sep=comma] {\DATAFILE};
				\addplot table [x=h, y=LinfP4, col sep=comma] {\DATAFILE};
				\addplot table [x=h, y=h5Linf, col sep=comma] {\DATAFILE};
				\addplot table [x=h, y=LinfP5, col sep=comma] {\DATAFILE};
				\addplot table [x=h, y=h6Linf, col sep=comma] {\DATAFILE};
			\end{loglogaxis}
		\end{tikzpicture}
		\label{fig:Airfoil_LInf_C2}
	\end{subfigure}
	\caption{Relative error for the airfoil domain with $C^2$-constraints at the vertices.}
	\label{fig:Airfoil_C2}
\end{figure}
\tikzstyle{Linea1}=[thick,mark=*]
\tikzstyle{Linea2}=[thick,dashed]


\pgfplotscreateplotcyclelist{Lista1}{%
	{Linea1,Red},
	{Linea2,Red},
	{Linea1,Green},
	{Linea2,Green},
	{Linea1,Cyan},
	{Linea2,Cyan},
	{Linea1,Violet},
	{Linea2,Violet}}

\def \DATAFILE {ErrQuartercircle.csv}

\begin{figure}[h]
	\centering
	\hspace*{\fill}
	\begin{subfigure}[t]{0.45\linewidth}
		\centering
		\begin{tikzpicture}[font=\footnotesize, trim axis left]
			\begin{loglogaxis}[
				cycle list name=Lista1,
				xmajorgrids=true,
				ymajorgrids=true,
				width=\linewidth,
				height=1.5\linewidth,
				xlabel={$h$},
				xminorticks=false,
				yminorticks=false,
				ymin=1e-11,
				ylabel={$\nicefrac{\|u-u_h\|_{\mathcal{H}^2}}{\|u\|_{\mathcal{H}^2}}$},
				legend columns=2,
				legend pos=south east,
				legend style={at={(0.99,0.01)},anchor=south east},
				legend entries={$\Sp^2/\Sp^0_{\textsc{m}}$,$O(h)$,$\Sp^3/\Sp^1_{\textsc{m}}$,$O(h^2)$,$\Sp^4/\Sp^2_{\textsc{m}}$,$O(h^3)$,$\Sp^5/\Sp^{3}_{\textsc{m}}$,$O(h^4)$}]
				\addplot table [x=h, y=H2P2, col sep=comma] {\DATAFILE};
				\addplot table [x=h, y=h1H2, col sep=comma] {\DATAFILE};
				\addplot table [x=h, y=H2P3, col sep=comma] {\DATAFILE};
				\addplot table [x=h, y=h2H2, col sep=comma] {\DATAFILE};
				\addplot table [x=h, y=H2P4, col sep=comma] {\DATAFILE};
				\addplot table [x=h, y=h3H2, col sep=comma] {\DATAFILE};
				\addplot table [x=h, y=H2P5, col sep=comma] {\DATAFILE};
				\addplot table [x=h, y=h4H2, col sep=comma] {\DATAFILE};
			\end{loglogaxis}
		\end{tikzpicture}
		\label{fig:Quartercircle_H2_C2}
	\end{subfigure}
	\hfill
	\begin{subfigure}[t]{0.45\linewidth}
		\centering
		\begin{tikzpicture}[font=\small, trim axis right]
			\begin{loglogaxis}[
				cycle list name=Lista1,
				xmajorgrids=true,
				ymajorgrids=true,
				width=\linewidth,
				height=1.5\linewidth,
				xlabel={$h$},
				xminorticks=false,
				yminorticks=false,
				ymin=1e-14,
				ylabel={$\nicefrac{\|u-u_h\|_{H^1}}{\|u\|_{H^1}}$},
				legend columns=2,
				legend pos=south east,
				legend style={at={(0.99,0.01)},anchor=south east},
				legend entries={$\Sp^2/\Sp^0_{\textsc{m}}$,$O(h^2)$,$\Sp^3/\Sp^1_{\textsc{m}}$,$O(h^3)$,$\Sp^4/\Sp^2_{\textsc{m}}$,$O(h^4)$,$\Sp^5/\Sp^{3}_{\textsc{m}}$,$O(h^5)$}]
				\addplot table [x=h, y=H1P2, col sep=comma] {\DATAFILE};
				\addplot table [x=h, y=h2H1, col sep=comma] {\DATAFILE};
				\addplot table [x=h, y=H1P3, col sep=comma] {\DATAFILE};
				\addplot table [x=h, y=h3H1, col sep=comma] {\DATAFILE};
				\addplot table [x=h, y=H1P4, col sep=comma] {\DATAFILE};
				\addplot table [x=h, y=h4H1, col sep=comma] {\DATAFILE};
				\addplot table [x=h, y=H1P5, col sep=comma] {\DATAFILE};
				\addplot table [x=h, y expr=0.1*\thisrow{h5H1}, col sep=comma] {\DATAFILE};
			\end{loglogaxis}
		\end{tikzpicture}
		\label{fig:Quartercircle_H1_C2}
	\end{subfigure}
\\
\hspace*{\fill}
	\begin{subfigure}[t]{0.45\linewidth}
		\centering
		\begin{tikzpicture}[font=\small, trim axis left]
			\begin{loglogaxis}[
				cycle list name=Lista1,
				xmajorgrids=true,
				ymajorgrids=true,
				width=\linewidth,
				height=1.5\linewidth,
				xlabel={$h$},
				xminorticks=false,
				yminorticks=false,
				ymin=1e-18,
				ylabel={$\nicefrac{\|u-u_h\|_{L^2}}{\|u\|_{L^2}}$},
				legend columns=2,
				legend pos=south east,
				legend style={at={(0.99,0.01)},anchor=south east},
				legend entries={$\Sp^2/\Sp^0_{\textsc{m}}$,$O(h^3)$,$\Sp^3/\Sp^1_{\textsc{m}}$,$O(h^4)$,$\Sp^4/\Sp^2_{\textsc{m}}$,$O(h^5)$,$\Sp^5/\Sp^{3}_{\textsc{m}}$,$O(h^6)$}]
				\addplot table [x=h, y=L2P2, col sep=comma] {\DATAFILE};
				\addplot table [x=h, y=h3L2, col sep=comma] {\DATAFILE};
				\addplot table [x=h, y=L2P3, col sep=comma] {\DATAFILE};
				\addplot table [x=h, y=h4L2, col sep=comma] {\DATAFILE};
				\addplot table [x=h, y=L2P4, col sep=comma] {\DATAFILE};
				\addplot table [x=h, y=h5L2, col sep=comma] {\DATAFILE};
				\addplot table [x=h, y=L2P5, col sep=comma, restrict expr to domain={\coordindex}{0:3}] {\DATAFILE};
				\addplot table [x=h, y expr=0.001*\thisrow{h6L2}, col sep=comma] {\DATAFILE};
			\end{loglogaxis}
		\end{tikzpicture}
		\label{fig:Quartercircle_L2_C2}
	\end{subfigure}
	\hfill
	\begin{subfigure}[t]{0.45\linewidth}
		\centering
		\begin{tikzpicture}[font=\small, trim axis right]
			\begin{loglogaxis}[
				cycle list name=Lista1,
				xmajorgrids=true,
				ymajorgrids=true,
				width=\linewidth,
				height=1.5\linewidth,
				xlabel={$h$},
				xminorticks=false,
				yminorticks=false,
				ymin=1e-17,
				ylabel={$\nicefrac{\|u-u_h\|_{L^{\infty}}}{\|u\|_{L^{\infty}}}$},
				legend columns=2,
				legend pos=south east,
				legend style={at={(0.99,0.01)},anchor=south east},
				legend entries={$\Sp^2/\Sp^0_{\textsc{m}}$,$O(h^3)$,$\Sp^3/\Sp^1_{\textsc{m}}$,$O(h^4)$,$\Sp^4/\Sp^2_{\textsc{m}}$,$O(h^5)$,$\Sp^5/\Sp^{3}_{\textsc{m}}$,$O(h^6)$}]
				\addplot table [x=h, y=LinfP2, col sep=comma] {\DATAFILE};
				\addplot table [x=h, y=h3Linf, col sep=comma] {\DATAFILE};
				\addplot table [x=h, y=LinfP3, col sep=comma] {\DATAFILE};
				\addplot table [x=h, y=h4Linf, col sep=comma] {\DATAFILE};
				\addplot table [x=h, y=LinfP4, col sep=comma] {\DATAFILE};
				\addplot table [x=h, y=h5Linf, col sep=comma] {\DATAFILE};
				\addplot table [x=h, y=LinfP5, col sep=comma, restrict expr to domain={\coordindex}{0:3}] {\DATAFILE};
				\addplot table [x=h, y expr=0.01*\thisrow{h6Linf}, col sep=comma] {\DATAFILE};
			\end{loglogaxis}
		\end{tikzpicture}
		\label{fig:Quartercircle_LInf_C2}
	\end{subfigure}
	\caption{Relative error for the quarter circle domain with $C^2$-constraints at the vertices.}
	\label{fig:Quartercircle_C2}
\end{figure}
\tikzstyle{Linea1}=[thick,mark=*]
\tikzstyle{Linea2}=[thick,dashed]


\pgfplotscreateplotcyclelist{Lista1}{%
	{Linea1,Red},
	{Linea2,Red},
	{Linea1,Green},
	{Linea2,Green},
	{Linea1,Cyan},
	{Linea2,Cyan},
	{Linea1,Violet},
	{Linea2,Violet}}

\def \DATAFILE {ErrQuartercircleR1.csv}

\begin{figure}[h]
	\centering
	\hspace*{\fill}
	\begin{subfigure}[t]{0.45\linewidth}
		\centering
		\begin{tikzpicture}[font=\small, trim axis left]
			\begin{loglogaxis}[
				cycle list name=Lista1,
				xmajorgrids=true,
				ymajorgrids=true,
				width=\linewidth,
				height=1.5\linewidth,
				xlabel={$h$},
				xminorticks=false,
				yminorticks=false,
				ymin={1e-13},
				ylabel={$\nicefrac{\|u-u_h\|_{\mathcal{H}^2}}{\|u\|_{\mathcal{H}^2}}$},
				legend columns=2,
				legend pos=south east,
				legend style={at={(0.99,0.01)},anchor=south east},
				legend entries={$\Sp^2/\Sp^0_{\textsc{m}}$,$O(h)$,$\Sp^3/\Sp^1_{\textsc{m}}$,$O(h^2)$,$\Sp^4/\Sp^2_{\textsc{m}}$,$O(h^3)$,$\Sp^5/\Sp^{3}_{\textsc{m}}$,$O(h^4)$}]
				\addplot table [x=h, y=H2P2, col sep=comma] {\DATAFILE};
				\addplot table [x=h, y=h1H2, col sep=comma] {\DATAFILE};
				\addplot table [x=h, y=H2P3, col sep=comma] {\DATAFILE};
				\addplot table [x=h, y=h2H2, col sep=comma] {\DATAFILE};
				\addplot table [x=h, y=H2P4, col sep=comma] {\DATAFILE};
				\addplot table [x=h, y=h3H2, col sep=comma] {\DATAFILE};
				\addplot table [x=h, y=H2P5, col sep=comma, restrict expr to domain={\coordindex}{0:3}] {\DATAFILE};
				\addplot table [x=h, y expr=0.1*\thisrow{h4H2}, col sep=comma] {\DATAFILE};
			\end{loglogaxis}
		\end{tikzpicture}
		\label{fig:Quartercircle_H2_C2R1}
	\end{subfigure}
	\hfill
	\begin{subfigure}[t]{0.45\linewidth}
		\centering
		\begin{tikzpicture}[font=\small, trim axis right]
			\begin{loglogaxis}[
				cycle list name=Lista1,
				xmajorgrids=true,
				ymajorgrids=true,
				width=\linewidth,
				height=1.5\linewidth,
				xlabel={$h$},
				xminorticks=false,
				yminorticks=false,
				ymin={1e-15},
				ylabel={$\nicefrac{\|u-u_h\|_{H^1}}{\|u\|_{H^1}}$},
				legend columns=2,
				legend pos=south east,
				legend style={at={(0.99,0.01)},anchor=south east},
				legend entries={$\Sp^2/\Sp^0_{\textsc{m}}$,$O(h^2)$,$\Sp^3/\Sp^1_{\textsc{m}}$,$O(h^3)$,$\Sp^4/\Sp^2_{\textsc{m}}$,$O(h^4)$,$\Sp^5/\Sp^{3}_{\textsc{m}}$,$O(h^5)$}]
				\addplot table [x=h, y=H1P2, col sep=comma] {\DATAFILE};
				\addplot table [x=h, y=h2H1, col sep=comma] {\DATAFILE};
				\addplot table [x=h, y=H1P3, col sep=comma] {\DATAFILE};
				\addplot table [x=h, y=h3H1, col sep=comma] {\DATAFILE};
				\addplot table [x=h, y=H1P4, col sep=comma] {\DATAFILE};
				\addplot table [x=h, y=h4H1, col sep=comma] {\DATAFILE};
				\addplot table [x=h, y=H1P5, col sep=comma, restrict expr to domain={\coordindex}{0:3}] {\DATAFILE};
				\addplot table [x=h, y expr=0.001*\thisrow{h5H1}, col sep=comma] {\DATAFILE};
			\end{loglogaxis}
		\end{tikzpicture}
		\label{fig:Quartercircle_H1_C2R1}
	\end{subfigure}
\\
\hspace*{\fill}
	\begin{subfigure}[t]{0.4\linewidth}
		\centering
		\begin{tikzpicture}[font=\small, trim axis left]
			\begin{loglogaxis}[
				cycle list name=Lista1,
				xmajorgrids=true,
				ymajorgrids=true,
				width=\linewidth,
				height=1.5\linewidth,
				xlabel={$h$},
				xminorticks=false,
				yminorticks=false,
				ymin=1e-18,
				ylabel={$\nicefrac{\|u-u_h\|_{L^2}}{\|u\|_{L^2}}$},
				legend columns=2,
				legend pos=south east,
				legend style={at={(0.99,0.01)},anchor=south east},
				legend entries={$\Sp^2/\Sp^0_{\textsc{m}}$,$O(h^3)$,$\Sp^3/\Sp^1_{\textsc{m}}$,$O(h^4)$,$\Sp^4/\Sp^2_{\textsc{m}}$,$O(h^5)$,$\Sp^5/\Sp^{3}_{\textsc{m}}$,$O(h^6)$}]
				\addplot table [x=h, y=L2P2, col sep=comma] {\DATAFILE};
				\addplot table [x=h, y=h3L2, col sep=comma] {\DATAFILE};
				\addplot table [x=h, y=L2P3, col sep=comma] {\DATAFILE};
				\addplot table [x=h, y=h4L2, col sep=comma] {\DATAFILE};
				\addplot table [x=h, y=L2P4, col sep=comma, restrict expr to domain={\coordindex}{0:3}] {\DATAFILE};
				\addplot table [x=h, y expr=0.005*\thisrow{h5L2}, col sep=comma] {\DATAFILE};
				\addplot table [x=h, y=L2P5, col sep=comma, restrict expr to domain={\coordindex}{0:2}] {\DATAFILE};
				\addplot table [x=h, y expr=0.0000025*\thisrow{h6L2}, col sep=comma] {\DATAFILE};
			\end{loglogaxis}
		\end{tikzpicture}
		\label{fig:Quartercircle_L2_C2R1}
	\end{subfigure}
	\hfill
	\begin{subfigure}[t]{0.4\linewidth}
		\centering
		\begin{tikzpicture}[font=\small, trim axis right]
			\begin{loglogaxis}[
				cycle list name=Lista1,
				xmajorgrids=true,
				ymajorgrids=true,
				width=\linewidth,
				height=1.5\linewidth,
				xlabel={$h$},
				xminorticks=false,
				yminorticks=false,
				ymin=1e-17,
				ylabel={$\nicefrac{\|u-u_h\|_{L^{\infty}}}{\|u\|_{L^{\infty}}}$},
				legend columns=2,
				legend pos=south east,
				legend style={at={(0.99,0.01)},anchor=south east},
				legend entries={$\Sp^2/\Sp^0_{\textsc{m}}$,$O(h^3)$,$\Sp^3/\Sp^1_{\textsc{m}}$,$O(h^4)$,$\Sp^4/\Sp^2_{\textsc{m}}$,$O(h^5)$,$\Sp^5/\Sp^{3}_{\textsc{m}}$,$O(h^6)$}]
				\addplot table [x=h, y=LinfP2, col sep=comma] {\DATAFILE};
				\addplot table [x=h, y=h3Linf, col sep=comma] {\DATAFILE};
				\addplot table [x=h, y=LinfP3, col sep=comma] {\DATAFILE};
				\addplot table [x=h, y=h4Linf, col sep=comma] {\DATAFILE};
				\addplot table [x=h, y=LinfP4, col sep=comma, restrict expr to domain={\coordindex}{0:3}] {\DATAFILE};
				\addplot table [x=h, y expr=0.005*\thisrow{h5Linf}, col sep=comma] {\DATAFILE};
				\addplot table [x=h, y=LinfP5, col sep=comma, restrict expr to domain={\coordindex}{0:2}] {\DATAFILE};
				\addplot table [x=h, y expr=0.000025*\thisrow{h6Linf}, col sep=comma] {\DATAFILE};
			\end{loglogaxis}
		\end{tikzpicture}
		\label{fig:Quartercircle_LInf__C2R1}
	\end{subfigure}
	\caption{Relative error for the quarter circle domain with $C^2$-constraints at the vertices and splines globally $C^1$-continuous.}
	\label{fig:Quartercircle_C2R1}
\end{figure}

\subsection{{ Tests without} $C^2$-continuity at the vertices}
Additionally, we tested the pairing $\Sp^p/\Sp^{p-2}_{\textsc{m}}$ with splines of maximal global regularity, $r = p - 1$, for the biharmonic problem on the quarter-circle domain, without imposing strong $C^2$-continuity at the vertex. As before, Figure~\ref{fig:Quartercircle_C0} presents the convergence plots in $\mathcal{H}^2$, $H^1$, $L^2$, and $L^{\infty}$ norms for primal spaces with polynomial degrees $p = 2, \dots, 5$. Numerical tests show that the pairing $\Sp^p/\Sp^{p-2}_{\textsc{m}}$ achieves the optimal rate of convergence even in this case. Figure~\ref{fig:Quartercircle_solutions} demonstrates that the difference between the numerical solutions obtained with and without $C^2$ continuity is localized in a neighborhood around the vertex.
\tikzstyle{Linea1}=[thick,mark=*]
\tikzstyle{Linea2}=[thick,dashed]


\pgfplotscreateplotcyclelist{Lista1}{%
	{Linea1,Red},
	{Linea2,Red},
	{Linea1,Green},
	{Linea2,Green},
	{Linea1,Cyan},
	{Linea2,Cyan},
	{Linea1,Violet},
	{Linea2,Violet}}

\def \DATAFILE {ErrQuartercircleRmC0.csv}

\begin{figure}[h]
	\centering
	\hspace*{\fill}
	\begin{subfigure}[t]{0.45\linewidth}
		\centering
		\begin{tikzpicture}[font=\small, trim axis left]
			\begin{loglogaxis}[
				cycle list name=Lista1,
				xmajorgrids=true,
				ymajorgrids=true,
				width=\linewidth,
				height=1.5\linewidth,
				xlabel={$h$},
				xminorticks=false,
				yminorticks=false,
				ymin=1e-12,
				ylabel={$\nicefrac{\|u-u_h\|_{\mathcal{H}^2}}{\|u\|_{\mathcal{H}^2}}$},
				legend columns=2,
				legend pos=south east,
				legend style={at={(0.99,0.01)},anchor=south east},
				legend entries={$\Sp^2/\Sp^0_{\textsc{m}}$,$O(h)$,$\Sp^3/\Sp^1_{\textsc{m}}$,$O(h^2)$,$\Sp^4/\Sp^2_{\textsc{m}}$,$O(h^3)$,$\Sp^5/\Sp^{3}_{\textsc{m}}$,$O(h^4)$}]
				\addplot table [x=h, y=H2P2, col sep=comma] {\DATAFILE};
				\addplot table [x=h, y=h1H2, col sep=comma] {\DATAFILE};
				\addplot table [x=h, y=H2P3, col sep=comma] {\DATAFILE};
				\addplot table [x=h, y=h2H2, col sep=comma] {\DATAFILE};
				\addplot table [x=h, y=H2P4, col sep=comma] {\DATAFILE};
				\addplot table [x=h, y=h3H2, col sep=comma] {\DATAFILE};
				\addplot table [x=h, y=H2P5, col sep=comma] {\DATAFILE};
				\addplot table [x=h, y=h4H2, col sep=comma] {\DATAFILE};
			\end{loglogaxis}
		\end{tikzpicture}
		\label{fig:Quartercircle_H2_C0}
	\end{subfigure}
	\hfill
	\begin{subfigure}[t]{0.45\linewidth}
		\centering
		\begin{tikzpicture}[font=\small, trim axis right]
			\begin{loglogaxis}[
				cycle list name=Lista1,
				xmajorgrids=true,
				ymajorgrids=true,
				width=\linewidth,
				height=1.5\linewidth,
				xlabel={$h$},
				xminorticks=false,
				yminorticks=false,
				ymin=1e-14,
				ylabel={$\nicefrac{\|u-u_h\|_{H^1}}{\|u\|_{H^1}}$},
				legend columns=2,
				legend pos=south east,
				legend style={at={(0.99,0.01)},anchor=south east},
				legend entries={$\Sp^2/\Sp^0_{\textsc{m}}$,$O(h^2)$,$\Sp^3/\Sp^1_{\textsc{m}}$,$O(h^3)$,$\Sp^4/\Sp^2_{\textsc{m}}$,$O(h^4)$,$\Sp^5/\Sp^{3}_{\textsc{m}}$,$O(h^5)$}]
				\addplot table [x=h, y=H1P2, col sep=comma] {\DATAFILE};
				\addplot table [x=h, y=h2H1, col sep=comma] {\DATAFILE};
				\addplot table [x=h, y=H1P3, col sep=comma] {\DATAFILE};
				\addplot table [x=h, y=h3H1, col sep=comma] {\DATAFILE};
				\addplot table [x=h, y=H1P4, col sep=comma] {\DATAFILE};
				\addplot table [x=h, y=h4H1, col sep=comma] {\DATAFILE};
				\addplot table [x=h, y=H1P5, col sep=comma] {\DATAFILE};
				\addplot table [x=h, y=h5H1, col sep=comma] {\DATAFILE};
			\end{loglogaxis}
		\end{tikzpicture}
		\label{fig:Quartercircle_H1_C0}
	\end{subfigure}
\\
\hspace*{\fill}
	\begin{subfigure}[t]{0.45\linewidth}
		\centering
		\begin{tikzpicture}[font=\small, trim axis left]
			\begin{loglogaxis}[
				cycle list name=Lista1,
				xmajorgrids=true,
				ymajorgrids=true,
				width=\linewidth,
				height=1.5\linewidth,
				xlabel={$h$},
				xminorticks=false,
				yminorticks=false,
				ymin=1e-18,
				ylabel={$\nicefrac{\|u-u_h\|_{L^2}}{\|u\|_{L^2}}$},
				legend columns=2,
				legend pos=south east,
				legend style={at={(0.99,0.01)},anchor=south east},
				legend entries={$\Sp^2/\Sp^0_{\textsc{m}}$,$O(h^3)$,$\Sp^3/\Sp^1_{\textsc{m}}$,$O(h^4)$,$\Sp^4/\Sp^2_{\textsc{m}}$,$O(h^5)$,$\Sp^5/\Sp^{3}_{\textsc{m}}$,$O(h^6)$}]
				\addplot table [x=h, y=L2P2, col sep=comma] {\DATAFILE};
				\addplot table [x=h, y=h3L2, col sep=comma] {\DATAFILE};
				\addplot table [x=h, y=L2P3, col sep=comma] {\DATAFILE};
				\addplot table [x=h, y=h4L2, col sep=comma] {\DATAFILE};
				\addplot table [x=h, y=L2P4, col sep=comma] {\DATAFILE};
				\addplot table [x=h, y=h5L2, col sep=comma] {\DATAFILE};
				\addplot table [x=h, y=L2P5, col sep=comma, restrict expr to domain={\coordindex}{0:3}] {\DATAFILE};
				\addplot table [x=h, y expr=0.0025*\thisrow{h6L2}, col sep=comma] {\DATAFILE};
			\end{loglogaxis}
		\end{tikzpicture}
		\label{fig:Quartercircle_L2_C0}
	\end{subfigure}
	\hfill
	\begin{subfigure}[t]{0.45\linewidth}
		\centering
		\begin{tikzpicture}[font=\small, trim axis right]
			\begin{loglogaxis}[
				cycle list name=Lista1,
				xmajorgrids=true,
				ymajorgrids=true,
				width=\linewidth,
				height=1.5\linewidth,
				xlabel={$h$},
				xminorticks=false,
				yminorticks=false,
				ymin=1e-17,
				ylabel={$\nicefrac{\|u-u_h\|_{L^{\infty}}}{\|u\|_{L^{\infty}}}$},
				legend columns=2,
				legend pos=south east,
				legend style={at={(0.99,0.01)},anchor=south east},
				legend entries={$\Sp^2/\Sp^0_{\textsc{m}}$,$O(h^3)$,$\Sp^3/\Sp^1_{\textsc{m}}$,$O(h^4)$,$\Sp^4/\Sp^2_{\textsc{m}}$,$O(h^5)$,$\Sp^5/\Sp^{3}_{\textsc{m}}$,$O(h^6)$}]
				\addplot table [x=h, y=LinfP2, col sep=comma] {\DATAFILE};
				\addplot table [x=h, y=h3Linf, col sep=comma] {\DATAFILE};
				\addplot table [x=h, y=LinfP3, col sep=comma] {\DATAFILE};
				\addplot table [x=h, y=h4Linf, col sep=comma] {\DATAFILE};
				\addplot table [x=h, y=LinfP4, col sep=comma] {\DATAFILE};
				\addplot table [x=h, y=h5Linf, col sep=comma] {\DATAFILE};
				\addplot table [x=h, y=LinfP5, col sep=comma, restrict expr to domain={\coordindex}{0:3}] {\DATAFILE};
				\addplot table [x=h, y expr=0.01*\thisrow{h6Linf}, col sep=comma] {\DATAFILE};
			\end{loglogaxis}
		\end{tikzpicture}
		\label{fig:Quartercircle_LInf_C0}
	\end{subfigure}
	\hspace*{\fill}
	\caption{Relative error for the quarter circle domain with $C^0$-constraints at the vertices.}
	\label{fig:Quartercircle_C0}
\end{figure}
\begin{figure}[h]
	\centering
	\includegraphics[width=\linewidth]{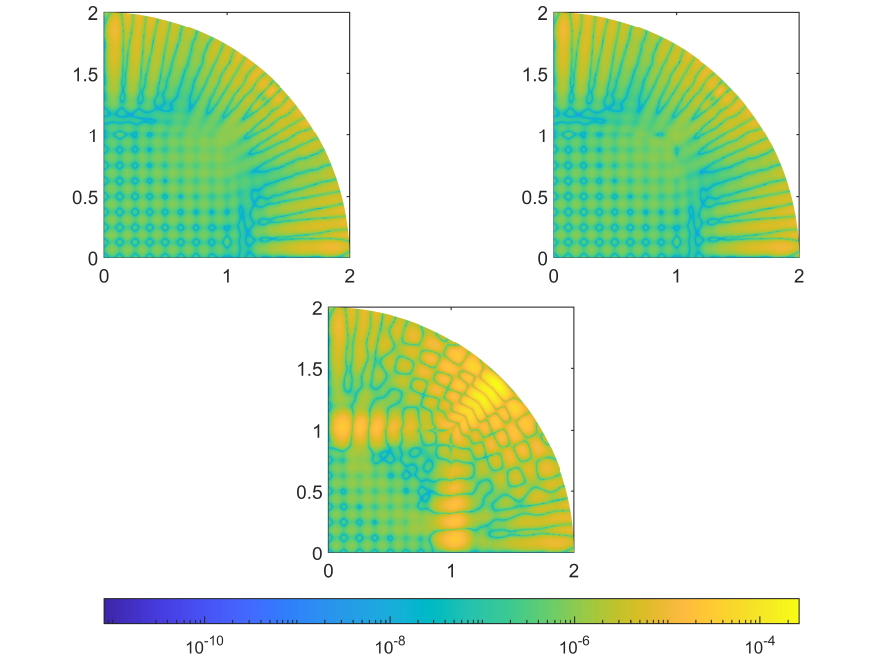}
	\caption{Error for $h \approx 0.155$ on the quarter circle domain  with: $C^2$-continuity at the vertices and $\Sp^3_2/\Sp^{1}_{0,\textsc{m}}$ pairing (top left), $C^0$-continuity at the vertices and $\Sp^3_2/\Sp^{1}_{0,\textsc{m}}$ pairing (top right) and $C^0$-continuity at the vertices and $\Sp^3_2/\Sp^{1}_0$ pairing (bottom).}
	\label{fig:Quartercircle_solutions}
\end{figure}

\subsection{{ Further} tests}
In the context of numerical experiments where strong $C^2$-continuity at the vertices is not imposed, we have tested the bilaplacian problem for the pairing $\Sp^p / \Sp^{p-2}$, i.e., without knot vector reduction for the Lagrange multiplier spaces. Although these tests extend beyond the theoretical framework, the numerical orders of convergence are optimal (see Figure~\ref{fig:QuartercircleKNRC0}). However, as shown in Figure~\ref{fig:Quartercircle_solutions}, this choice leads to larger errors, particularly along $\Sigma$.
\tikzstyle{Linea1}=[thick,mark=*]
\tikzstyle{Linea2}=[thick,dashed]


\pgfplotscreateplotcyclelist{Lista1}{%
	{Linea1,Red},
	{Linea2,Red},
	{Linea1,Green},
	{Linea2,Green},
	{Linea1,Cyan},
	{Linea2,Cyan},
	{Linea1,Violet},
	{Linea2,Violet}}

\def \DATAFILE {ErrQuartercircleNRmC0.csv}

\begin{figure}[h]
	\centering
	\hspace*{\fill}
	\begin{subfigure}[t]{0.45\linewidth}
		\centering
		\begin{tikzpicture}[font=\small, trim axis left]
			\begin{loglogaxis}[
				cycle list name=Lista1,
				xmajorgrids=true,
				ymajorgrids=true,
				width=\linewidth,
				height=1.5\linewidth,
				xlabel={$h$},
				xminorticks=false,
				yminorticks=false,
				ymin=1e-9,
				ylabel={$\nicefrac{\|u-u_h\|_{\mathcal{H}^2}}{\|u\|_{\mathcal{H}^2}}$},
				legend columns=2,
				legend pos=south east,
				legend style={at={(0.99,0.01)},anchor=south east},
				legend entries={$\Sp^2/\Sp^0$,$O(h)$,$\Sp^3/\Sp^1$,$O(h^2)$,$\Sp^4/\Sp^2$,$O(h^3)$,$\Sp^5/\Sp^{3}$,$O(h^4)$}]
				\addplot table [x=h, y=H2P2, col sep=comma] {\DATAFILE};
				\addplot table [x=h, y=h1H2, col sep=comma] {\DATAFILE};
				\addplot table [x=h, y=H2P3, col sep=comma] {\DATAFILE};
				\addplot table [x=h, y=h2H2, col sep=comma] {\DATAFILE};
				\addplot table [x=h, y=H2P4, col sep=comma] {\DATAFILE};
				\addplot table [x=h, y=h3H2, col sep=comma] {\DATAFILE};
				\addplot table [x=h, y=H2P5, col sep=comma] {\DATAFILE};
				\addplot table [x=h, y=h4H2, col sep=comma] {\DATAFILE};
			\end{loglogaxis}
		\end{tikzpicture}
		\label{fig:Quartercircle_H2_C0_NME}
	\end{subfigure}
	\hfill
	\begin{subfigure}[t]{0.45\linewidth}
		\centering
		\begin{tikzpicture}[font=\small, trim axis right]
			\begin{loglogaxis}[
				cycle list name=Lista1,
				xmajorgrids=true,
				ymajorgrids=true,
				width=\linewidth,
				height=1.5\linewidth,
				xlabel={$h$},
				xminorticks=false,
				yminorticks=false,
				ymin=1e-12,
				ylabel={$\nicefrac{\|u-u_h\|_{H^1}}{\|u\|_{H^1}}$},
				legend columns=2,
				legend pos=south east,
				legend style={at={(0.99,0.01)},anchor=south east},
				legend entries={$\Sp^2/\Sp^0$,$O(h^2)$,$\Sp^3/\Sp^1$,$O(h^3)$,$\Sp^4/\Sp^2$,$O(h^4)$,$\Sp^5/\Sp^{3}$,$O(h^5)$}]
				\addplot table [x=h, y=H1P2, col sep=comma] {\DATAFILE};
				\addplot table [x=h, y=h2H1, col sep=comma] {\DATAFILE};
				\addplot table [x=h, y=H1P3, col sep=comma] {\DATAFILE};
				\addplot table [x=h, y=h3H1, col sep=comma] {\DATAFILE};
				\addplot table [x=h, y=H1P4, col sep=comma] {\DATAFILE};
				\addplot table [x=h, y=h4H1, col sep=comma] {\DATAFILE};
				\addplot table [x=h, y=H1P5, col sep=comma] {\DATAFILE};
				\addplot table [x=h, y=h5H1, col sep=comma] {\DATAFILE};
			\end{loglogaxis}
		\end{tikzpicture}
		\label{fig:Quartercircle_H1_C0_NME}
	\end{subfigure}
\\
\hspace*{\fill}
	\begin{subfigure}[t]{0.45\linewidth}
		\centering
		\begin{tikzpicture}[font=\small, trim axis left]
			\begin{loglogaxis}[
				cycle list name=Lista1,
				xmajorgrids=true,
				ymajorgrids=true,
				width=\linewidth,
				height=1.5\linewidth,
				xlabel={$h$},
				xminorticks=false,
				yminorticks=false,
				ymin=1e-15,
				ylabel={$\nicefrac{\|u-u_h\|_{L^2}}{\|u\|_{L^2}}$},
				legend columns=2,
				legend pos=south east,
				legend style={at={(0.99,0.01)},anchor=south east},
				legend entries={$\Sp^2/\Sp^0$,$O(h^3)$,$\Sp^3/\Sp^1$,$O(h^4)$,$\Sp^4/\Sp^2$,$O(h^5)$,$\Sp^5/\Sp^{3}$,$O(h^6)$}]
				\addplot table [x=h, y=L2P2, col sep=comma] {\DATAFILE};
				\addplot table [x=h, y=h3L2, col sep=comma] {\DATAFILE};
				\addplot table [x=h, y=L2P3, col sep=comma] {\DATAFILE};
				\addplot table [x=h, y=h4L2, col sep=comma] {\DATAFILE};
				\addplot table [x=h, y=L2P4, col sep=comma] {\DATAFILE};
				\addplot table [x=h, y=h5L2, col sep=comma] {\DATAFILE};
				\addplot table [x=h, y=L2P5, col sep=comma, restrict expr to domain={\coordindex}{0:3}] {\DATAFILE};
				\addplot table [x=h, y expr=0.05*\thisrow{h6L2}, col sep=comma] {\DATAFILE};
			\end{loglogaxis}
		\end{tikzpicture}
		\label{fig:Quartercircle_L2_C0_NME}
	\end{subfigure}
	\hfill
	\begin{subfigure}[t]{0.45\linewidth}
		\centering
		\begin{tikzpicture}[font=\small, trim axis right]
			\begin{loglogaxis}[
				cycle list name=Lista1,
				xmajorgrids=true,
				ymajorgrids=true,
				width=\linewidth,
				height=1.5\linewidth,
				xlabel={$h$},
				xminorticks=false,
				yminorticks=false,
				ymin=1e-13,
				ylabel={$\nicefrac{\|u-u_h\|_{L^{\infty}}}{\|u\|_{L^{\infty}}}$},
				legend columns=2,
				legend pos=south east,
				legend style={at={(0.99,0.01)},anchor=south east},
				legend entries={$\Sp^2/\Sp^0$,$O(h^3)$,$\Sp^3/\Sp^1$,$O(h^4)$,$\Sp^4/\Sp^2$,$O(h^5)$,$\Sp^5/\Sp^{3}$,$O(h^6)$}]
				\addplot table [x=h, y=LinfP2, col sep=comma] {\DATAFILE};
				\addplot table [x=h, y=h3Linf, col sep=comma] {\DATAFILE};
				\addplot table [x=h, y=LinfP3, col sep=comma] {\DATAFILE};
				\addplot table [x=h, y=h4Linf, col sep=comma] {\DATAFILE};
				\addplot table [x=h, y=LinfP4, col sep=comma] {\DATAFILE};
				\addplot table [x=h, y=h5Linf, col sep=comma] {\DATAFILE};
				\addplot table [x=h, y=LinfP5, col sep=comma] {\DATAFILE};
				\addplot table [x=h, y =h6Linf, col sep=comma] {\DATAFILE};
			\end{loglogaxis}
		\end{tikzpicture}
		\label{fig:Quartercircle_LInf_C0_NME}
	\end{subfigure}
	\caption{Relative error for the quarter circle domain with $C^0$-constraints at the vertices and without knot vector reduction for the Lagrange multiplier spaces.}
	\label{fig:QuartercircleKNRC0}
\end{figure}
    \FloatBarrier
	\section{Conclusions}\label{sec:conclusions}
In this work, we propose an isogeometric method for solving the biharmonic equation on $C^0$-conforming multi-patch domains. The required $C^1$-continuity across patch interfaces is enforced weakly through a mortar formulation. We use a spline space of degree $p$ for the primal space, and a spline space of degree $p-2$ for the Lagrange multiplier space, which is built on the same underlying mesh as the primal space but with a reduced knot vector near the patch vertices. This approach ensures that the method satisfies a discrete inf-sup stability condition. As a result, we have established well-posedness for the method and proven optimal a priori error estimates in a patch-wise $H^2$ broken norm. Numerical examples confirm these error estimates and demonstrate optimal-order convergence in $L^2$, $H^1$, and $L^{\infty}$ norms.
{ An interesting direction for future work is the extension of our mortar approach to fourth-order problems posed on multi-patch surfaces. This extension is expected to be particularly useful when the surface has only weak smoothness across patch interfaces, although in that case new analytical tools and a suitable generalization of the present theory would be required.}

	\section*{Acknowledgements}
	G. Loli and G. Sangalli are members of the Gruppo Nazionale per il Calcolo Scientifico - Istituto Nazionale di Alta Matematica (GNCS-INDAM). G. Sangalli acknowledges the support of the Italian Ministry of University and Research (MUR) through the PRIN 2022 PNRR project NOTES (No. P2022NC97R), funded by the European Union - NextGenerationEU. This research also received financial support from ICSC - the Italian Research Center on High-Performance Computing, Big Data, and Quantum Computing, funded by the European Union - NextGenerationEU. {This research was funded in whole or in part by the Austrian Science Fund (FWF) 10.55776/P37177, a research project titled ``Isogeometric Multi-Patch Shells and Multigrid Solvers'' led by T. Takacs.} These contributions are gratefully acknowledged.
	%
	
	\bibliography{Bibliography}
	\bibliographystyle{ws-m3as}
	
\end{document}